\documentclass[11pt]{amsart}

\usepackage{amsgen, amsmath, amsfonts, amsthm, amssymb,
enumerate,amscd, color, epsf}
\usepackage[all]{xy}

\usepackage[varg]{txfonts}

\usepackage{multicol}

\newtheorem{thm}{Theorem}[section]
\newtheorem{prop}[thm]{Proposition}
\newtheorem{cor}[thm]{Corollary}
\newtheorem{lem}[thm]{Lemma}

\theoremstyle{definition}
\newtheorem{defn}[thm]{Definition}

\theoremstyle{remark}
\newtheorem{rmk}[thm]{Remark}

\newcommand\hol{^{\mathrm{hol}}}
\newcommand\Hol{{\mathcal{H}}}

\newcommand\h{{\sf h}}

\newcommand\al{\alpha}
\newcommand\bt{\beta}
\newcommand\Gm{\Gamma}
\newcommand\Gf{\Gamma} 
\newcommand\gm{\gamma}
\newcommand\Dt{\Delta}
\newcommand\dt{\delta}
\newcommand\e{\varepsilon}
\newcommand\z{\zeta}
\newcommand\Th{\Theta}
\renewcommand\th{\vartheta}
\renewcommand\k{\kappa}
\newcommand\Ld{\Lambda}
\newcommand\ld{\lambda}
\newcommand\x{\xi}

\newcommand\s{\sigma}
\newcommand\ups{\upsilon}

\newcommand\ph{\varphi}
\newcommand\ch{\chi}
\newcommand\ps{\psi}
\newcommand\om{\omega}

\newcommand\CC{\mathbb{C}}
\newcommand\NN{\mathbb{N}}
\newcommand\QQ{\mathbb{Q}}
\newcommand\RR{\mathbb{R}}
\newcommand\ZZ{\mathbb{Z}}

\newcommand\EE{{\mathbf{E}}}
\newcommand\HH{{\mathbf{H}}}
\newcommand\WW{{\mathbf{W}}}
\newcommand\XX{{\mathbf{X}}}
\newcommand\YY{{\mathbf{Y}}}
\newcommand\bb{{\mathbf b}}

\newcommand\ff{{\mathbf f}}
\renewcommand\gg{{\mathbf g}}
\newcommand\ii{{\mathbf i}}
\newcommand\jj{{\mathbf j}}
\newcommand\mm{{\mathbf m}}

\newcommand\C{{\mathcal C}}
\newcommand\D{{\mathcal D}}

\newcommand\M{{\mathcal M}}

\newcommand\V{{\mathcal V}}
\newcommand\W{{\mathcal W}}

\renewcommand\={\;=\;}
\newcommand\SL{{\mathrm{SL}}}
\newcommand\PSL{{\mathrm{PSL}}}
\newcommand\uhp{{\mathfrak H}}
\newcommand\npar{{n_{\mathrm{par}}}}
\newcommand\nell{{n_{\mathrm{ell}}}}
\newcommand\ngen{\ngeng{\Gm}}
\newcommand\ngeng[1]{{t(#1)}}
\newcommand\Gmod{{\Gm_{\!\mathrm{mod}}}}
\newcommand\tGmod{\tilde\Gm_{\!\mathrm{mod}}}
\newcommand\Gcom{{\Gm_{\!\mathrm{com}}}}
\newcommand\tGcom{\tilde\Gm_{\!\mathrm{com}}}
\newcommand\tG{{\tilde G}}
\newcommand\tGm{{\tilde\Gm}}
\newcommand\tgm{{\tilde\gm}}
\newcommand\tDt{{\tilde\Dt}}
\newcommand\tZ{{\tilde Z}}
\newcommand\oh{{\mathrm{O}}}
\newcommand\re{\mathrm{Re}\,}
\newcommand\im{\mathrm{Im}\,}
\newcommand\sign{\,\mathrm{Sign}}
\newcommand\supp{\mathrm{Supp}}
\newcommand\hypg[2]{{}_{#1}\!F_{\!#2}}
\newcommand\fd{{\mathfrak F}}

\newcommand\proj[1]{{\mathbb P}^1_{#1}}
\newcommand\map{\mathrm{Map}}
\newcommand\cp{\mathrm{cpt}}
\newcommand\divides{\mathrel{\bigm|}}

\newcommand\ml{{\sf m}}
\newcommand\pr{{\mathrm{pr}}}
\newcommand\glie{{\mathfrak g}}
\renewcommand\setminus{\smallsetminus}

\newcommand\Ft[1]{F_{\!#1}\,}
\newcommand\discr{{\mathrm{discr}}}
\newcommand\cont{{\mathrm{cont}}}
\newcommand\mult{\mathrm{Mult}}
\newcommand\A{{\mathcal E}}
\newcommand\coker{{\mathrm{coker}}}

\newcommand\txtfrac[2]{{\textstyle\frac{#1}{#2}}}

\newcommand\iv[2]{#1^{\Gm\!,#2}}
\newcommand\ivv[3]{#1^{#2,#3}}
\newcommand\ivt[2]{#1^{\tGm,#2}}

\newcommand\summ[2]{ \mathchoice
{\mathop{{\sum}^{#1}}_{#2\hspace*{.2em}}} {\sum^{#1}_{#2}}
{\sum^{#1}_{#2}} {\sum^{#1}_{#2}} }

\newcommand\rmrk[1]{\smallskip\par\noindent\emph{#1. }}

\makeatletter
\newcommand\hmatc[4]{\left[ {#1\@@atop #3}{#2\@@atop #4}\right]}
\newcommand\hmatr[4]{\left[ {\hfill #1\@@atop\hfill #3}{\hfill
#2\@@atop\hfill #4}\right]}
\newcommand\matc[4]{\left( {#1\@@atop #3}{#2\@@atop #4}\right)}
\newcommand\matr[4]{\left( {\hfill #1\@@atop\hfill #3}{\hfill
#2\@@atop\hfill #4}\right)}
\makeatother

\numberwithin{equation}{section}
\newcommand\equationstandard{\arabic{section}.\arabic{equation}}
\renewcommand\theequation{\equationstandard}

\newcounter{eqaux}
\newcommand\beqabc%
 { \stepcounter{equation}%
   \edef\tempeq{\theequation}%
   \setcounter{eqaux}{\arabic{equation}}%
   \renewcommand\theequation{\tempeq\alph{equation}}%
   \setcounter{equation}{0}%
 }

\newcommand\beqabclbl[1]%
 { \refstepcounter{equation}\label{#1}%
   \edef\tempeq{\theequation}%
   \setcounter{eqaux}{\arabic{equation}}%
   \renewcommand\theequation{\tempeq\alph{equation}}%
   \setcounter{equation}{0}%
 }
\makeatother

\begin{document}
\title{Higher order Maass forms}
\author{Roelof Bruggeman}
\address{Mathematisch Instituut Universiteit Utrecht, Postbus 80010,
NL-3508 TA Utrecht, Nederland}
\email{r.w.bruggeman@uu.nl}
\author{Nikolaos Diamantis}
\address{School of Mathematical Science, University of Nottingham,
Nottingham NG7 2RD, UK}
\email{Nikolaos.Diamantis@maths.nottingham.ac.uk}

\begin{abstract} The linear structure of the space of Maass forms of
even weight and of arbitrary order is determined.
\end{abstract}


\maketitle

\section{Introduction} In this work, the structure of the space of
Maass forms of general order and integral weight as a linear vector
space is determined. It is proved that, under suitable conditions,
this space is as large as one would expect it to be.

There are mainly two objects and associated problems that suggest the
study of specifically this type of higher-order form. The first is
Eisenstein series modified with modular symbols defined by
\begin{equation}
\label{eismod}
E^*(z, s)\=\sum_{\gm \in \Gm_{\infty} \backslash \Gm_0(N)} \langle f,
\gm \rangle \;\text{Im}(\gm z)^s\,,\end{equation}
where $\Gm_{\infty}$ is the subgroup of translations of the congruence
group $\Gm_0(N)$, $f$ a weight $2$ newform and $\langle f, \gm 
\rangle :=-2\pi i\int_{\infty}^{\gm \,  \infty} f(w) dw.$ The study
of this function has led to important results, such as the proof that
the suitably normalised modular symbols follow the normal
distribution (\cite{PR}). 
The function $E^*(-, s)$ is not automorphic but transforms as a
second-order automorphic form.

We recall that, for a group $\Gm$ of motions on the upper
half-plane~$\uhp$, a function is said to be $\Gm$-invariant of
\emph{order~$q\in \NN$} and weight $0$, if it satisfies
\begin{equation}
f|(\gm_1-1) (\gm_2-1) \cdots (\gm_q-1)\=0\quad \text{ for all
}\gm_1,\gm_2,\ldots,\gm_n \in \Gm\,. \label{hof}
\end{equation}
Here, the action $|$ of $\Gm$ on functions on $\uhp$ is given by
$$f|\gm (z)\;:=\; f(\gm z)\,.$$
and it is extended linearly to an action of the group ring $\CC[\Gm]$.

Clearly, several types of conditions on holomorphicity, growth etc.\
can be imposed on functions of general order. The function
$E^*(-, s)$ in particular, is an eigenfunction 
of the Laplacian and therefore we view it as a Maass form of order
$2$.

The second object leading to functions that are $\Gm$-invariant of
second-order arises from considerations related to values of
derivatives of $L$-functions of cusp forms: In \cite{G} and \cite{Di}
certain ``period integrals'' are associated to derivatives of
$L$-functions of weight $2$ cusp forms in a way analogous to the link
between values of $L$-functions and modular integrals
(\cite{Man}). Specifically, let $f$ be a newform of weight $2$ for
$\Gamma_0(N)$ and let $L_f(s)$ be its $L$-function. If $L_f(1)=0$,
then, for each prime $p$, $(p, N)=1$, $L'_f(1)$ can be written as a
linear combination of integrals of the form
\vskip -2mm
\begin{equation}\label{fu}
\int_0^{\gamma(0)} f(z)\, u(z)\, dz\,, \qquad \gamma \in \Gamma_0(N)
\end{equation}
plus some ``lower order terms". Here $u(z):=
\log\eta(z)+\log\eta(Nz)$, where $\eta$ is the Dedekind
$\eta$-function. The differential $f(z) u(z) \, dz$ is not
$\Gamma_0(N)$-invariant. It does satisfy a transformation law which
is reminiscent of (\ref{hof}), but is not quite $\Gm_0(N)$-invariant
of order $2$ in the narrow sense. If it were, the value of the
derivative at $1$ would be expressed as the value of the actual
$L$-function of second-order $\Gm_0(N)$ at $1$. That could be
advantageous for the study of $L_f'(1)$ in terms of the outstanding
conjectures, especially since there is now evidence that a motivic
structure underlies higher order forms (see \cite{DSr} and
\cite{Sr}).

Here we show that it is indeed possible to obtain a second-order
$\Gm_0(N)$-invariant function from $u(z)$ provided we move to a
different domain. This domain is the universal covering group which
we will be defining in detail in \S\ref{Ucg}.

As will become apparent in the sequel, it is natural, in higher
orders, to unify the study of Maass forms and that of forms on
universal covering groups. The full definition of the
\emph{higher-order Maass forms with generalised weight on the
universal covering group} is discussed in \S\ref{sect-gw-uc}. Theorem
\ref{thm-tAi} then allows us to translate results on the universal
covering group to the analogous results on the upper-half plane.

A fundamental question is how ``large'' this space is. In the case of
\it holomorphic \rm higher-order cusp forms, the corresponding spaces
are finite-dimensional and the answer can be given by computing the
dimensions (\cite{DO} and \cite{DSi}). In the present case, where the
relevant space is not finite dimensional, a different
characterisation of ``size'' is required. Such a characterisation is
proposed in \S\ref{sect-hoi}.

Although our results imply that there are ``many" higher order Maass
forms, the proofs are highly inductive and do not easily lead to
explicit examples. In \S\ref{sect-ex1} and \S\ref{sect-ex2} we
address this problem, by illustrating various methods that lead to
explicit examples of higher order Maass. Surprisingly, these examples
are derived very naturally from the theory which was developed in a
completely different context in \cite{Br86, Br94}.

Finally, a particular aspect of the proof that deserves to be singled
out because of its independent interest is the definition of
genuinely higher-order Fourier expansions. Higher order automorphic
forms need not be invariant under the group fixing a cusp, so there
is no obvious Fourier expansion. To date, to address this problem one
had to partially revert to the classical setting by imposing the
somewhat unnatural extra condition of invariance under the parabolic
elements of the group. In \S\ref{sect-Ft}, appropriate higher-order
Fourier terms are constructed, thus avoiding additional invariance
conditions.

\vskip -7mm

\section{Structure of the paper} \label{sect-struct}
In \S\ref{sect-hoi} we first discuss higher-order invariants for
general groups and modules. This allows a precise definition of the
concept of ``as large as possible" (\it maximally perturbable\rm). A
first maximal perturbability result for a general space of maps is
also proved.

In \S\ref{sect-mf}, Maass forms on $\uhp$ (both general and
holomorphic) are defined and the first two main theorems of the paper
(\ref{thm-Mf-mp} and~\ref{thm-hol-mp})
are stated. The section includes an extended discussion of concrete
examples of low-order forms on $\uhp$.

In the next section the universal covering group $\tG$ is introduced
and the basic facts about $\tG$ are given.

Maass forms on the universal covering group are defined in
\S\ref{sect-gw-uc} and the counterparts of Theorems \ref{thm-Mf-mp}
and~\ref{thm-hol-mp} for forms on the universal covering group are
stated. The section concludes with concrete examples of low-order
forms on $\tG$.

Section ~\ref{sect-Ft} is of independent interest. A theory of Fourier
expansions for higher-order forms is developed.

The proof of Theorems \ref{thm-Mf-mp} and~\ref{thm-hol-mp} is the
content of \S\ref{sect-hoM}. The proof involves the construction of
two spaces with support conditions. To deduce their maximal
perturbability we employ spectral techniques.

\section{Higher order invariants} \label{sect-hoi}
In this section, we discuss higher order invariants in general and
then specialise their study to discrete cofinite subgroups
$\Gm \subset \PSL_2(\RR)$. We introduce the concept of a ``maximally
perturbable'' $\Gm$-module to make precise the statement that there
are as many higher order invariants of a given type as one can
expect. A first maximal perturbability result in a general context is
proved. 

\subsection{Higher order invariants on general groups} The concept of
higher order invariant functions on the upper half plane is a special
case of the concept ``higher order invariants'' for any group~$\Gm$
and any $\Gm$-module~$V$. We work with \emph{right} $\Gm$-modules, an
write the action as $v\mapsto v|\gm$. It should be clear from the
context when we refer to this general meaning of $|$ and when to the
more narrow meaning given in the Introduction. We define the
\emph{higher order invariants} inductively:
\begin{equation}\label{invdef}
\begin{aligned}
\iv V 1 &\= V^\Gm \= \bigl\{ v\in V\;:\; \forall \gm\in \Gm, \,
v|\gm=v\bigr\}\,, \\
\iv V {q+1} &\= \bigl\{ v\in V\;:\; \forall \gm\in \Gm, \, v|(\gm-1)
\in \iv V q\bigr\}\,.
\end{aligned}
\end{equation}
We set $\iv V 0=\{0\}$.

Let now $\Gm$ be finitely generated and let $I$ be the augmentation
ideal in the group ring $\CC[\Gm]$, generated by $\gm-1$ with
$\gm\in \Gm$. A fundamental role in the paper will be played by the
map
$$\ml_q: V^{\Gm, 
q+1} \to \text{hom}_{\mathbf C [\Gm]}(I^{q+1}\backslash 
I^q, V^{\Gm}).$$
To define it we first quote from \cite{Dei} (before Proposition~1.2):
\begin{equation}\label{iv-I}
\iv V q \;\cong\; \hom_{\CC[\Gm]}(I^q\backslash \CC[\Gm],V)\,.
\end{equation}
Next, we note that $I^{q+1}\backslash I^q$ is generated by
\[ I^{q+1}+(\gm_1-1)\cdots(\gm_q-1)\,,\]
with $\gm_i \in \Gm$. To each $v \in V^{\Gm, q+1}$ we associate the
map on $I^{q+1}\backslash I^q$ sending this element to
$v|(\gm_1-\nobreak 1) \cdots (\gm_q-\nobreak 1)$.
This map is well-defined because
$v|(\gm_1-\nobreak 1) \cdots (\gm_{q+1}-\nobreak 1)=0$.
In this way, we obtain a map $\ml_q$ from $V^{\Gm, q+1}$ to
$$
\hom_{\CC[\Gm]}(I^{q+1}\backslash I^q, V) \cong
\hom_{\CC[\Gm]}(I^{q+1}\backslash I^q, V^{\Gm})$$
(since the action induced on $I^{q+1}\backslash I^q$ by the operation
of $\Gm$ is trivial). It is easy to see that the kernel of $\ml_q$ is
$V^{\Gm, q}$ and thus we obtain the exact sequence
\begin{equation}\label{exs-iv}
0 \longrightarrow \iv V q \longrightarrow \iv V {q+1} \stackrel{\ml_q}
\longrightarrow \hom_{\CC[\Gm]}(I^{q+1}\backslash I^q,V)
\end{equation}
The map $\ml_q$ may or may not be surjective and we will interpret the
phrase ``as large as possible'' as surjectivity of~$\ml_q$ for all
$q\in\NN$.

\begin{defn} Let $\Gm$ be a finitely generated group. We will call a
$\Gm$-module $V$ \emph{maximally perturbable} if the linear map
$\ml_q: \iv V{q+1} \rightarrow 
\hom_{\CC[\Gm]}(I^{q+1}\backslash I^q,V^\Gm)$ is surjective for all
$q\geq 1$.\end{defn}

A reformulation of this definition which is occasionally easier to
use, uses the finite dimension
\begin{equation}
\label{dimen}
n(\Gm,q):=\dim_{\CC}(I^{q+1}\backslash I^q).
\end{equation}
$V$ is maximally perturbable if and only if
$\iv V{q+1}/\iv V q \cong (V^\Gm)^{n(\Gm,q)}$ for all $q\in \NN$.

In \cite{DSi} higher order cusps forms of weight~$k$ for a
discrete group $\Gm$ are considered in the space of holomorphic
functions on $\uhp$ with exponential decay at the cusps that moreover
are invariant under the parabolic transformations. The dimensions of
these spaces are computed and generally turn out to be strictly
smaller than $n(\Gm,q)$. So the corresponding $\Gm$-module is not
maximally perturbable. 
\medskip

A useful definition is based on the isomorphism
$\hom_{\CC[\Gm]}(I^{q+1}\backslash 
I^q, V^\Gm) \cong \mult^q(\Gm,V^\Gm)$, the space of maps $\Gm^q 
\rightarrow V^\Gm$ inducing group homomorphisms $\Gm\rightarrow \CC$
on each of their coordinates. For a finitely generated group $\Gm$,
$\mult^q(\Gm,V^\Gm) \cong \mult^q(\Gm,\CC) \otimes_\CC V^\Gm$ where
$\mult^q(\Gm,\CC)$ is the $q$-th tensor power of the abelianised
group $\Gm^{\mathrm{ab}}=\Gm / [\Gm,\Gm]$. With this notation we
define
\begin{defn}\label{pert-def}Let $q\in \NN$. For any group $\Gm$ and
any $\Gm$-module~$V$ we call $f\in \iv V q$ a \emph{perturbation} of
$\ph\in V^\Gm$ if there exists $\mu_f \in \mult^q(\Gm,\CC)$ such that
for all $\gm_1,\cdots,\gm_q\in \Gm$:
\begin{equation}
\label{pert} f|(\gm_1-1)\cdots(\gm_q-1) \= \mu_f(\gm_1,\ldots,\gm_q)
\,\ph\,.
\end{equation}

We call a perturbation \emph{commutative} if $\mu_f$ is invariant
under all permutations of its arguments. If not, we call it
non-commutative.
\end{defn}

\subsection{Higher order invariants on subgroups of $\PSL_2(\RR)$}
\subsubsection{Canonical generators} \label{canon-gener}
In this section we recall the ``canonical generators'' of cofinite
discrete subgroups of~$\PSL_2(\RR)$, and use them to show that
certain modules are maximally perturbable.
\medskip

Let $\Gm\subset \PSL_2(\RR)$ be a cofinite discrete group of motions
in the upper half-plane~$\uhp$. A system of \emph{canonical
generators} for $\Gm$ consists of
\begin{itemize}
\item Parabolic generators $P_1,\ldots,P_{\npar}$, each conjugate in
$\PSL_2(\RR)$ to $\pm\matc 1101$. We shall assume that $\Gm$ has
cusps: $\npar\geq 1$.
\item Elliptic generators $E_1,\ldots,E_\nell$, with $\nell\geq 0$.
Each $E_j$ is conjugate to
$\pm\matr{\cos(\pi/v_j)}{\sin(\pi/v_j)}{-\sin(\pi/v_j)}{\cos(\pi/v_j)}$
in $\PSL_2(\RR)$ for some $v_j\geq 2$.
\item Hyperbolic generators $H_1,\ldots,H_{2g}$, with $g\geq 0$, each
conjugate in $\PSL_2(\RR)$ to the image $\pm \matc t00{t^{-1}}$,
$t>1$, $t\neq 1$, of a diagonal matrix
\end{itemize}
See, {\sl e.g.}, \cite{Leh}, Chap.~VII.4, p.~241, or \cite{Pe48}, \S3.
The relations are given by the condition that each $E_j^{v_j}=I$ for
$j=1,\ldots, \nell$, and one large relation
\begin{equation}
P_1\cdots P_\npar E_1 \cdots E_\nell
\,[H_1,H_2]\cdots[H_{2g-1}H_{2g}]\=\mathrm{Id}\,.
\end{equation}
The choice of canonical generators is not unique, but the numbers
$\npar$, $\nell$ and $g$, and the elliptic orders
$v_1,\ldots,v_{\nell}$ are uniquely determined by~$\Gm$.

Each group homomorphism $\Gm\rightarrow\CC$ vanishes on the $E_j$, and
is determined by its values on $H_1,\ldots,H_{2g}$,
$P_1,\ldots,P_{\npar-1}$, hence
\begin{equation}\dim\hom(\Gm,\CC)\=\npar-1+2g\,.\end{equation}

We put $\ngen=\npar+2g$, and denote $A_1=P_1,\ldots,\allowbreak
A_{\npar-1}=P_{\npar-1},\allowbreak  A_{\npar}=H_1,\allowbreak
A_{\npar+1}=H_2,\ldots,\allowbreak  A_{\ngen-1}=H_{2g}$. The group
$\Gm$ is generated by $E_1,\ldots,E_{\nell}$ and
$A_1,\ldots, A_{\ngen-1}$.

For the modular group we have $\npar=1$, $P_1=\pm\matc 1101$,
$\nell=2$, $E_1=\pm\matr11{-1}0 $, $E_2=\pm S:=\pm\matr0{-1}10$,
$g=0$, and hence $\hom(\Gmod,\CC)=\{0\}$ and $\ngeng\Gmod=1$.
\medskip

In the sequel, we will need a basis for $I^{q+1}\backslash I^q$.
Arguing as in Lemma~2.1 in~\cite{Dei} we can deduce that the elements
\begin{equation}\label{bi-def0}
\bb(\ii) = (A_{\ii(1)}-1)\cdots(A_{\ii(q)}-1)\,,
\end{equation}
where $\ii$ runs over all $(\ngen\!-\!1)$-tuples of elements of $\{
1,\ldots,\ngen\!-\!1\}$ form a basis of $I^{q+1}\backslash I^q$. We
do not give a proof here, since it follows from the more general
result Proposition~\ref{prop-augpow}.

\subsubsection{A first maximal perturbability result} We view the
space $\map(\Gm,\CC)$ of all maps from $\Gm$ to~$\CC$ as a right
$\Gm$-module for the action $|$ by left translation.
\begin{prop}\label{prop-Gmfct-mp}If $\Gm$ is a discrete cofinite
subgroup of~$\PSL_2(\RR)$ with cusps, then $\map(\Gm,\CC)$ is
maximally perturbable.
\end{prop}
\begin{proof}We construct functions $\gg_\ii\in \map(\Gm,\CC)$ for
$n$-tuples $\ii$ from $\{1,\ldots,\ngen-\nobreak 1\}$. Firstly, let
$\Gm_0$ be the free subgroup of $\Gm$ which is generated by the
elements $A_j$, $1\leq j \leq \ngen-1$. It is clear that there is 
a unique system of functions $\{\gg_\ii\}$ on $\Gm_0$ such that
\begin{equation}\label{gidef}
\begin{aligned}
\gg_{()}&\=1\,,\\
\gg_{(j,\ii)}|(A_j-1)&\= \gg_\ii\,,\\
\gg_\ii|(A_j-1)&\=0\quad\text{ if }\ii(1)\neq j,\\
\gg_\ii(1)&\=0\quad\text{ if }|\ii|\geq 1\,.
\end{aligned}
\end{equation}
By $|\ii|$ we denote the length of the tuple~$\ii$.

We next set $\gg_\ii(\gm) = \gg_\ii\bigl(\ph_0(\gm) \bigr)$ for
$\gm\in \Gm$, where $\ph_0$ is the homomorphism defined by
$\ph_0(E_j)=I$, $\ph_0(A_j)=A_j$. With the map $\ml_q$
in~\eqref{exs-iv} and for for $|\ii|=|\jj|$ we have on the basis
elements in~\eqref{bi-def0},
\begin{equation}
\begin{aligned}(\ml_q \gg_\ii )&\bigl(\bb(\jj) \bigr)
\= \gg_\ii|(A_{\jj(1)}-1) \cdots (A_{\jj(q)}-1)\\
&\= \dt_{\ii(1),\jj(1)} \gg_{\ii'} | (A_{\jj(2)}-1) \cdots
(A_{\jj(q)}-1)\,,
\end{aligned}
\end{equation}
where $\ii'$ is the tuple $(\ii(2),\ldots,\ii(q))$. Inductively we
obtain
\begin{equation}(\ml_q \gg_\ii )\bigl(\bb(\jj) \bigr)\=
\dt_{\ii,\jj}:= \prod_{l=1}^q \dt_{\ii(l), \jj(l)}\,.
\end{equation}
Hence the $\gg_\ii$ with $|\ii|=q$ form a dual system for the
generators $\bb(\ii)$. This implies that the image $\ml_q \iv V{q+1}$
has maximal dimension~$n(\Gm,q)$.
\end{proof}

\section{Maass forms} \label{sect-mf}
We turn to spaces of functions on the upper half-plane that
contain the classical holomorphic automorphic forms and the more
general Maass forms. The first main results of this paper are stated
in Theorems~\ref{thm-Mf-mp} and~\ref{thm-hol-mp}. In \S\ref{sect-ex1}
we give some explicit examples of higher order Maass forms.

\subsection{General Maass forms} Let $\Gm$ be a cofinite discrete
subgroup~$\Gm$ of the group $G=\PSL_2(\RR)$. For each cusp $\k$, we
choose $g_\k \in\PSL_2(\RR)$ such that
\begin{equation}
\label{scal}
\k=g_\k\infty \qquad \text{and} \, \, g_\k^{-1} \Gm_{\k} g_\k\=\Bigl\{
\pm\matc 1n01\;:\; n\in \ZZ\Bigr\}
\end{equation}
Here, $\Gm_{\k}$ is the set of elements of $\Gm$ fixing $\k$.
The elements $g_\k$ are determined up to right multiplication
by elements $\pm \matc ab0{a^{-1}}\in G$. We choose the $g_\k$ for
cusps in the same $\Gm$-orbit so that
$g_{\gm\k} \in \gm g_\k \Gm_\infty$.

We further consider a generalisation of the action $|$ considered in
the last section. For a fixed $k$ and for a $f: \uhp \to \CC$ we set
\begin{equation}\label{action}
f \Bigm|_k \matc abcd (z) \= (cz+d)^{-k}\, f\bigl(
(az+b)/(cz+d)\bigr).
\end{equation}We finally set
\begin{equation}
\label{Lkdef}L_k \= -y^2\partial_x^2-y^2\partial_y^2 +ik y
\partial_x-ky
\partial_y + \frac k2\bigl(1-\frac k2 \bigr) \,.\end{equation}

With this notation we have
\begin{defn}\label{MAkdef} Let $k\in 2\ZZ$ and $\ld\in \CC$.
\begin{enumerate}
\item[i)] $\M_k(\Gm, \ld)$ denotes the space of smooth functions
$f:\uhp \rightarrow \CC$ such that $L_kf=\ld f$ and for which there
is some $a\in \RR$ such that
\begin{equation}
\begin{aligned}
 f(g_\k&(x+iy))\=\oh(y^a)\quad(y\rightarrow \infty)
 \\
&\text{ uniformly for $x$ in compact sets in $\RR$, for \it all \rm
cusps $\k$ of $\Gm$ }\,.
\end{aligned}
\end{equation}
\item[ii)] $\A_k(\Gm, \ld)$ denotes the space of smooth functions $f$
such that $L_kf=\ld f$ and for which there is some $a\in \RR$ such
that
\begin{equation}
\begin{aligned}
 f(g_\k&(x+iy))\=\oh(e^{ay})\quad(y\rightarrow \infty)\\
&\text{ uniformly for $x$ in compact sets in $\RR$, for \it all \rm
cusps $\k$ of $\Gm$.}
\end{aligned}\end{equation}
\item[iii)] We denote the invariants in these spaces by
\begin{equation}
\label{Maass}
E_k(\Gm,\ld)\,:=\,\A_k(\Gm,\ld)^\Gm \quad \text{and} \;
M_k(\Gm,\ld)\,:=\,\M_k(\Gm,\ld)^\Gm\,.
\end{equation}
We call the elements of $E_k(\Gm,\ld)$ (resp. $M_k(\Gm,\ld)$)
\it Maass forms of polynomial (resp. exponential)
growth of weight $k$ and eigenvalue $\ld\in \CC$ for~$\Gm$. \rm
\end{enumerate}
\end{defn}

\noindent\it
Remarks. \rm
\begin{enumerate}
\item[i)] Since $L_k$ is elliptic, all its eigenfunctions are
automatically real-analytic.
(See, {\sl e.g.}, \cite{La74}, \S5 of App.~A4, and the references
therein.)
If $f$ is holomorphic, then it is an eigenfunction of~$L_k$ with
eigenvalue $\frac k2\bigl(1-\frac k2  \bigr)$.
\item[ii)] The space $M_k(\Gm,\ld)$ is known to have finite dimension.
The space $E_k(\Gm,\ld)$ has, for groups $\Gm$ with cusps, infinite
dimension. The subspace of $E_k(\Gm,\ld)$ corresponding to a fixed
value of~$a$ in the bound $\oh(e^{ay})$ has finite
dimension.
\item[iii)] In an alternative definition, suitable for functions not
necessarily holomorphic, one replaces the Maass forms $f$ as defined
above by $h(z) = y^{k/2} f(z)$. Then invariance under \eqref{action}
becomes invariance under the action
\begin{equation}
\label{actionb}
f \Bigm|_k \matc abcd (z) \= e^{-ik\arg(cz+d)}\, f\bigl(
(az+b)/(cz+d)\bigr)\,
\end{equation}
and the eigenproperty in the terms of the Laplacian
\begin{equation}\label{efC}
\bigl( - y^2\partial_x^2-y^2\partial_y^2+ i k y \partial_x \bigr) h \=
\ld h. \end{equation}
The formulation of the growth conditions remains unchanged. Now
antiholomorphic automorphic forms $a(z)$ of weight~$k$ give Maass
forms $h(z) = y^{k/2} a(z)$ of weight~$-k$.
\end{enumerate}

Our main result for general Maass forms on $\uhp$ is
\begin{thm}\label{thm-Mf-mp}Let $\Gm$ be a cofinite discrete group of
motions in $\uhp$ with cusps. Then the $\Gm$-module $\A_k(\Gm,\ld)$
is maximally perturbable for each $k\in 2\ZZ$ and each $\ld \in \CC$.
\end{thm}
In the course of the proof in \S\ref{sect-hoM} we will see that even
if we start with Maass forms with polynomial growth the construction
of higher order invariants will lead us to functions that have
exponential growth.

\subsection{Holomorphic automorphic forms} \label{haf} For even $k$
the space $\A_k(\Gm,\ld_k)$, with
$\ld_k=\frac k2(1-\nobreak\frac k2)$ contains the subspace
$\A_k\hol(\Gm,\ld_k)$ where the condition $L_kf=\lambda_k f$ is
replaced by the stronger condition that $f$ is holomorphic. In the
alternative definition, condition~\eqref{efC} is replaced by the
condition that $z\mapsto y^{-k}f(z)$ is holomorphic. The space
$\A_k\hol(\Gm,\ld_k)$ is a $\Gm$-submodule of $\A_k(\Gm,\ld_k)$. We
also have the $\Gm$-submodule
$\M_k\hol(\Gm,\ld_k)=\M_k(\Gm,\ld_k) \cap \A_k\hol(\Gm,\ld_k)$ of
$\M_k(\Gm,\ld_k)$.

The space $\M_k\hol(\Gm,\ld_k)^\Gm$ is the usual space of entire
weight $k$ automorphic forms for~$\Gm$, and $\A_k\hol(\Gm,\ld_k)^\Gm$
is the space of meromorphic automorphic forms with singularities only
at cusps. Sometimes, e.g. in~\cite{BOR}, the elements of
$\A_k\hol(\Gm,\ld_k)^\Gm$ are called \emph{weakly holomorphic}. There
the elements of $\A_k(\Gm,\ld_k)^\Gm$ are called \emph{harmonic weak
Maass forms}. We prefer to use the term \emph{harmonic} for Maass
forms in $\A_k(\Gm,0)^\Gm$. (Note that $\ld_k\neq 0$ for $k \ne 
0, 2$.)

Our main result for holomorphic automorphic forms on $\uhp$ is:
\begin{thm}\label{thm-hol-mp}Let $\Gm$ be a cofinite discrete group of
motions in~$\uhp$ with cusps. Then $\A_k\hol(\Gm,k/2-\nobreak k^2/4)$
is maximally perturbable for each $k\in 2\ZZ$.
\end{thm}

\subsection{Examples of harmonic and holomorphic forms of order two
and three}\label{sect-ex1}
According to Theorems \ref{thm-Mf-mp} and~\ref{thm-hol-mp} there are
plenty of examples of higher order Maass forms for cofinite groups
with cusps for which $\dim_\CC \hom(\Gm,\CC)\geq 1$. It is, however,
not very easy to exhibit explicit examples.

For the modular group $\Gmod=\PSL_2(\ZZ)$ the space $\hom(\Gmod,\CC)$
is zero. Hence it does not accept higher order invariants. For the
commutator subgroup $\Gcom=[\Gmod,\Gmod]$ we will employ three
different approaches to exhibit full sets of perturbations of~$1$ (as
defined in~Definition~\ref{pert-def}) of orders two and three. A
reader only interested in the existence of higher order forms may
prefer to skip this subsection.

\subsubsection{Holomorphic perturbation of~$1$}\label{sect-hol-per1}
In \cite{Leh}, Chap.~XI, \S3E, p.~362, one finds various facts
concerning~$\Gcom$. It is freely generated by $D=\pm \matc 2111$ and
$C=\pm \matr2{-1}{-1}1$. It has no elliptic elements, and one
cuspidal orbit $\Gcom\;\infty =\proj\QQ$. The group $(\Gcom)_\infty$
fixing $\infty$ is generated by $\pm \matc1601$. We have
$\ngeng{\Gcom}=3$.\smallskip

The space of holomorphic cusp forms of weight~$2$ has dimension $g=1$.
We use the basis element $\eta^4$ (power of the Dedekind
eta-function). The map
\begin{equation}\label{Hdef}
H(z) \= -2\pi i\int_\infty^z \eta(\tau)^4\, d\tau \= -6 e^{\pi i z/3}
+ \oh\bigl( e^{7\pi i z/3}\bigr)
\end{equation}
induces an embedding of $\Gcom\backslash\uhp$ into an elliptic curve,
which can be described as $\CC/\Ld$, with
\begin{equation}
\Ld\= \varpi\, \ZZ[\rho]\,,\qquad
\varpi\=\pi^{1/2}\Gf(1/6)/\bigl(6\sqrt 3\,\Gf(2/3)\bigr)\,,\quad
\rho\=e^{\pi i/3}\,.
\end{equation}
(See computations in \S15.2--3 in~\cite{Br94}.) The map $H$ maps
$\uhp$ onto $\CC\setminus\Ld$, and satisfies for $\gm\in \Gcom$
\begin{equation}\label{H-ld-def}
H(\gm z)\= H(z) + \ld(\gm)\,,\qquad \ld(\gm) \= -2\pi i
\int_\infty^{\gm\infty}\eta(\tau)^4\, d\tau\,,
\end{equation}
where $\ld(C) = \rho\varpi$ and $\ld(D) = \bar\rho\varpi$. So the
lattice $\Ld$ is the image of $\ld : \Gcom\rightarrow\CC$, and
$\hom(\Gcom,\CC) = \mult^1(\Gcom,\CC)$ has $\ld$, $\bar\ld$ as a
basis. We note that the kernel $\ker(\ld)$ is a subgroup with
infinite index in $\Gcom$; it is in fact the commutator subgroup
of~$\Gcom$. The element $\pm\matc 1601$ generating the subgroup
of~$\Gcom$ fixing $\infty$ is in $\ker(\ld)$. Since $\ker(\ld)$ has
no elliptic elements, composition with $H$ gives a bijection from the
holomorphic functions on $\CC\setminus \Ld$ to the holomorphic
$\ker(\ld)$-invariant functions on~$\uhp$.

Clearly, $H$ is a holomorphic second order perturbation of~$1$ with
linear form~$\ld$. It is also a \emph{harmonic perturbation} of~$1$,
i.e., a perturbation which is harmonic as a function. By conjugation
we obtain the antiholomorphic harmonic perturbation of~$1$ with
linear form $\bar\ld$.

According to Theorem~\ref{thm-hol-mp} there should also be a
holomorphic second order perturbation of~$1$ with a linear form
that is linearly independent of~$\ld$. Here we can use the
Weierstrass zeta-function
\begin{equation}
\z(u;\Ld)\= \frac{1}{u}+\summ'{\om \in \Ld}\biggl(\frac1{u-\om}+\frac
1 \om
+\frac{u}{\om^2}\biggr)\,.
\end{equation}
See, {\sl e.g.}, \cite{KK}, Chap.~I, \S6. It is holomorphic on
$\CC\setminus\Ld$ and satisfies $\z(u+\nobreak\om;\Ld) = \z(u;\Ld) +
\h(\om)$ for all $\om\in \Ld$, where $\h \in \hom(\Ld,\CC)$ is
linearly independent of~$\om\mapsto \om$. (The classical notation for
$\h$ is $\eta$. We write $\h$ to avoid confusion with the Dedekind
eta function.)
Pulling back this zeta-function to~$\uhp$ we get a second order
holomorphic perturbation of~$1$
\begin{equation}
 W(z)\= \z\bigl(H(z);\Ld\bigr)\end{equation}
with the linear form $\gm\mapsto\h\bigl(\ld(\gm)\bigr)$. The Laurent
expansion of the Weierstrass zeta-function at $0$ starts with
$\z(u;\Ld) = u^{-1} + \oh(u^3)$. Hence $W$ has a Fourier expansion
at~$\infty$ starting with
\begin{equation}
W(z) \= \frac{-1}6\, e^{-\pi i z/3} + \oh\bigl( e^{\pi i z} \bigr)\,.
\end{equation}
This shows that $W$ has exponential growth at the cusps.

We may carry this out also for holomorphic forms of order three, to
obtain the following commutative perturbations of~$1$ of order~$3$:
\begin{equation}
\renewcommand\arraystretch{1.2}
\begin{array}{|c|c|c|c|}
\hline
f& H(z)^2 &H(z)\,W(z) &W(z)^2\\
\mu_f& 2\,\ld \otimes \ld&
\ld\otimes(\h\circ\ld)+(\h\circ\ld)\otimes\ld
& 2(\h\circ\ld)\otimes(\h\circ\ld)\\
\hline
\end{array}
\end{equation}

We know that there also exist non-commutative holomorphic
perturbations of order~$3$. To find an explicit example, we have to
work on~$\uhp$, since the group $\Ld$ acting on~$\CC$ is
abelian.\smallskip

The closed holomorphic $1$-forms
\[ \om=-2\pi i\, \eta(\tau)^4\,d\tau\quad\text{ and }\quad \om_1 =
-2\pi i \, W(\tau)\, \eta(\tau)^4\, d\tau\] on $\uhp$ transform as
follows under~$\Gcom$:
\begin{equation}
\om|\gm\=\om\,,\qquad \om_1|\gm \= \om_1 + \h\bigl(\ld(\gm) \bigr)\,
\om\,.
\end{equation}
For an arbitrary base point $z_0\in \uhp$ we put
\begin{equation}\label{iiK}
K(z) \= \int_{z_0}^z \om_1\,.
\end{equation}
This defines a holomorphic function on~$\uhp$ that satisfies for
$\gm\in \Gcom$:
\[ K|(\gm-1)(z) \= \int_z^{\gm z} \om_1\,,\]
and hence for $\gm,\dt\in \Gcom$:
\begin{align*} K|(\gm-1)(\dt-1)(z) &\= \biggl( \int_{\gm z}^{\gm\dt z}
- \int_z^{\dt z} \biggr)\, \om_1 \= \int_{z}^{\dt z} \om_1|\gm
-\int_z^{\dt z}\,\om_1 \\
&\= \h\bigl( \ld(\gm) \bigr)\,\int_z^{\dt z} \om \= \h\bigl(
\ld(\gm)\bigr)\, \ld(\dt)\,.
\end{align*}
Thus, we have a holomorphic third order non-commutative perturbation
$K$ of~$1$ with non-symmetric multilinear form
$(\h\circ\ld) \otimes\ld$. Since holomorphic forms are harmonic in
weight zero these perturbations are also harmonic perturbations
of~$1$.

\subsubsection{Iterated integrals} The construction of the third order
form $K$ in \eqref{iiK} is closely related to the iterated integrals
used in \cite{DSr} to prove maximal perturbability of spaces of
smooth functions.

The idea is that we have two closed $\Gcom$-invariant differential
forms on~$\uhp$, $dH(z)=\om=-2\pi i \,\eta(z)^4\, dz$, and
\[ \om_0 \= dW(z) \= -\wp\bigl( H(z) \bigr)\, d \bigl( H(z)
\bigr)\,,\]
where $\wp(u;\Ld) = - \frac d{du}\z(u;\Ld)$ is the Weierstrass
$\wp$-function. If $t\mapsto z(t)$, $0\leq t\leq 1$ is a path
in~$\uhp$ from $z_0$ to $z_1$, then
\begin{align*}
\int_{t_2=0}^1 \int_{t_1=0}^{t_2} &\om_0\bigl( z(t_1) \bigr) \,
\om\bigl( z(t_2)
\bigr)
\= \int_{t_2=0}^1 \bigl( W\bigl( z(t_2) \bigr) -W(z_0)\bigr) \,
dH\bigl(z(t_2) \bigr)\\
&\=-2\pi i \int_{t=0}^1 W\bigl( z(t)\bigr) \, \eta\bigl( z(t)
\bigr)^4\, z'(t)\, dt - W(z_0) \bigl(H(z_1)-H(z_0)
\bigr) \\
&\= K(z_1)- W(z_0) \bigl(H(z_1)-H(z_0)
\bigr)
\end{align*}
depends only on $z_0$ and $z_1$, not on the actual path. For a fixed
base point $z_0$ the holomorphic function $z_1\mapsto W(z_0)\bigl(
H(z_1)-H(z_0) \bigr)$ is invariant of order two. So up to lower order
terms the invariant $K$ is given by an iterated integral, as in (3)
of~\cite{DSr}; see also~\cite{Ch}.

\subsubsection{Differentiation of families} We start by considering a
general finitely generated group $\Gm $ acting on a space~$X$. We
will use the notation $f|\gm(x) = f(\gm x)$ for the action induced on
functions defined on $X$. We consider a family of characters of~$\Gm$
of the form $\ch_r(\gm) = e^{i r \cdot \al(\gm)}$, where
$r\cdot\al(\gm)= r_1\al_1(\gm)+\cdots+r_n\al_n(\gm)$,
$\al_1,\ldots,\al_n\in \hom(\Gm,\RR)$, $r$ varying over an open set
$U$ in~$\RR^n$. In this way $\ch_r$ is a family of unitary
characters.

We consider a $C^\infty$ family $r\mapsto f_r$ on a
neighborhood~$U\subset\RR^n$ of~$0$ of functions $X\rightarrow\CC$
that satisfy
\begin{equation}\label{tach}
f_r(\gm x) \= \ch_r(\gm)\, f_r(x)\qquad(\gm\in \Gm)\,.
\end{equation}
We assume that $\ch_0$ is the trivial character and that $f_0$ is a
$\Gm$-invariant function $f$.

We now set $h(x) = \partial_{r_j} f_r(x) \bigr|_{r=0}$, for one of the
coordinates of~$r$. The transformation behaviour gives
$h(\gm x) \= i  \al_j(\gm)\,f(x) + h(x)$, or, rewritten,
\[ h| \gm - h \= i \al_j(\gm)\, f\,.\]
The function $h$ is a second order perturbation of $f$, with $i\al_j$
as the corresponding element of $\hom(\Gm,\CC)$. This can be
generalised:
\begin{prop}\label{prop-diff}For all multi-indices $a\in \NN^n$ the
derivative
\[f^{(a)}(x) :=
\partial_r^a f_r(x) \bigr|_{r=0}\]
is a commutative perturbation of $f$ with order $1+|a|$.
\end{prop}
We use the notations $\partial^a_r =
\partial_{r_1}^{a_1}\cdots\partial_{r_n}^{a_n}$ and
$|a|=a_1+a_2+\cdots+a_n$.
\begin{proof}
We use induction on the length~$|a|$ of the multi-index. The case
$|a|=1$ has already been handled above. For $|a|>1$ we have
\[ f^{(a)} (\gm x) \= \sum_{0\leq b \leq a} (i\al(\gm))^{a-b}\, \binom
a b\, f^{(b)}(x)\,,\]
where $b$ runs over the multi-indices with $0\leq b_j \leq a_j$, where
$\binom b a = \prod_j \binom{a_j}{b_j}$, and where
$\al(\gm)^c = \prod_j \al_j(\gm)^{c_j}$. Hence
\begin{equation}\label{a-dec} f^{(a)} | (\gm-1) \= \sum_{0\leq b < a}
(i \al(\gm))^{a-b} \, \binom ab\, f^{(b)}\end{equation}
is a linear combination of higher order forms $f^{(b)}$ of orders
$1, \dots, |a|$. So $f^{(a)}$ is an invariant of order at
most~$1+|a|$. Furthermore
\begin{equation} \label{a-dec1}
 f^{(a)}|(\gm_1-1)\cdots (\gm_{|a|}-1)
\= \sum_{0\leq b < a} (i \al(\gm_1))^{a-b} \, \binom ab\, f^{(b)} |
(\gm_2-1)\cdots
(\gm_{|a|}-1)\,.
\end{equation}
By induction assumption, the $f^{(b)}|(\gm_2-\nobreak 1)\cdots
(\gm_{|a|}-\nobreak 1)$ are multiples of $f$ (zero if $|b|<|a|-1$). So
$f^{(a)}$ is a perturbation of~$f$.

For the commutativity of the perturbation we note by induction that,
for all $g_1, \dots, g_s \in \Gm$
\begin{align*}
(g_1-1)(g_2-1)\cdots(g_s-1)=\sum_{l=0}^s (-1)^{s-l}
\sum_{i_1<i_2<\cdots
<i_l}( g_{i_1}g_{i_2}\cdots
g_{i_l}-1)\,,
\end{align*}
where the $i_j$ run through the set $\{1,\ldots,s\}$. Application of
\eqref{a-dec} leads to
\begin{align*}
f^{(a)} \bigr|(\gm_1-1) \cdots (\gm_{|a|}-1)
= \sum_{l=0}^{|a|} (-1)^{|a|-l} \sum_{i_1<i_2<\cdots<i_l}
\sum_{0\leq b <|a|} \bigl( i
\al(\gm_{i_1}\gm_{i_2}\cdots\gm_{i_l}))^{a-b} \, \binom a b \,
f^{(b)}\,.
\end{align*}
Since $\al$ is a homomorphism, the factor
$\al(\gm_{i_1}\gm_{i_2}\cdots\gm_{i_l})$ does not depend on the order
of the $\gm_{i_j}$. Hence we may rewrite the expression as follows.
\begin{align*}
f^{(a)} \bigr|(\gm_1-1) \cdots (\gm_{|a|}-1)
= \sum_{l=0}^{|a|} \frac{(-1)^{|a|-l}}{l!} \sum_{\ii} \sum_{0\leq b
<|a|} \bigl( i \al(\gm_{i_1}\gm_{i_2}\cdots\gm_{i_l}))^{a-b} \, \binom
a b \, f^{(b)}\,,
\end{align*}
where $\ii$ in the sum $\sum_{\ii}$ runs over the subsets of
$\{1,\ldots,|a|\}$ with $l$ elements. This is an expression that is
invariant under permutations of the $\gm_j$, which shows that
$f^{(a)}$ is a commutative perturbation.
\end{proof}
\noindent \emph{Remark. } Proposition~\ref{prop-diff} shows that
commutative perturbations can arise 
as infinitesimal perturbations of a family of automorphic forms. That
is our motivation to use the word \emph{perturbation} in
Definition~\ref{pert-def}.\medskip

\noindent\bf Application to harmonic perturbations of~$1$. \rm We use
the method of differentiation of families to produce explicit
harmonic higher order forms for~$\Gcom$ of order~$3$. We employ
families studied in~\cite{Br94}.
\medskip

Since $\Gcom$ is free on the generators $C=\pm\matr2{-1}{-1}1$ and
$D=\pm\matc2111$, the character group of $\Gcom$ is isomorphic to
$\CC^\ast\times\CC^\ast$. We can parametrise the characters by
\begin{equation}\label{ch-vw}
\ch_{v,w}(\gm) \= e^{i v \ld(\gm) + i w \overline{\ld(\gm)}}\,,
\end{equation}
where $(v,w)$ runs through~$\CC^2$, and where $\ld\in \hom(\Gcom,\CC)$
is as defined in \eqref{H-ld-def}. We are interested only in $(v,w)$
in a neighborhood of~$0\in \CC^2$.

In \cite{Br94}, \S15.5 it is shown that there is a meromorphic
Eisenstein family $E(v,w,s)$ of automorphic forms for $\Gcom$, with
the character~$\ch_{v,w}$ and eigenvalue $\frac14-s^2$ for
$\om_0=-y^2\left( \partial_x^2+\partial_y^2\right)$.
(In \cite{Br94} the discussion of the family $E$ is made in the
context of families of automorphic forms of varying weight which are
thus defined on the covering group~$\tGcom$. However, in \S15.5 the
weight is zero, and the automorphic forms are, in effect, on the
discrete group~$\Gcom$.)
The restriction to $s=\frac12$ exists (\cite{Br94}, \S15.6) and forms
a meromorphic family $(v,w)\mapsto f(v,w;z)$ on $\CC^2$ such that
$f(v,w;\gm z) = \ch_{v,w}(\gm)\, f(v,w;z)$, and
$L_0 f(v,w;z)  =0$ 
for the dense set of $(v,w)$ at which $f$ is holomorphic. There is a
meromorphic family $(v,w)
\mapsto h(v,w;\cdot)$ on~$\CC$, such that
$f(v,w;z) = h\bigl(v,w;H(z)\bigr)$, satisfying
$h(v,w;u+\nobreak\ld) = e^{iv\ld+iw\bar\ld} \, h(v,w;u)$
(\cite{Br94}, \S15.1--6). Chapter~15 of~\cite{Br94} gives a
complicated but explicit construction
(obtained with the help of D.Zagier)
of such a family $h$ with Jacobi theta-functions.

Specifically, in \S15.6.11 the function $h$ is expressed as a sum
\begin{equation}\label{h-def}
h(v,w;u) \= G_{(v+w)\varpi/2\pi}(u,w)+ G_{-(v+w)\varpi/2\pi}(-\bar
u,-v)\,,
\end{equation}
where the function $G_\mu(u,w)$, for $\mu\not \in \ZZ$ and $0<\im u<
\frac12\varpi\sqrt 3$ is given by
\begin{equation}\label{G-def}
G_\mu(u,w) \= \sum_{m=-\infty}^\infty \frac 1{\mu+m}\;
\frac{\x^{\mu+m}}{\eta\, q^m-1}\,,
\end{equation}
with $q=-e^{-\pi \sqrt 3}$, $\x=e^{2\pi i u/\varpi}$, and
$\eta=e^{-w\varpi\sqrt 3}$ We consider this for $u$, $w$, and $\mu$
near zero, but not equal to zero. Hence $\eta\approx 1$ but
$\eta\neq 1$, and $|q|<|\x|<1$. The latter inequalities imply
absolute convergence of the series. We shall derive the Taylor
expansion of $\tilde h(v,w;u):=vw\, h(v,w;u)$ in terms of $(v,w)$
near zero up to order two, from which we can obtain higher order
forms by Proposition~\ref{prop-diff}.

The term of $G_{\mu}(u, w)$ with $m=0$
\begin{equation}
\frac 1\mu \, \frac{ \x^\mu }{\eta-1}\,,
\end{equation}
has singularities at $\mu=0$, and, due to $\frac1{\eta-1}$, also at
$w=0$. This term has the following contribution to $h(v,w;u)$
in~\eqref{h-def}.
\begin{equation}
\frac{2\pi}{\varpi(v+w)} \frac{e^{ i u\,(v+w)}} {e^{-w\varpi\sqrt3}-1}
- \frac{2\pi}{\varpi(v+w)} \; \frac{e^{ i \bar u \, (v+w)}}{
e^{v\varpi\sqrt 3}-1}
\end{equation}
We write the corresponding contribution to $\tilde h(v,w;u)=vw\,
h(v,w;u)$ as follows.
\begin{align*}
\frac{2\pi}{\varpi}&
\frac{v\,w}{(e^{-w\varpi\sqrt3}-1)(e^{v\varpi\sqrt 3}-1)}\;\biggl(
\frac{e^{-w\varpi\sqrt 3}(e^{(v+w)\varpi\sqrt 3}-1)}{v+w}
\\
&\qquad\hbox{} + \frac{
(e^{ i u (v+w)}-1)(e^{v\varpi\sqrt 3}-1)}{v+w} - \frac{(e^{ i \bar
u(v+w)}-1)(e^{-w\varpi\sqrt 3}-1)
}{v+w} \biggr) \,.
\end{align*}
The last three quotients are holomorphic as a function of $v+w$ in a
neighborhood of~$0$. We replace them by their Taylor expansion up to
the term $(v+\nobreak w)^2$ and after that the Taylor expansion in
both $v$ and $w$ up to order~$2$ is computed. This gives
\begin{equation}
\begin{aligned}
\frac{-2\pi}{\varpi^2\sqrt 3} & \biggl( 1 + iu\, v + i \bar u\, w -
\frac12 u^2\, v^2- \frac12 \bar u^2 w^2 \\
&\qquad\hbox{}- \frac{\sqrt 3}{2\pi}\bigl( -\frac{\pi \varpi^2}{2\sqrt
3} -\pi i \varpi u +\pi i \varpi \bar u + \frac\pi{\sqrt 3}(u^2+\bar
u^2)
\bigr)\, vw \biggr)+ \cdots\,.
\end{aligned}
\end{equation}
\medskip

In the terms with $m\neq 0$ in \eqref{G-def} we write $\x=e^{2\pi i
u/\varpi}$, $\tilde \x=e^{-2\pi i \bar u/\varpi}$,
$\eta_1=e^{-w\varpi\sqrt 3}$, $\eta_2=e^{v\varpi\sqrt 3}$,
$q=-e^{-\pi\sqrt 3}$, and $\mu=(v+w)\varpi/2\pi$. We find the
following contribution to $\tilde h(v,w;u)$:
\[ \sum_{m=1}^\infty \biggl( \frac{vw}{m+\mu} \,
\frac{\x^{m+\mu}}{\eta_1 q^m-1}
+ \frac{vw}{\mu-m} \; \frac{\x^{\mu-m}}{\eta_1 q^{-m}-1} +
\frac{vw}{m-\mu}\, \frac{\tilde \x^{m-\mu}}{\eta_2 q^m-1} +
\frac{vw}{-m-\mu}\; \frac{\tilde \x^{-m-\mu}}{\eta_2 q^{-m}-1}
\biggr)\,.\]
This contribution is holomorphic near $v=w=0$. Its expansion starts
with the term $vw$. So for third order forms we need only the
contribution to $h(0,0;u)$:
\begin{equation}
\sum_{m=1}^\infty \frac 1 m \biggl( \frac{\x^m}{q^m-1} + \frac{
(q/\x)^m}{q^m-1} + \frac{\tilde\x^m}{q^m-1} +
\frac{(q/\tilde\x)^m}{q^m-1} \biggr)\,.
\end{equation}
Each of these terms gives a convergent series on the region $0<\im u <
\frac12\varpi\sqrt 3$.

\rmrk{Commutative perturbations}In this expansion we find various
higher order harmonic forms that we have seen above. Denoting
$f=\frac{-2\pi}{\varpi^2\sqrt 3}$ we find:
\begin{equation}\label{com-tab}
\begin{array}{|c|c|l|}
\hline
\text{term of}&\text{ on $\CC$} & \text{ on $\uhp$}\\ \hline
1& f & f \ (\text{constant function})
\\
v& if\, u & if\,H(z)
\\
w& i f\, \bar u & if\,\overline{H(z)}
\\
v^2& \frac{-f}2 \,u^2 &\frac{-f}2\, H(z)^2
\\
w^2& \frac{-f}2\, \bar u^2 &\frac{-f}2\, \overline{H(z)}^2
\\ \hline
\end{array}
\end{equation}
The coefficient of $vw$ gives a third order form
\begin{equation}\label{b11}
\begin{aligned}
b_{1,1} (u)&\;:=\; \frac\pi{\sqrt 3} \biggl( \bigl( \frac u\varpi -
\frac{i\sqrt 3}2\bigr)^2 + \bigl( \frac{\bar u}\varpi+\frac{i\sqrt
3}2\bigr)^2+1\biggr)\\
&\qquad\hbox{} + S(u) + S(\varpi\rho- u) + S(-\bar u) + S(\varpi\rho
+\bar u)\,,
\\
\text{with}\quad S(u)&\;:=\; \sum_{m=1}^\infty \frac{e^{2\pi i m
u/\varpi}}{m\,(q^m-1)}\,,\qquad \rho\= \frac12+\frac i2\sqrt 3\,.
\end{aligned}
\end{equation}
By $B_{1,1}(z)=b_{1,1}\bigl(H(z)\bigr)$ we denote the corresponding
harmonic third order perturbation of~$1$ on~$\uhp$. The way $B_{1,1}$
has been derived, together with the proof of
Proposition~\ref{prop-diff}, ensures that it is a perturbation of~$1$
with a multilinear form that is a multiple of
$\ld \otimes \bar\ld+\bar\ld\otimes\ld$.

However, $b_{1,1}(u)$ is represented by \eqref{b11} only on the region
$0<\im u
< \frac12\varpi\sqrt 3$. In \cite{Br94}, \S15.3.5, the image under $H$
of the fundamental domain
\vskip -4mm
\[\bigcup_{n=-2}^3 \matc 1n01 \fd_{\mathrm{mod}}\]
\vskip -1mm
\noindent
(where $\fd_{\mathrm{mod}}$ is the standard fundamental domain of the
modular group)
is shown to be the regular hexagon with centre $0$ and one corner at
$-\frac13(e^{\pi i/3}+\nobreak 1)\varpi$. 
Only the upper half of this hexagon is in the region where we have an
expression for~$b_{1,1}$. We shall continue this function to the
entire~$\CC$.

We first note that the series in \eqref{b11} defining $S(u)$ converges
absolutely for $\im u>0$ yielding a holomorphic function in that
region. To extend $S(u)$ to other values we use the following
identity, valid for $\im u>\frac12\varpi\sqrt 3$:
\begin{equation}\label{Sext}
S(u) \= \sum_{m=1}^\infty \frac{e^{2\pi i m (u/\varpi +
\rho)}}{m(q^m-1)} - \sum_{m=1}^\infty \frac{e^{2\pi i m u/\varpi}}m
\= S(u+\varpi\rho) + \log\bigl(1-e^{2\pi i u/\varpi}\bigr)\,.
\end{equation}
Via this identity, we can define $S(u)$ in the region $\im
u>-\frac12\varpi\sqrt 3$. This extension of $S$ is multivalued, since
it depends on the way in which we extend the function $u\mapsto
\log\bigl(1-\nobreak e^{2\pi i u/\varpi}\bigr)$, which is given by
the second series in~\eqref{Sext} only for $\im u>0$. However, the
sum
\begin{equation} \label{exten}
S(u)+S(-\bar u) \= S( u + \varpi\rho) + S\bigl( -
\overline{u}+\varpi\rho \bigr)
+ 2\log |1-e^{2\pi iu/\varpi}|
\end{equation}
is single-valued on $\im u>-\frac12\varpi\sqrt 3$, with logarithmic
singularities at~$u=\varpi n$, $n\in \ZZ$. Applying \eqref{exten}
repeatedly, we can extend $S(u)+S(-\bar u)$ to all $\CC$ to obtain a
harmonic function with singularities at the points in $\Ld=\varpi
\,\ZZ[\rho]$ which have non-positive imaginary part.

Via \eqref{b11}, we then obtain the continuation of the function
$b_{1,1}$. It is harmonic on $\CC\setminus\Ld$, with logarithmic
singularities at all points of $\Ld$.

Let us explicitly check the transformation behaviour: Since $S$ is
periodic with period $\varpi$
(and, equivalently, $S(u+\varpi \rho) \= S(u-\varpi \bar \rho)$ ),
\begin{align*}
b_{1,1}(u&+\varpi) -b_{1,1}(u) \= \frac\pi{\sqrt 3} \biggl( \bigl(
\frac u\varpi+1-\frac{i\sqrt 3}2\bigr)^2
-\bigl( \frac u\varpi-\frac{i\sqrt 3}2\bigr)^2 \\
&\qquad\hbox{}+ \bigl( \frac{\bar u}\varpi+1+\frac{i\sqrt 3}2 \bigr)^2
- \bigl( \frac{\bar u}\varpi+\frac{i\sqrt 3}2 \bigr)^2 \biggr)
\displaybreak[0]\\
&\= \frac\pi{\sqrt 3} \biggl( \frac {2u}\varpi +1-i\sqrt 3
+ \frac{2\bar u}\varpi+1
+ i \sqrt 3 \biggr) \= \frac{2\pi}{\varpi\sqrt 3}(u+\bar
u+\varpi)\,;\\
b_{1,1}(u&+\varpi\rho)-b_{1,1}(u)
\= \frac\pi{\sqrt 3}\biggl( \bigl( \frac u\varpi +\frac12\bigr)^2 -
\bigl( \frac u\varpi - \frac{i\sqrt 3}2 \bigr)^2 \\
&\qquad\qquad\hbox{}
+ \bigl( \frac{\bar u}\varpi
+\frac12\bigr)^2 - \bigl( \frac{\bar u}\varpi + \frac{i\sqrt
3}2\bigr)^2 \biggr) \\
&\qquad\hbox{}- 2 \log\bigl| 1-e^{2\pi i u/\varpi}\bigr| + S(-u) +
S(\bar u+\varpi) - S(\varpi \rho-u) - S(\varpi\rho+\bar
u)\displaybreak[0]\\
&\= \frac\pi{\sqrt 3} \biggl(
\rho\bigl(\frac{2u}\varpi+\rho^{-1}\bigr) + \rho^{-1} \bigl(
\frac{2\bar u}\varpi + \rho\bigr) \biggr)\\
&\qquad\hbox{}- 2 \log\bigl| 1-e^{2\pi i u/\varpi}\bigr| + 2
\log\bigl| 1- e^{-2\pi i u/\varpi}\bigr|
\displaybreak[0]\\
&\= \frac{2\pi}{\sqrt 3} \bigl( 1 + \frac{\rho u+\overline{\rho
u}}\varpi\bigr)-2\pi i(u-\bar u)/\varpi \= \frac{2\pi}{\sqrt 3}
\bigl( 1+\rho^{-1}\, \frac u\varpi+\rho\, \frac{\bar
u}\varpi\bigr)\,.
\end{align*}
Let us denote by $T_\om$ the translation by $\om \in \Ld$, and use the
notations $b_{1,0}(u)=u$, $b_{0,1}(u) = \bar u$. With the notations
$f=\frac{-2\pi}{\varpi^2\sqrt 3}$ and $a= \frac{2\pi}{\varpi\sqrt 3}
= -f\varpi$ we have 
\begin{align}\label{b11-trf}
b_{1,0}|(T_{\varpi} -1)&\= \om\,,\qquad b_{0,1}|(T_\om-1)\=
\bar \om\,,
\displaybreak[0]\\
\nonumber
b_{1,1}|(T_\varpi-1)&\= a\,\bigl( b_{1,0}+b_{0,1}+\varpi)\,,
\displaybreak[0]\\
\nonumber
b_{1,1}|(T_{\rho\varpi}-1)&\= a\,\bigl( \bar\rho b_{1,0} + \rho
b_{0,1} + \varpi)\,,
\displaybreak[0]\\
\nonumber
b_{1,1}|(T_\varpi-1)^2&\= 2a\, \varpi \= -2f\, \varpi^2\,,
\displaybreak[0]\\
\nonumber
b_{1,1}|(T_\varpi-1)(T_{\rho\varpi}-1)&\=
a\,\bigl(\rho\varpi+\bar\rho\varpi \bigr)
\=-f \,\bigl( \bar \varpi\cdot \rho\varpi + \varpi \cdot
\overline{\rho\varpi}\bigr)\,,
\displaybreak[0]\\
\nonumber
b_{1,1}|(T_{\rho\varpi}-1)^2&\=2a\, \varpi \= - 2f (\rho\varpi)
(\overline{\rho\varpi})\,.
\end{align}
Since $\Ld$ is commutative we need not consider
$b_{1,1}|(T_{\rho\varpi}-\nobreak1)(T_\varpi-\nobreak1)$. We conclude
that the pull-back $-f^{-1} B_{1, 1}=-f^{-1}b_{1, 1}\circ H$ is a
harmonic commutative perturbation of~$1$ for the multilinear form
$\mu$ determined by the following values at the generators $C$ and
$CD$ of $\Gcom$:
\begin{equation} \mu(g, h)\= \begin{cases} 2 \varpi^2 & \text{ if }
g=h=C \text{ or } CD\,, \\
\varpi^2 &\text{ if } g=C, \; h=CD\,, \; \text{ or if } g=CD, \; h=C
\,.
\end{cases}
\end{equation}
We have used the values of $\ld$ given below~\eqref{H-ld-def}.)
With these values at the generators, $\mu$ coincides with
$\ld\otimes\bar\ld+\bar\ld\otimes\ld$ as predicted above by the way
$B_{1, 1}$ was constructed.

\rmrk{Non-commutative perturbation} Proposition~\ref{prop-diff} shows
that differentiation of families produces only commutative
perturbations. However, by Theorem~\ref{thm-Mf-mp}, there are
non-commutative third order harmonic perturbations of~$1$. We can
obtain such perturbations from $B_{1,1}$ upon decomposing it as
$B_{1,1}=A+B$ for a holomorphic function $A$ and an anti-holomorphic
function $B$.

Specifically, in view of \eqref{b11}, for those $z\in \uhp$ for which
$H(z)$ is in the upper half of the fundamental hexagon for~$\CC/\Ld$,
we can set
\begin{equation}\label{AB}
\begin{aligned}
A(z) &\= \frac\pi{2\sqrt 3} + \frac\pi{\sqrt 3}\biggl( \frac
{H(z)}\varpi - \frac{i\sqrt 3}2\biggr)^2 + S\bigl(H(z)\bigr) +
S\bigl(\varpi \rho-H(z) \bigr)\,,\\
B(z) &\= \frac\pi{2\sqrt 3}+\frac\pi{\sqrt 3}\biggl(
\frac{\overline{H(z)}}\varpi+\frac{i\sqrt 3}2\biggr)^2 + S\bigl(
-\overline{H(z)} \bigr) + S \bigl( \varpi\rho+\overline{H(z)}\bigr)\,.
\end{aligned}
\end{equation}

As shown above, $B_{1, 1}|(\gm-1)(\dt-1)=-f \, \ld \otimes \bar \ld -
f\, \bar\ld\otimes \ld$. Hence,
\[ A|(\gm-1)(\dt-1) \= -B|(\gm-1)(\dt-1) -f \, \ld \otimes \bar \ld -
f\, \bar\ld\otimes \ld. \]
gives an equality between a holomorphic and an antiholomorphic
function, and therefore, there is $\nu:\Gm^2\rightarrow \CC$ such
that
\[ A|(\gm-1)(\dt-1)\= \nu(\gm,\dt)\,, \, \, B|(\gm-1)(\dt-1) \=
-f \, \ld \otimes \bar \ld - f\, \bar\ld\otimes \ld
-\nu(\gm,\dt) \]
for all $\gm,\dt\in \Gm$. This implies that $A$ and $B$ are third
order invariants, and that $\nu \in \mult^2(\Gm,\CC)$.

To determine the bilinear form $\nu$, we recall that $\ld(C)= \rho
\varpi$ and $\ld(D) = \bar\rho \varpi=(1-\nobreak\rho)\varpi$. We
consider the following four functions:
\begin{equation}
\label{AB-tab}
\begin{array}{ccc}
A|(C-1)+f\bigl( \bar\rho\varpi \, H + \frac{\varpi^2}2\bigr)\,,
&& B|(C-1)+ f\bigl(\rho\varpi \, \bar H + \frac{\varpi^2}2\bigr)\,,
\\
A|(D-1)+f\bigl(\rho\varpi \, H + \frac{\varpi^2}2\bigr)\,,
&& B|(D-1)+ f\bigl(\bar \rho\varpi\, \bar H+ \frac{\varpi^2}2\bigr)\,,
\end{array}\end{equation}
The functions on the left are holomorphic, and those on the right are
antiholomorphic. We consider the sum of the two functions on the
first row, and denote $u=H(z)$. With~\eqref{b11-trf}:
\begin{align*}
B_{1,1}|(C-1)(z) &+ f\bar\rho \varpi H(z) + f\rho\varpi\overline{H(z)}
+ f\varpi^2
\\
&\= b_{1,1}|(T_{\rho\varpi}-1)(u)+f\varpi (\bar \rho u + \rho \bar u)
+ f\varpi^2\\
&\= -f\varpi\bigl( \bar \rho u +\rho \bar u+\varpi\bigr)+ f\varpi
(\bar \rho u + \rho \bar u) + f\varpi^2\=0\,.
\end{align*}
Similarly the sum of the two functions on the second row gives
\begin{align*}
B_{1,1}|&(D-1)(z) + f\bigl( \rho\varpi H(z) + \bar \rho \varpi
\overline{H(z)} + \varpi^2\bigr)\\
&\= b_{1,1}|(T_{\bar\rho\varpi}-1)(u) + f\varpi(\rho u +
\overline{\rho u}+\varpi)\\
&\= b_{1,1}|(T_\varpi-1)T_{\rho\varpi}^{-1} -
b_{1,1}|(T_{\rho\varpi}-1)T_{\rho\varpi}^{-1}+ f\varpi(\rho u +
\overline{\rho u}+\varpi)\\
&\= \bigl( -f\varpi(u+\bar u+\varpi) +f \varpi (\bar\rho u+\rho \bar
u+\varpi)
\bigr)|T_{\rho\varpi}^{-1}+ f\varpi(\rho u + \overline{\rho
u}+\varpi)\\
&\= f\varpi\bigl( -u+\rho\varpi -\bar u + \bar\rho\varpi -\varpi +
\bar\rho( u - \rho\varpi) + \rho(\bar u -\bar \rho\varpi) + \varpi \\
&\qquad\hbox{}+\rho u + \overline{\rho u}
+ \varpi\bigr) \= 0\,.
\end{align*}
The sums of the rows in~\eqref{AB-tab} are zero, so the individual
functions are constant. We do not try to determine these constants.

For~$A$ we have
\begin{align*}
A|(C-1)(C-1)&\= -f\,\bigl( \bar \rho \varpi H|(C-1) +0\bigr)\= -f
\bar\rho \varpi \ld(C) \= -f \overline{\ld(C)}\, \ld(C)\,,\\
A|(C-1)(D-1)&\= -f \bigl( \bar \rho \varpi \ld (D) \bigr) \= -f
\overline{\ld(C)}\, \ld(D)\,,\\
A|(D-1)(C-1) &\= - f \overline{\ld(D)}\, \ld(C)\,,\\
A|(D-1)(D-1) &\= -f \overline{\ld(D)}\, \ld(C)\,.
\end{align*}

We conclude that $-f^{-1}\, A$ is a non-commutative holomorphic third
order holomorphic perturbation of~$1$ with multilinear form $ \bar\ld
\otimes \ld$. Then the multilinear form of the anticommutative third
order perturbation of~$1$ given by
$-f^{-1}\,B=-f^{-1}(B_{1,1}-\nobreak A)$ is $\bigl( \ld  \otimes
\bar \ld + \bar \ld\otimes \ld\bigr) - \bar\ld\otimes \ld = \ld
\otimes \bar\ld $.

\section{Universal covering group}
\subsection{Universal covering group of $\SL_2(\RR)$} \label{Ucg}
To define the universal covering group of $\SL_2(\RR)$, which is
also the universal covering group of $G=\PSL_2(\RR)$, we first note
that, as an analytic variety, $\SL_2(\RR)$ is isomorphic to
$\uhp\times \bigl( \RR/2\pi\ZZ\bigr)$, by the Iwasawa decomposition
expressing each element of $\SL_2(\RR)$ uniquely as a product
\[ \matc{\sqrt y}{x/\!\!\sqrt y}0{1/\!\!\sqrt y} \,
\matr{\cos\th}{\sin\th}{-\sin\th}{\cos\th}\,,\]
with $x+iy\in \uhp$ and $\th \in \RR/2\pi\ZZ$. Left multiplication by
$\matc abcd \in \SL_2(\RR)$ amounts to
\begin{equation}
\label{g-def}
(z,\th+2\pi\ZZ) \mapsto \Bigl( \frac{az+b}{cz+d}, \th-\arg\bigl(
j\bigl(\matc abcd, z\bigr)\bigr)+2\pi \ZZ\Bigr)\,.\end{equation}
Here, $j\bigl(\matc abcd, z\bigr):=cz+d$. This describes an action of
$\SL_2(\RR)$ on $\uhp\times\bigl(  \RR/2\pi \ZZ\bigr)$.

We define for each $\matc abcd \in \SL_2(\RR)$ 
the operator
\begin{equation}\label{tildeg-def}
\widetilde{\matc abcd} : (z,\th) \mapsto \bigl( \frac{az+b}{cz+d},
\th-\arg \bigl( j\bigl(\matc abcd, z\bigr)\bigr) \bigr)
\end{equation}
from $\uhp\times\RR$ to itself, where we choose the argument such that
$-\pi<\arg(cz+\nobreak d)\leq \pi$. We note that the map $g\mapsto
\tilde g$ is injective.

\begin{defn} \label{ucg} The universal covering group $\tG$ of $G$ is
the group of operators $\uhp\times\RR\rightarrow\uhp\times\RR$
generated by the operators $\tilde g$ in~\eqref{tildeg-def} for all
$g\in \SL_2(\RR)$.
\end{defn}

A lengthy but routine calculation shows
\begin{lem}\label{lem-compat} If the vertical maps in the diagram
$$\begin{CD}
\uhp \times \RR @>>> \uhp \times \RR/\ZZ @= \SL_2(\RR)\\
@V\tilde g VV @VVV @VgVV\\
\uhp \times \RR @>>> \uhp \times \RR/\ZZ @= \SL_2(\RR)\\
\end{CD}$$
are given by \eqref{tildeg-def}, \eqref{g-def} and by left
multiplication by $g$ respectively, then the diagram is commutative.
(The last horizontal maps are defined by the Iwasawa decomposition.)
\end{lem}

Suppose now that $\tilde g_1 \tilde g_2\cdots\tilde g_n $ is the
identity as an operator on $\uhp\times\RR$. Then $z\mapsto
g_1g_2\cdots g_n z$ is the identity operator on~$\uhp$. So
$g_1g_2\cdots g_n \in \{I,-I\}\subset \SL_2(\RR)$.
By Lemma~\ref{lem-compat}, it is impossible that $g_1g_2\cdots
g_n=-I$ while $\tilde g_1\tilde g_2\cdots \tilde g_n $ is the identity
operator. So $g_1g_2\cdots g_n=I$. This implies that the map $\tilde
g \mapsto g$ on the generators extends to a group homomorphism
\[ \pr_2:\tG\longrightarrow\SL_2(\RR)\,.\]
The composition of $\pr_2$ with the natural projection $\SL_2(\RR) \to
\PSL_2(\RR)$ gives a map
\[ \pr:\tG\longrightarrow\PSL_2(\RR)\,.\]

\medskip

We single out the following following families of elements of $\tG$.
\begin{enumerate} \item[a)] For $x\in \RR$ we put $n(x) =
\widetilde{\matc1x01}$ in~$\tG$. This induces an injective group
homomorphism $n:\RR\rightarrow\tG$.
\item[b)] For $y \in \RR_+^\ast$ we set $a(y) =
\widetilde{\matc{y^{1/2}}00{y^{-1/2}}}$. This induces an injective
group homomorphism $a:\RR_+^\ast \rightarrow \tG$.
\item[c)] For $\th \in \RR$, we set
\begin{equation} \label{k(th)}
k(\th) (z,\th_1) \= \bigl( \frac{z\cos\th+\sin \th}{-z\sin
\th+\cos\th},\th_1+\th-\arg(e^{i \th}(-z\sin \th+\cos\th))
\bigr).
\end{equation}
This defines $k(\th) \in \tG$ satisfying $\pr_2 k(\th) =
\matr{\cos\th}{\sin\th}{-\sin\th}{\cos\th}$. For fixed $(z,\th_1)\in
\uhp\times\RR$, the quantity $k(\th)(z,\th_1)$ is real-analytic
in~$\th$. If both $\th$ and $\th'$ have values near zero then
 $k(\th+\nobreak \th')=k(\th)k(\th')$, since $\pr_2$ is locally an
isomorphism. By analyticity this relation extends to all
$\th,\th'\in\RR$. So we have a group homomorphism
$k:\RR\rightarrow\tG$. The kernel of the composition $\pr_2\circ k$
is $2\pi \ZZ$. For each $n\in \ZZ$ the element $k(n\pi)$ acts as
 $(z,\th_1)\mapsto (z,\th_1+\pi n)$. This implies that $k$ is an
 injective group homomorphism. Although it satisfies $\pr_2 k(\th) =
\matr{\cos\th}{\sin\th}{-\sin\th}{\cos\th}$ for all $\th\in \RR$,
the relation 
$\widetilde{\matr{\cos\th}{\sin\th}{-\sin\th}{\cos\th}}=k(\th)$ holds
only for $\th\in [-\pi,\pi)$.
\end{enumerate}

With these definitions and notations we deduce some basic facts about
$\tG$.
\smallskip

\noindent\it Centre of $\tG$: \rm The elements $k(\pi n)$ with $n\in
\ZZ$ form the centre $\tZ$ of~$\tG$.
\smallskip

\noindent\it Transitivity of action of $\tG$ on $\uhp \times \RR$: \rm
This is implied by $n(x)a(y)k(\th) \, (i,0) = (x+\nobreak iy,\th)$
for all $x+iy\in \uhp$ and $\th\in \RR$.
\smallskip

\noindent
\it Generators of $\tG$: \rm The elements $n(x)$, $a(y)$ and $k(\th)$
generate $\tG$, and each element of $\tG$ can be written uniquely as
$n(x)a(y)k(\th)$. This is a consequence of the relations
\begin{eqnarray}
a(y)n(x)&=&n(y^2x)a(y) \qquad \text{and} \\
k(\th)n(x)a(y) &=& n(x_\th)a(y_\th)k\bigl(\th-\arg\bigl(
e^{i\th}(-z\sin\th+\cos\th)\bigr)\,\bigr)
\end{eqnarray}
with $z=x+iy$ and $x_\th+iy_\th =
\frac{z\cos\th+\sin\th}{-z\sin\th+\cos\th}$.
\smallskip

\noindent
\it $\tG \cong \uhp\times\RR$. \rm Because of the last two facts, we
can identify $\tG$ with $\uhp\times\RR$ as analytic varieties.
Furthermore, the group operations are analytic with respect to the
structure of $\uhp\times\RR$ as an analytic variety. So $\tG$ is a
Lie group. The maps $\pr$ and $\pr_2$ are covering maps. One can show
that any covering of $\SL_2(\RR)$ factors through $\tG$, hence $\tG$
is the universal covering group of~$\SL_2(\RR)$.
\smallskip

\noindent
\it Section $g \to \tilde g$: \rm This is a homeomorphism for $g$ near
the unit element of $\SL_2(\RR)$, but it is discontinuous at $\matc
abcd \in \SL_2(\RR)$ with $c=0$ and $d<0$. This section is not a
group homomorphism but instead there is a $\ZZ$-valued $2$-cocycle
$w$ on $\SL_2(\RR)$ such that $\tilde g \tilde g_1 = \widetilde{gg_1}
k\bigl(2\pi w(g,g_1) \bigr)$ for all $g,g_1\in \SL_2(\RR)$. See
Theorem~16 on p.~115 of~\cite{Ma} for an explicit description of this
cocycle. Each element of $\tG$ has a unique decomposition as
$\tilde g\, k(2\pi n)$ with $g\in \SL_2(\RR)$ and $n\in \ZZ$. In this
paper we will not use this description of the group structure of
$\tG$. We work with the interpretation as a group of operators in
$\uhp\times\RR$, and occasionally use the ``one-parameter subgroups''
$n$, $a$ and~$k$.
\smallskip

\noindent
\it The action of $\tG$ \rm on $\uhp^*:=\uhp \cup \{\text{cusps}\}$ is
given by $\gm z:= \pr(\gm)\, z.$

\vskip -5mm

\subsection{The Lie algebra of the universal covering group}
\label{Lie} The direction of the three one-parameter subgroups
$n$, $a$ and $k$ 
at the origin determines elements of the (real) Lie algebra
$\glie_\RR$ of~$\tG$. The groups $\tG$, $\SL_2(\RR)$ and
$\PSL_2(\RR)$ have the same Lie algebra, since they are locally
isomorphic. The Lie algebra elements corresponding to $n$, $a$ and
$k$ are, respectively,
\begin{equation}
\XX\= \matc0100\,,\qquad \frac12 \HH\=\matr{1/2}00{-1/2}\,,\quad
\text{ and }\quad \WW\=\matr01{-1}0\,.
\end{equation}
The Lie algebra acts on the functions on $\tG$ by differentiation on
the right: $\YY F (g) = \partial_t F(g\exp(t\YY))|_{t=0}$ for
$\YY\in \glie_\RR$. This action can be extended to the complexified
Lie algebra $\glie=\CC\otimes_\RR \glie_\RR$, and to the universal
enveloping algebra of~$\glie$. All the resulting differential
operators commute with the action of~$\tG$ by left translation. With
the identification of $\tG$ as $\uhp\times \RR$ we have in the
coordinates given by $(x+\nobreak iy,\th)$:
\begin{equation}\label{Lie-diff}
\begin{aligned}\XX&\= \partial_x\,,\qquad \HH\= 2y\partial_y\,,\qquad
\WW\=\partial_\th\,,\\
\EE^+& \=\HH+i(2\XX-\WW)\=
e^{2i\th}(2iy\partial_x+2y\partial_y-i\partial_\th)\,,\\
\EE^-&\=\HH-i(2\XX-\WW)\=
e^{-2i\th}(-2iy\partial_x+2y\partial_y+i\partial_\th)\,,\\
\om&\= - \frac14\EE^\pm\EE^\mp+\frac14\WW^2\mp \frac i2 \WW \=
-y^2\partial_y^2-y^2\partial_x^2+y\partial_x\partial_\th\,.
\end{aligned}
\end{equation}
The \emph{Casimir operator} $\om$ generates the centre of the
enveloping algebra of~$\glie$. The corresponding differential
operator commutes with left and right translations in~$\tG$.

\subsection{Cofinite discrete subgroups}\label{sect-cdsg}To a cofinite
discrete subgroup $\Gm$ of $\PSL_2(\RR)$ we associate its full
original $\tGm:=\pr^{-1}\Gm$ in $\tG$. This gives a bijective
correspondence between cofinite discrete subgroups of $\PSL_2(\RR)$
and cofinite discrete subgroups of~$\tG$ that contain the centre
$\tZ=\langle \z \rangle,$ where $\z:=k(\pi)$. The projection $\pr $
induces an isomorphism $\Gamma \cong \tGm/\tZ.$

As an example we consider the \emph{modular group}
$\Gmod=\PSL_2(\ZZ)$, with corresponding group $\tGmod\subset \tG$. It
is known that $\PSL_2(\ZZ)$ is presented by the generators
$S=\pm\matr0{-1}10$ and $T=\pm\matc1101$ and relations
$S^2=(TS)^2=I$.

Set $s:=k(-\pi/2) = \widetilde{\matr0{-1}10}$ and
$t:=n(1)=\widetilde{\matc1101}$ with $\pr (s)=S$ and $\pr (t)=T$.
Then $s^2 = k(-\pi)=\z^{-1} \in \tZ$, so $s$ and $t$
generate~$\tGmod$. The relation $S^2=I$ is replaced by the centrality
of $s^2$. We have
\begin{align*} ts(i,0) &\= t \widetilde{\matr0{-1}10}\,(i,0)
\= \widetilde{\matc1101} \,\bigl(i,-\arg i)
\= \bigl( \frac{i-1}i, -\pi/2-\arg 1\bigr)
= \widetilde{\matr1{-1}10}\, (i,0)\= (i+1,-\pi/2)\,.
\end{align*}
So $ts=\widetilde{\matr1{-1}10}$, and it corresponds to $(i+\nobreak
1,-\pi/2)$ in $\uhp\times\RR\cong \tG$. Hence
\begin{align*}
(ts)^3 &\= ts \widetilde{\matr1{-1}10} \,(i+1,-\pi/2)
\= ts \,\bigl( \frac i{i+1}, -\frac\pi 2- \arg(i+1)\bigr)\\
&\= \widetilde{\matr1{-1}10}\, \bigl( \frac{i+1}2,-\frac{3\pi}4\bigr)
\= \bigl( \frac{i-1}{i+1}, -\frac{3\pi}4-\arg(i-1)\bigr)\\
&\= \bigl(i,-\pi) \= \z^{-1}\= s^2\,.
\end{align*}
The conclusion is that $\tGmod$ has the presentation with generators
$s$ and $t$ and relations $s^2t=ts^2$ and $tstst=s$. This implies
that the linear space $\hom(\tGmod,\CC)$ has dimension~$1$, and is
generated by $\al:t\mapsto \frac\pi 6$, $\al:s\mapsto\frac{-\pi}2$.
For reasons that will become clear later, we take this basis element,
and not an integral-valued one.

\subsection{Canonical generators} \label{cgen}
The canonical generators of $\Gm$ induce canonical generators
of~$\tGm$:
\begin{itemize}
\item Elements $\pi_1,\ldots,\pi_{\npar}$ of the form $\pi_j = \tilde
g_{\k_j} n(1) \tilde g_{\k_j}^{-1}$ fixing a system of
representatives $\k_1,\ldots,\k_{\npar}$ of the $\tGm$-orbits of
cusps.
\item Elements $\e_1,\ldots,\e_{\nell}$ conjugate in $\tG$ to
$k(\pi/v_j)$ with $v_j\geq 2$.
\item Elements $\eta_1,\ldots,\eta_{2g}$ conjugate in~$\tG$ to
elements $a(t_j)$ with $t_j>1$.
\item The generator $\z=k(\pi)$ of the centre $\tZ$ of~$\tGm$.
\end{itemize}
The relations are:
\begin{equation}
\begin{aligned}
&\z\text{ is central}\,,\\
&\e_j^{v_j} = \z\text{ for }1\leq j \leq \nell\,,\\
&\pi_1\cdots\pi_{\npar}\e_1\cdots
\e_{\nell}[\eta_1,\eta_2]\cdots[\eta_{2g-1},\eta_{2g}]\=
\z^{2g-2+\npar+\nell}\,.
\end{aligned}
\end{equation}
The integer $2g-2+\npar+\nell$ is always positive. For these facts
see~\cite{Br94}, \S3.3.

If $\nell>0$ or if $2g-2+\npar=1$ and $\nell=0$, we do not need $\z$
as a generator. If $\nell=0$ the group $\tGm$ is free on
$\pi_1,\ldots,\pi_{\npar-1},\eta_1,\ldots,\eta_{2g},\z$.

Among the canonical generators we single out the following elements:
$\al_1=\pi_1,\ldots,\al_{\npar-1}=\pi_{\npar-1}$, $\al_\npar
=\eta_1,$ $\ldots,\al_{\ngen-1}=\eta_{2g}, \al_\ngen=\z$.
(We recall that $\ngen=\npar+2g$.)
The $\al_j$ together with the $\e_j$ generate~$\tGm$, with
$\e_j^{v_j}=\z$ and the centrality of~$\z$ as the sole relations.
\smallskip

For the \emph{modular group} $\tGmod$ we have $\npar=1$, $\nell=2$,
$g=0$, and $\ngeng\Gmod=1$. We may take $\pi_1=t=n(1)$,
$\e_1=t^{-1}s^{-1}$, and $\e_2=s^{-1}=k(\pi/2) = p^{-1} k(\pi/3)p$,
with $p=n(-1/2)a(\sqrt 3/2)$.
\medskip

By $I$ we now denote the augmentation ideal of the group ring
$\CC[\tGm]$. In $\CC[\tGm]$ we have the elements
\begin{equation}\label{bi-def}
\bb(\ii) \= (\al_{\ii(1)}-1) \cdots (\al_{\ii(q)}-1)\qquad \ii \in
\{1,\ldots,\ngen\}^q\,.
\end{equation}
We allow ourselves to use the same notation as in~\eqref{bi-def0},
since from now on we will use~$\tGm$. The centrality of $\z$ allows
us to move $(\z-\nobreak1)$ through the product. So it suffices to
consider only $q$-tuples~$\ii$ for which all $\ii(l)=\ngen$ occur at
the end. Such $q$-tuples we will call \emph{$\tGm$-$q$-tuples}.
\begin{prop}\label{prop-augpow}
A $\CC$-basis of $I^{q+1}\backslash I^q$ is induced by the elements
\begin{equation}
\bb(\ii) \= (\al_{\ii(1)}-1) \cdots (\al_{\ii(q)}-1)\,,
\end{equation}
where $\ii$ runs over the $\tGm$-$q$-tuples.
\end{prop}
\begin{proof}
The ideal $I^q$ is generated by the products of the form
$(\gm_1-\nobreak1)\cdots(\gm_q-\nobreak 1)$ with
$\gm_1,\ldots,\gm_q\in \tGm$.
(Lemma~1.1 in \cite{Dei}.) With the relation \[(\gm\dt-\nobreak 1)
=(\gm-\nobreak 1)(\dt-\nobreak 1)+(\gm-\nobreak
1)+(\dt-\nobreak1)\,,\]
we can take the $\gm_j$ in a system of generators, for instance
$\al_1,\ldots,\allowbreak\al_\ngen,\allowbreak\e_1,\ldots,\e_\nell$.
For the elliptic elements $\e_j$ we use
$\z-1=\sum_{k=0}^{v_j-1}\e_j^k(\e_j-\nobreak 1)\equiv
v_j(\e_j-\nobreak 1)\bmod I^2$ to see that the $\al_j$ suffice. (Note
that $v_j$ is invertible in~$\CC$.) Since $\al_\ngen=\z$ is central,
we can move all occurrences of $\z-1$ to the right to see that the
$\bb(\ii)$ in the proposition generate $I^{q+1}\backslash I^q$.

To see that the $\bb(\ii)$ are linearly independent over~$\CC$ we
proceed in rewriting terms $\x(\al_{\ii(1)}-\nobreak
1)\cdots(\al_{\ii(q)}-\nobreak 1)$ by replacing $\x\in R:=\CC[\tGm]$
by $n+\eta$ with $n\in \CC$ and $\eta\in I$. In this way, we express
each element of $I^q$ as a $\CC$-linear combination of products of
$q$ factors $\al_j-1$ plus a term in $I^N$, with $N>q$. 
To eliminate $I^N$ we consider the $I$-adic completion $\hat R$ of
$\CC[\tGm]$, with closure $\hat I^q$ of $I^q$. Each element of
$\hat I\supset I$ is a countable sum of products of a complex number
and finitely many factors $\al_j-1$. Since
$\hat I^{q+1}\backslash \hat I^q$ and $I^{q+1}\backslash I^q$ are
isomorphic, it suffices to prove that the $\bb(\ii)$ are linearly
independent as elements of $\hat I^{q+1}\backslash \hat I^q$.

We suppose that there are $x_\ii\in\CC$ for all $q$-tuples~$\ii$ such
that
\begin{equation}
\label{sumxiali}
\sum_\ii x_\ii (\al_{\ii(1)}-1) \cdots (\al_{\ii(q)}-1) \;\in \; \hat
I^{q+1}\,.\end{equation}
We can write this element of $\hat I^{q+1}$ as $\sum_\jj c_\jj\,
\x_\jj$ with $c_\jj\in \CC$, and $\x_\jj$ running over the countably
many products $ (\al_{\jj(1)}-\nobreak 1)\cdots
(\al_{\jj(m)}-\nobreak 1)$ with $m$-tuples from $\{1,\ldots,\ngen\}$
for all $m>q$.

We form the ring $N=\CC\langle\Xi_1,\ldots,\Xi_\ngen\rangle$ of power
series in the non-commuting, algebraically independent (over $\CC$)
variables $\Xi_1,\ldots,\Xi_t$, and the two-sided ideal $Z$ in~$N$
generated by the commutators
\[\Xi_j \,\Xi_\ngen-\Xi_\ngen\,\Xi_j\qquad \text{for }1\leq j \leq
\ngen\,.\]

The quotient ring $M:=N/Z$ is non-commutative if $\ngen\geq 3$. The
relations between the generators imply that there is a group
homomorphism $\ph : \tGm\rightarrow M^\ast$ given by
$\ph(\al_j)=1+\Xi_j$ for $1\leq j \leq \ngen$, and
\[ \ph(\e_j) \= (1+\Xi_\ngen)^{1/v_j} \= \sum_{l\geq 0} \binom{1/v_j}l
\, \Xi_\ngen^l\,. \]
This group homomorphism induces a ring homomorphism $\hat \ph:\hat R
\longrightarrow M$, for which \[ \hat \ph(\x_\ii) =
\hat\ph(\al_{\ii(1)}-1)\, \hat\ph(\al_{\ii(2)}-1)\cdots
\hat \ph(\al_{\ii(|\ii|)}-1)
 \=\Xi^\ii:= \Xi_{\ii(1)}\,\Xi_{\ii(2)}\cdots\Xi_{\ii(|\ii|)}\,.
\]
Now we have
\[ \sum_\ii x_\ii \Xi^\ii \= \hat \ph\biggl(\sum_\ii x_\ii\,
\xi_\ii\biggr)
\= \hat \ph \biggl( \sum_\jj c_\jj \x_\jj\biggr) \= \sum_\jj c_\jj
\Xi^\jj\,,\]
where $\ii$ runs over $q$-tuples, and $\jj $ runs over countably many
tuples with length strictly larger than~$q$. Hence all $x_\ii$ (and
$c_\jj$) vanish.
\end{proof}

So for $\tGm$ with cusps the trivial $\tGm$-module $I^{q+1}\backslash
I^q$ is always non-trivial. The dimension is equal to the number of
all $\tGm$-$q$-tuples. Thus we have
\begin{equation}
\label{formula}
\dim_\CC (I^{q+1}\backslash I^q)
\= n(\tGm,q) \= \sum_{m=0}^q
(\ngen\!-\!1)^m
\= \begin{cases}1&\text{ if }\ngen=1\,,\\
q+1&\text{ if }\ngen=2\,,\\
\frac{(\ngen-1)^{q+1}-1}{\ngen-2}&\text{ if }\ngen\geq 3\,.
\end{cases}
\end{equation}

We obtain for each $\tGm$-module~$V$ an exact sequence
\[ 0 \longrightarrow \ivt V q \longrightarrow \ivt V {q+1}
\stackrel{\ml_q}\longrightarrow \bigl( V^\tGm\bigr)^{n(\tGm,q)}\]
with
\begin{equation}
\label{tmldef}
\bigl(\ml_q f \bigr)_\ii \= f|(\al_{\ii(1)}-1)\cdots
(\al_{\ii(q)}-1)\,.
\end{equation}

For the modular group, we have $\npar=1$, $\nell=2$ and $g=0$, hence
$\ngeng{\Gmod}=1$, and $n(\tGmod,q)=1$ for all~$q$. So in contrast to
$\Gmod$, for $\tGmod$ we may hope for non-trivial higher order
automorphic forms.

\section{Maass forms with generalised weight on the universal covering
group}\label{sect-gw-uc}
\subsection{The logarithm of the Dedekind eta function} In the
introduction we mentioned that one of the motivating objects for the
study of higher order forms on the universal covering group is the
logarithm of the Dedekind eta function. Its branch is fixed by the
second of the following expressions:
\begin{equation} \label{FeLdef}
\log\eta(z) \= \frac{\pi i z}{12}+\sum_{n=1}^\infty
\log\bigl(1-e^{2\pi i n z}\bigr) \= \frac{\pi i z}{12} -
\sum_{n=1}^\infty \s_{-1}(n)\,e^{2\pi i n z} \,.
\end{equation}
where $\s_u(n)=\sum_{d\divides n} d^u$. One can show that its
behaviour under $\Gmod$ is given by
\begin{equation}\label{logeta-trf}
\log\eta(z+1)\= \log\eta(z) + \frac{\pi i}{12}\,,\quad\log\eta(-1/z)
\= \log\eta(z) + \frac12\log z-\frac{\pi i}4\,.
\end{equation}
Except for the term $\frac12\log z$ this looks like a second order
holomorphic modular form of weight zero. In the next few sections we
make this precise by generalizing the concept ``weight'' of Maass
forms, and replacing the group $\Gmod$ by the discrete subgroup
$\tGmod$ of the universal covering group of~$\SL_2(\RR)$, using the
notation we introduced in the last section.

We first define the following function on~$\uhp\times\RR$:
\begin{equation}
\label{Ldef}
L(z,\th) \= \frac12\log y + 2\log\eta(z) + i \th\,.
\end{equation}
With \eqref{logeta-trf} we check easily that
$L\bigl( \gm(z,\th)) \=L(z,\th) + i \al(\gm)$ for $\gm=t$ and
$\gm=s$, where $\al:\tGmod \rightarrow\frac\pi 6\ZZ$ is the group
homomorphism at the end of~\S\ref{sect-cdsg}. Thus $L$ has the
transformation behaviour of a second order invariant in the functions
on~$\tG$ for the action by left translation.

Routine computations show that $L$ satisfies $\EE^-L=0$, $\WW L = i$
and $\om L =\frac12$.

\subsection{General Maass forms on the universal covering group} The
considerations on the function $L$ on $\tG$ induced by the logarithm
of the eta functions lead us to the definition of Maass forms on
$\tG$.

We first establish appropriate notions of weight and holomorphicity.
We say that a function $f$ on $\tG$ has \emph{(strict) weight}~$r\in
\CC$ if $f(z,\th) = e^{ir\th}\,f(z,0)$. Such a function is completely
determined by the function $f_r(z)=f(z,0)$ on~$\uhp$ and satisfies
$\WW f =ir f$.

The left translation of $f$ by $\tilde g$, with $g=
\matr{a}{b}{c}{d}\in \SL_2(\ZZ)$, induces an action $|$ of $\tG$ on
the space of functions of strict weight on $\tG$. On the other hand,
$\tG$ acts on the space of corresponding functions $f_r$ on $\uhp$
via
\[ f_r| \tilde g (z) \= e^{-ir\arg(cz+d)}\, f_r\bigl(
\frac{az+b}{cz+d}\bigr)\,,\]
The latter action corresponds to~\eqref{actionb} when $r\in \ZZ$. In
general, this is an action of~$\tG$, not of $\SL_2(\RR)$. The map
$f \to f_r$ defined above on the space of functions of strict weight
is then equivariant in terms of these actions.

Many important functions on $\tG$, such as $L$, are not eigenfunctions
of the operator $\WW$, but they are annihilated by a power of $\WW$.
This suggests the following definition.
\begin{defn} An $f\in C^\infty(\tG)$ has \emph{generalised
weight}~$r\in \CC$ if $(\WW-ir)^n f=0$ for some~$n\in \NN$.
\end{defn}
Thus, $L$ and all its powers have generalised weight~$0$.

Next, holomorphy of $F_r=y^{-r/2}f_r$ corresponds to the property
$\EE^- f=0$.

\begin{defn} We call any differentiable function $f$ on~$\tG$
\emph{holomorphic} (resp. \emph{antiholomorphic}) if $\EE^- f=0$,
(resp. $\EE^+f=0$). We call any twice differentiable function $f$
on~$\tG$ \emph{harmonic} if it satisfies $\om f=0$.
\end{defn}
Note that, for functions of non-zero weight, this definition of
harmonicity does not correspond to the use of the word harmonic in
``harmonic weak Maass forms'' in, 
{\sl e.g.},~\cite{BOR}.

With these definitions we set
\begin{defn}\label{defTa} Let $k,\ld\in \CC$. Let $\tGm$ be a discrete
cofinite subgroup of~$\tG$. \\
i. The space $\tilde \A_k(\tGm,\ld)$ consists of the smooth functions
$f:\uhp\times\RR \rightarrow\CC$ that satisfy:
\begin{enumerate}
\item[a)] \emph{(Eigenfunction Casimir operator) } $\om f = \lambda
f$.
\item[b)] \emph{(Generalised weight) } $\bigl( \WW -\nobreak ik
\bigr)^n f=0$ for some $n\in \NN$.
\item[c)] \emph{(Exponential growth) } There exists $a\in \RR$ such
that for all compact sets $X$ and $\Th\subset\RR$ and for all cusps
$\k$ of $\tGm$ we have
\begin{equation} \label{expgrow}
f\bigl(\tilde g_\k (x+\nobreak iy,\th )
\bigr) = \oh(e^{ay})
\end{equation}
as $y\rightarrow\infty$ uniformly in $x\in X$ and $\th\in \Th$.
\end{enumerate}
\noindent ii.
$$\tilde E_k(\tGm ,\ld):=\tilde \A_k(\tGm,\ld)^{\tGm}$$
(where $\tGm$ acts by left translation). The elements of
$\tilde E_k(\tGm ,\ld)$ are called \it Maass forms on $\tG$ of
generalised weight $k$ and eigenvalue $\ld$ for~$\tGm$.
\end{defn}
The space $\tilde E_r(\tGm, \ld)$ is infinite dimensional. Further,
since $\om$ and $\WW$ commute with left translations in $\tG$, the
space $\tilde\A_k(\tGm,\ld)$ is invariant under left translation by
elements of~$\tGm$.

When $k\in 2\ZZ$, the space $E_k(\Gm,\ld)$ can be identified with
$\tilde E_k(\tGm ,\ld)$. We prove the following slightly stronger
statement.
\begin{thm}\label{thm-tAi}Let $\tGm$ be a cofinite discrete subgroup
of~$\tG$, and let $k,\ld\in \CC$. If $\tilde \A_k(\tGm,\ld)^\tZ$
contains a non-zero element~$f$, then $k\in 2\ZZ$ and $ \partial_\th
f\bigl( z, \th \bigr)=ik f\bigl( z,\th \bigr)$.

If $k\in 2\ZZ$, then the elements $f\in \tilde E_k(\tGm ,\ld)$
correspond bijectively to the Maass forms $F\in E_k(\Gm,\ld)$ by
$$f\bigl(z, \th \bigr)
= y^{k/2} \,F(z)\, e^{ik\th}.$$
\end{thm}
So the condition of $\tZ$-invariance implies that the weight $k$ is
even, and that the weight is \emph{strict}, {i.e.}, condition~b)
holds with $n=1$.
\begin{proof}[Proof of Theorem~\ref{thm-tAi}.] Any smooth function
$f\in C^\infty( \uhp\times\nobreak\RR)$ satisfying b) in
Definition~\ref{defTa} can be written in the form $f\bigl( z,\th
\bigr) = \sum_{j=0}^{n-1} \ph_j(z) \, e^{ik\th} \, \th^j$, with
$\ph_j \in C^\infty(\uhp)$.

If such a function is left-invariant under $\tZ$, then the action of
$k(\pi m)\in \tZ\subset\tGm$, implies for each~$m\in \ZZ$:
\[e^{\pi i k m} \sum_j \ph_j(z)e^{ik\th} (\th+\nobreak \pi m)^j =
\sum_j \ph_j(z) e^{ik\th} \th^j\qquad\text{ for all }m\in \ZZ\,.\]
With induction this gives $k\in 2\ZZ$ and $\ph_j=0$ for $j\geq 1$,
hence $f(z,\th) = \ph_0(z) e^{ik\th}$. Moreover, the stronger
condition $f\in \tilde E_k(\tGm,\ld) = \tilde \A_k(\tGm,\ld)^\tGm$ can be
checked to be equivalent to $F_k\in E_k(\Gm,\ld)$ for
$F_k(z) = y^{-k/2} f(z,0)$.
\end{proof}

We have the following generalisation of Theorem~\ref{thm-Mf-mp}.
\begin{thm}\label{thm-tMf-mp}Let $\tGm$ be a cofinite discrete
subgroup of~$\tG$ with cusps. Then the $\tGm$-module $\tilde
\A_k(\tGm,\ld)$ is maximally perturbable for each $k\in 2\ZZ$ and
each $\ld\in \CC$.
\end{thm}
In Section~\ref{sect-hoM} 
we will prove this theorem. In this section we will show that it
implies the corresponding result for $E_k(\Gm,\ld)$. We first give
some facts that are of more general interest.

The map identifying $E_k(\Gm, \lambda)$ and $\tilde E_k(\tGm ,\ld)$
can be extended to an isomorphism
$$\mu: \A_k(\Gm, \ld) \longrightarrow \tilde\A_k(\tGm,\ld)^{\tZ}.$$
Since the centre $\tZ$ of $\tGm$ acts trivially on
$\tilde\A_k(\tGm,\ld)^{\tZ}$, it can be considered as a $\Gm$-module.
With this interpretation we obtain an identification of the
$\Gm$-modules $\A_k(\Gm,\ld)$ and $\tilde\A_k(\tGm,\ld)^{\tZ}$.
Specifically, for $F\in \A_k(\Gm,\ld)$,
$g\in\tilde\A_k(\tGm,\ld)^\tZ$ we have
\begin{equation}
\begin{aligned}
(\mu f)(z,\th)&\= y^{k/2}\, F(z)\, e^{ik\th}\,,\\
(\mu^{-1}g)(z)&\= y^{-k/2}\, g(z,0)\,,\\
\mu(F|_k\gm)&\=\mu(F)|\nu(\gm)&&(\gm\in \Gm)\,,\\
\mu^{-1}(g|\tZ\dt)&\=\mu^{-1}(g)|_k \nu^{-1}(\tZ\dt)&&(\dt\in \tGm)\,,
\end{aligned}
\end{equation}
where $\nu$ denotes the isomorphism identifying $\Gm$ with
$\tZ\backslash\tGm$.

\begin{prop}\label{prop-reduction}Let $\Gm$ be a cofinite discrete
subgroup of~$G$ 
with cusps, and let $\tGm=\pr^{-1}\Gm$. If the $\tGm$-module $V$ is
maximally perturbable, then the subspace $V^{\tZ}$, considered as a
$\Gm$-module, is maximally perturbable.
\end{prop}
\begin{proof}The projection $\pr:\tGm\rightarrow \Gm$ induces linear
maps $\pr: \CC[\tGm]\rightarrow \CC[\Gm]$ between the group rings,
$\pr: I_\tGm\rightarrow I_\Gm$ between the augmentation ideals, and
$\pr: I_{\tGm}^{q+1}\backslash I_{\tGm}^q \rightarrow
I_{\Gm}^{q+1}\backslash I_{\Gm}^q$ for all $q\in \NN$. Since,
pr$(A_{\ii})=\alpha_{\ii}$, on the basis elements $\bb_{\tGm}(\ii)$
in Proposition~\ref{prop-augpow} and $\bb_\Gm(\ii)$
in~\eqref{bi-def0} we have for $\tGm$-$q$-tuples:
\begin{equation}
\pr\, \bb_\tGm(\ii) \= \begin{cases}
\bb_\Gm(\ii) &\text{ if }\ii(l)<\ngen \text{ for }l=1,\ldots,q\,,\\
0&\text{ if }\ii(q)=\ngen\,.
\end{cases}
\end{equation}This means that we have the commutative diagram
\[ \xymatrix{ 0 \ar[r] & V^{\tGm,q}_{} \ar[r] & V^{\tGm,q+1}_{}
\ar[r]^(.35){\ml_q} & \hom\bigl( I_\tGm^{q+1}\backslash I_\tGm^q,
V^\tGm\bigr) \ar[r]&0\\
0\ar[r] & (V^\tZ)^{\Gm,q}_{} \ar[r] & (V^\tZ)^{\Gm,q+1}_{}
\ar[r]^(.35){\ml_q} & \hom\bigl( I_\Gm^{q+1}\backslash I_\Gm^q,
(V^\tZ)^\Gm\bigr) \ar@{->}[u] } \]
where the vertical arrow sends $f:I_\Gm^{q+1}\backslash
I_\Gm^q\rightarrow
(V^\tZ)^\Gm=V^\tGm$ to $\tilde f:I_\tGm^{q+1}\backslash I_\tGm^q
\rightarrow V^\tGm$ such that $\tilde f\bigl( \bb_{\tGm}(\ii)\bigr) =
f\bigl( \bb_\Gm(\ii)\bigr)$ if $\ii\in \{1,\ldots,\ngen-1\}^q$, and
$\tilde f\bigl( \bb_{\tGm}(\ii)\bigr) = 0$ otherwise.

We want to write a given $f: I_\Gm^{q+1}\backslash I_\Gm^q
\rightarrow (V^\tZ)^\Gm$ as $\ml_q v_0$ with $v_0\in
(V^\tZ)^{\Gm,q+1}$. By assumption, there is an element
$v\in V^{\tGm,q+1}$ such that $\ml_q v = \tilde f$. If
$v|(\z-\nobreak 1)=0$, then $v\in V^{\tGm,q+1} \cap V^\tZ =
(V^\tZ)^{\Gm, q+1}$, and we are done.

Suppose that $w=v|(\z-\nobreak 1)\neq0$. Take $r\in [1,q]$ minimal
such that $w\in V^{\tGm,r}$. We will show that we can replace $v$ by
another element $v_1\in v+V^{\tGm,q}$ with $v_1|(\z-\nobreak 1)\in
V^{\tGm,r_1}$ and $r_1<r$. Repeating this process brings us
eventually to $v_j|(\z-\nobreak1)=0$. For this $v_j$ we will have
$\ml_q v_j = \tilde f$ and $v_j|(\z-\nobreak1)=0$ which, according to
the remark of the last paragraph suffices to prove the proposition.

From $w|(\gm_1-\nobreak 1)\cdots(\gm_{q-1}-\nobreak 1) =
v|(\gm_1-\nobreak 1)\cdots(\gm_{q-1}-\nobreak 1)(\z-\nobreak 1)=
\tilde f(\gm_1,\cdots,\gm_{q-1},\z)=0$ we conclude that $r\leq q-1$.
Define $\tilde g\in \hom(I_\tGm^{r+1}\backslash I_\tGm^r, V^\tGm)$ by
$\tilde g(\bb_\tGm(\jj)) = w|(\al_{\jj(1)}-\nobreak
1)\cdots(\al_{\jj(r-1)}-\nobreak
1)$ if the $\tGm$-$r$-tuple $\jj$ satisfies $\jj(r)=\ngen$ and $
\tilde g(\bb_\tGm(\jj)) =0$ otherwise. There is $u\in
V^{\tGm,r+1}\subset V^{\tGm,q}$ with $\ml_r u=\tilde g$. We take
$v_1=v-u\in v+V^{\tGm,q}$. We check that for all
$\tGm$-$(r\!-\!1)$-tuples~$\jj$
\begin{align*}
v_1|&(\z-1) ( \al_{\jj(1)}-1) \cdots (\al_{\jj(r-1)}-1)\\
& \= w | ( \al_{\jj(1)}-1) \cdots (\al_{\jj(r-1)}-1) - u | (
\al_{\jj(1)}-1) \cdots
(\al_{\jj(r-1)}-1)(\z-1)\\
&\=0\,.
\end{align*}
This shows that $v_1|(\z-\nobreak 1)$ has order less than~$r$.
\end{proof}

\noindent\it
Proof of Theorem~\ref{thm-Mf-mp}. \rm From Theorem~\ref{thm-tMf-mp},
$V=\tilde\A_k(\tGm,\ld_k)$ is maximally perturbable. Therefore, by
Proposition~\ref{prop-reduction}, the space
$\tilde\A_k(\tGm,\ld_k)^{\tilde Z} \cong \A_k(\Gm,\ld_k)$ is
maximally perturbable too. \qed

This proof illustrates the fact that, for groups with cusps, there are
really more higher order forms with generalised weight than with
strict weight: The basis in Proposition~\ref{prop-augpow} is for all
such discrete groups larger than the corresponding basis in
\S\ref{canon-gener}.

\subsection{Holomorphic forms on the universal covering group}
\label{holom}
\begin{defn}For $k\in 2\ZZ$ we define $\Hol_k(\tGm)$ as the space of
elements of $C^\infty(\uhp\times\RR)$ that satisfy
\begin{enumerate}
\item \emph{(Holomorphy)} \ $\EE^- f=0$.
\item \emph{(Generalised weight)} \ $(\WW-ik)^n f=0$ for some
$n\in \NN$.
\item \emph{(Exponential growth)} as described in condition~c) in
Definition~\ref{defTa}.
\end{enumerate}
\end{defn}
This is a $\tGm$-module for the action by left translation. We denote
by $\Hol^p_k(\tGm)$ (resp. $\Hol^c_k(\tGm)$) the space of $f \in 
\Hol_k(\tGm)$ satisfying $f\bigl(\tilde g_\k (x+\nobreak iy,\th )
\bigr) = \oh(y^C)$ for some $C\in\RR$
(resp. $f\bigl(\tilde g_\k (x+\nobreak iy,\th )
\bigr) = \oh(e^{ay})$ for some $a<0$) instead of
\eqref{expgrow}.

We will prove:
\begin{thm}\label{thm-thol-mp}Let $\tGm$ be a cofinite discrete
subgroup of~$\tG$ with cusps. Then the $\tGm$-module $\Hol_k(\tGm)$
is maximally perturbable for each $k\in 2\ZZ$.
\end{thm}

\noindent\it
Proof of Theorem~\ref{thm-hol-mp}. \rm As in the case of general Maass
forms, we can show that, for $k \in 2  \mathbb Z$,
$\A_k\hol(\Gm,\ld_k) \cong \Hol_k(\tGm)^\tZ$. Then,
Proposition~\ref{prop-reduction} implies Theorem~\ref{thm-hol-mp}.
\qed
\medskip

\noindent{\bf Second order forms and derivatives of
$L$-functions. } 
With this definition, $L$ is a second order invariant belonging to
$\Hol_0(\tGmod)^{\tGmod, 2}$. (Incidentally, this example shows that,
for generalised weight $k$, the space $\Hol_k(\tGm)$ need not be
contained in $\tilde \A_k(\tGm,\ld_k)$.)

Based on $L$ we can construct a second-order form which is related to
derivatives of classical modular forms. Specifically, for positive
integer $N$, denote by $G_N$ the group generated by $\tilde g$,
$g \in <\Gm_0(N), W_N>$ where
$W_N:=\matr{0}{-\sqrt N^{-1}}{\sqrt N}{0}.$ Set
$$L_1(z, \th)=L(z, \th)+L(Nz, \th).$$ Using the transformation law for
$L$ and the identity $\matr{N}{0}{0}{1} \matr{a}{b}{Nc}{d}=
\matr{a}{Nb}{c}{d} \matr{N}{0}{0}{1}$, a routine calculation implies
that, for some $\beta \in$ Hom$(G_N, \CC)$,
$$L_1\bigl( \gm(z,\th)) =L_1(z,\th) + i \beta(\gm), \qquad
\text{for all} \, \gm \in G_N.$$
Let now $f$ be a newform in the space $S_{\!2}$ of cusp forms of
weight $2$ for $\Gm_0(N)$ such that its $L$-function $L_f(s)$
vanishes at $1$. Then, $f(W_N w)d(W_N w)=f(w)dw$ and, for all
$\th \in \RR$,
\begin{equation} \label{Lfunction}
\begin{aligned}
\int_0^{\infty}f(iy)L_1(iy, \th)diy=
-\int_{W_N 0}^{W_N \infty}f(iy)L_1(iy, \th)diy&=
-\int_0^{\infty}f(W_Niy)L_1(W_N iy, \th)d(W_Niy)\\
&=-\int_0^{\infty}f(iy)L_1(W_N iy, \th)diy.
\end{aligned}
\end{equation}
Since $L_1(z, \th+x)=L_1(z, \th)+2ix$ and $L_f(1)=2 \pi
\int_0^{\infty}  f(iy)dy=0$, our integral is independent of $\th$. It
further equals
\begin{equation}
-\int_0^{\infty}f(iy)L_1(\tilde{W_N}(iy, 0))diy=
-\int_0^{\infty}f(iy)(L_1(iy, 0)+i\beta(\tilde{W_N}))diy=
-\int_0^{\infty}f(iy)L_1(iy, 0)diy
\end{equation}
Therefore, $\int_0^{\infty}f(iy)L_1(iy, 0)dy=
-\int_0^{\infty}f(iy)L_1(iy, 0)dy$, i.e.
$$\int_0^{\infty}f(iy)L_1(iy, 0)dy=0 \quad \text{and hence}$$
$$\int_0^{\infty}f(iy)\log y\,dy+2\int_0^{\infty}f(iy)u(iy)\,dy=0$$
where $u(z):=\log(\eta(z))+\log(\eta(Nz)).$ From this we see that,
since, $L'_f(s)= 2 \pi \int_0^{\infty} f(iy)\log(y)dy$, we can
retrieve, from a alternative perspective, the formula
$$L'_f(1)=-4 \pi \int_0^{\infty} f(iy) u(iy)dy$$
first derived in \cite{G}.

Thus, Goldfeld's expression of $L_f'(1)$ is equivalent to the
orthogonality of $L_1 \in  \Hol^p_0(G_N)^{G_N, 2}$ to
$S_{\!2} \hookrightarrow  \Hol^c_2(G_N)^{G_N}$ in terms of the
pairing
$$\langle\cdot, \cdot\rangle: \Hol^c_2(G_N)^{G_N} \times \Hol^p_0(G_N)^{G_N, 2} 
\to \CC$$
defined by
$$\langle g, h\rangle=\int_0^{\infty} g(iy, 0)h(iy, 0) \frac{dy}{y}.$$

\subsection{Examples of higher order forms for the full modular
group}\label{sect-ex2}
Theorems \ref{thm-tMf-mp} and~\ref{thm-thol-mp} show that there are
perturbations of~$1$ for the full original~$\tGmod$ of $\SL_2(\ZZ)$
in the universal covering group. Since $\ngeng{\Gmod}=1$ all these
perturbations are commutative (see \eqref{formula}).
\smallskip

\noindent
$1.$ The function $L$ can lead to second order \emph{harmonic
perturbations of~$1$}. Specifically, although
$L\not\in \ivt{\tilde \A_0(0)}{2}$ (because
$\om \,L = \frac12$), the imaginary part
$\im L:(z,\th) \mapsto 2\,\im \log\eta(z)+\th$ is harmonic, has
second order, and corresponds to the linear form
$\al\in \mult^1(\tGmod,\CC)$. It has generalised weight~$0$, and it
is not holomorphic.
\smallskip

\noindent
$2$. Set $\ch_r = e^{ir\al}$, $r\in \CC$, where $\al\in
\hom(\tGmod,\CC)$ is given by $\al\bigl(n(1)\bigr)=\frac\pi 6$ and
$\al\bigl(k(\pi/2)\bigr)=\frac\pi 2$. The family
\begin{equation}
r\mapsto e^{rL(z,\th)} \= y^{r/2} \, \eta(z)^{2r} \, e^{ir\th}
\end{equation}
consists of elements of $\Hol_r(\tGm)$ that are $\tGmod$-invariant
under the action given by
$$(f|\gm)(z)=f(\gm z) \overline{\chi_r(\gm)}.$$
By Proposition~\ref{prop-diff}, for $k\geq 1$ the derivative
\[ \partial_r^k e^{rL(z,\th)}\bigr|_{r=0}= L(z,\th)^k \]
is a holomorphic perturbation of~$1$ of order~$k+1$. The corresponding
element of $\mult^k(\tGmod,\CC)$ is $i^k\, k! \, \al ^{\otimes
k}$.\smallskip

\noindent
$3$. It is possible to obtain a more or less explicit description of a
harmonic perturbation of~$1$ of order~$3$. We sketch how this can be
done with the meromorphic continuation of the Eisenstein in weight
and spectral parameter jointly. This family is studied in
\cite{Br86}. In that work, automorphic forms are described as
functions on~$\uhp$ transforming according to a multiplier system
of~$\Gmod$. These correspond to functions on~$\tG$ that transform
according to a character of~$\tGmod$. Carrying out the reformulation,
we can rephrase \S2.18 in~\cite{Br86} as stating that there is a
meromorphic family of Maass forms on $U\times\CC$, where $U$ is some
neighborhood of $(-12,12)$ in~$\CC$. We retrieve the exact family
studied in \cite{Br86} by considering $z\mapsto E(r,s;z,0)$. For each
$(r,s)\in U\times\CC$ at which $E$ is not singular it is an
automorphic form of weight~$r$ for the character $\ch_r=e^{ir\al}$ of
$\tGmod$ with eigenvalue $\ld_s=\frac14-s^2$. It is a meromorphic
family of automorphic forms on $\tGmod$ with character $\ch_r$ with a
Fourier expansion of the form
\begin{equation} E(r,s) \=\mu_r(r/12,s) + C_0(r,s) \,\mu_r(r/12,-s) +
\sum_{n\neq 0} C_n(r,s)\,\om_r(n+r/12,s)\,,\end{equation}
where the $C_n(r,s)$ are meromorphic functions, and where we use the
following notations.
\begin{equation}\label{om-mu-expl}
\begin{aligned}
\om_r(\nu,s;z,\th)&\= e^{2\pi i \nu x}\,
W_{r\sign(\re\nu)/2,s}(4\pi\nu\sign(\re\nu)y)\, e^{ir\th}\,,\\
\mu_r(\nu,s;z,\th)&\= e^{2\pi i \nu z}\, y^{\frac12+s}\,
\hypg11\bigl(\txtfrac12+s-\txtfrac r2;1+2s;4\pi\nu y\bigr) \,
e^{ir\th}\,.
\end{aligned}
\end{equation}
This family and its Fourier coefficient $C_0$ satisfy the following
functional equations.
\begin{equation}
\label{f.e.}
\begin{aligned}
E(r,-s)&\;=\; C_0(r,-s)E(r,s)\,,\\
E\bigl(r,s;-x+iy,-\th\bigr) &\;=\; E\bigl(-r,s;x+iy,\th\bigr)\,.
\end{aligned}
\end{equation}
Further, the restriction of this family to the (complex) line $r=0$
exists, and gives a meromorphic family of automorphic forms depending
on one parameter~$s$. This is a family of weight zero, so it does not
depend on the parameter $\th$ on~$\tG$. The resulting family
on~$\uhp$ is the meromorphic continuation of the Eisenstein series
for~$\Gmod$ in weight~$0$, with Fourier expansion
\begin{equation}\label{Four-E}
\begin{aligned}
E(0,s) &\= \mu_0(0,s) + \frac{\sqrt \pi\, \Gf(s)\,
\z(2s)}{\Gf(s+\frac12)\,\z(2s+1)}\, \mu_0(0,-s)\\
&\qquad\hbox{} + \frac{\pi^{s+\frac12}}{\Gf(s+\frac12)\,\z(2s+1)}
\sum_{n\neq0} \frac{\s_{2s}(|n|)}{|n|^{s+\frac12}}\, \om_0(n,s)\,.
\end{aligned}
\end{equation}
where
\begin{align*}
\mu_0(0,s;z, \th)&\=y^{\frac12+s}\,,\\
\om_0(n,s;z, \th) &\= e^{2\pi i n x} \, W_{0,s}(4\pi |n|y) \= e^{2\pi
inx}\, 2|n|^{1/2}\, K_s(2\pi|n|y)\,.
\end{align*}

At $(0,-\frac12)$ the family $E$ is holomorphic in both variables $r$
and $s$, with a constant as its value at $(0,-\frac12)$.
(This is a consequence of Proposition~6.5~ii) in~\cite{Br86}.) So in
principle, we obtain higher order harmonic perturbations of~$1$ by
differentiating $r\mapsto E(r,-\frac12)$. Here we encounter the
problem that we have an explicit Fourier expansion~\eqref{Four-E}
only for $E(0,s)$ and thus we cannot describe the derivatives in the
direction of~$r$ directly. To overcome this problem we use the
fact that for $r$ near~$0$ we have
\begin{equation}\label{E-spec}
\begin{aligned}
E\bigl( r,-\frac{1-r}2 ;z, \th \bigr)&\= H_r(z,\th) \= e^{rL(z,\th)}
\,,\\
E\bigl( r,-\frac{1+r}2;z,\th) &\= H_{-r}(-\bar
z,-\th)\=e^{-r\overline{L(z,\th)}}\,.
\end{aligned}
\end{equation}
The proof of the first equality is contained in 6.10 in~\cite{Br86}.
The second one follows from the second functional equation
in~\eqref{f.e.}. Now we use the Taylor expansion of $E$ of degree~$2$
at $(r,s)=(0,-\frac12)$:
\begin{equation}\label{Eis-Tayl}
\begin{aligned}
E(r,s) &\= 1 + r\; A_{1,0}+\bigl(s+\frac12\bigr)\; A_{0,1} \\
&\qquad\hbox{} + \frac12 r^2\; A_{2,0}+ r\,\bigl(s+\frac12\bigr)\;
A_{1,1}+ \frac12 \bigl(s+\frac12\bigr)^2\; A_{0,2} + \cdots
\end{aligned}
\end{equation}
By Proposition~\ref{prop-diff}, the coefficients $A_{1,0}$ and
$A_{2,0}$ are harmonic perturbations of~$1$ of order $2$ and $3$,
respectively. From~\eqref{E-spec} we obtain the following results:
\begin{equation}
\begin{aligned}
A_{1,0}&\= i\, \im L \,, &
\quad A_{0,1}&\= 2\re L\,,\\
A_{2,0}+\frac14 A_{0,2}&\= \,\re L^2\,,& A_{1,1}&\= i \, \im L^2\,.
\end{aligned}
\end{equation}
This confirms that $\im L$ is a second order harmonic perturbation
of~$1$. Differentiation in the direction of $s$ preserves
$\tGmod$-invariance. So $A_{0,1}=2\re L$ and $A_{0, 2}$ are
$\tGmod$-invariant. However these functions are not in the kernel
of~$\om$.

Thanks to the identity $A_{2,0}+\frac14 A_{0,2}=\re L^2$, to determine
the third order harmonic perturbation $A_{2, 0}$ it suffices to
explicitly compute $A_{0, 2}$ because $\re L^2$ is known in a fairly
explicit way. The function $A_{0, 2}$ can be obtained as the
coefficient of $\frac12
(s+\frac12)^2$ in the Taylor expansion of $E(0, s)$ at
$s=-\frac{1}{2}$. As a by-product of this computation we will also
obtain the $\tGmod$-invariant function $A_{0, 1}$ as the coefficient
of $s+\frac12$ in the same expansion. We shall examine each term of
the expansion separately.

Set $\x:=s+\frac12$. The first term of our expansion is
\begin{equation}
\label{c0a}
\mu_0(0,s;z,0) \= y^{\frac12+s} \= 1 + \x \log y + \x^2 \frac12(\log
y)^2+\cdots
\end{equation}

For the next term $\frac{\Ld(2s)}{\Ld(2s+1)}\,
 \mu(0,-s;z,0)= \frac{\Ld(2-2\x)}{\Ld(1-2\x)}\, \mu(0,-s;z,0)$ with
$\Ld(u) = \pi^{-u/2}\,\Gf\bigl( \frac u2 \bigr)\, \z(u)=\Ld(1-u)$, we
define $a_0$ and $b_1$ by
\begin{equation}
\Ld(1+h) \= h^{-1}+a_0+\cdots\,,\qquad \Ld(2+h)\= \frac\pi 6+b_1 h +
\cdots\,.
\end{equation}
We get
\begin{equation}\label{c0b}
\frac{\Ld(2s)}{\Ld(2s+1)}\mu_0(0,-s;z,0) \= -\frac\pi3\, y\, \x +
\bigl( 4b_1 - \frac{2\pi a_0}3 + \frac\pi 3 \, \log y \bigr)\,y\,
\x^2+\cdots\,.
\end{equation}

For the other terms we use
\begin{align*}
W_{0,-s}(t)&\= W_{0,s}\=
\frac{e^{-\frac{t}{2}}}{\Gf(\frac12+s)}\int_0^\infty e^{-x} \bigl(
x(1+\frac{x}{t}) \bigr)^{s-\frac12}\, dx\,,\displaybreak[0]\\
W_{0,1/2}(t)&\= \frac{e^{-\frac{t}{2}}}1 \cdot 1 \=
e^{-\frac{t}{2}}\,,\displaybreak[0]\\
-\partial_s W_{0,s}(t) \bigr|_{s=-\frac12}&\=
\partial_s W_{0,s}(t) \bigr|_{s=\frac12}\=
-\frac{e^{-\frac{t}{2}}}{1^2}\Gf'(1)\cdot 1 + e^{-\frac{t}{2}}
\int_0^\infty e^{-x}\, \log\bigl(x(1+\frac{x}{t}) \bigr)\, dx\\
&\= e^{-\frac{t}{2}} \biggl( -\Gf'(1) + \Gf'(1) + \int_0^\infty e^{-x}
\,\log(1+\frac{x}{t})\, dx \biggr)\\
\text{(part.\ int.)}\quad&\= e^{-\frac{t}{2}} \int_0^\infty e^{-x}
\,\frac{dx}{x+t} \= e^{\frac{t}{2}} \int_t^\infty e^{-x}\,\frac
{dx}x= e^{\frac{t}{2}}\, \Gf(0,t)\,,
\end{align*}
with the incomplete gamma-function $\Gf(a,t) = \int_t^\infty e^{-x}\,
x^{a-1}\, dx$. With these ingredients:
\begin{equation}\label{cn}
\begin{aligned}
&\frac{\s_{2s}(|n|)}{\Ld(2s+1)\, |n|^{s+\frac12}} \om_0(n,s;z,0) \=
\sum_{d\bigm| |n|} \frac1d \biggl( -2 e^{-2\pi |n|y}\,\xi \\
&\quad\hbox{}
+ \bigl( 2 e^{2\pi|n|y}\Gf(0,4\pi|n|y)- 2e^{-2\pi|n|y} \log
\frac{d^2}{|n|}-4a_0\, e^{-2\pi|n|y }\bigr)\,\x^2+\ldots \biggr)\,
e^{2\pi i n x}\,.
\end{aligned}
\end{equation}

The results in \eqref{c0a}, \eqref{c0b}\, and \eqref{cn} confirm that
the constant term equals~$1$, and that
\begin{align*}A_{0,1}(z,0)= \log y - \frac\pi 3 y -
2\sum_{n\geq1}\sum_{d\mid n} \frac1d \bigl( q^n + \bar q^n \bigr)
= 2\re \bigl( \frac12\log y + \frac{\pi i}6 z -
\sum_{n=1}^\infty\s_{-1}(n)\, q^n \bigr) \= 2\re L(z,0)\,,
\end{align*}
with the notation $q=e^{2\pi i z}$. The term of order $2$ leads to:
\begin{equation}
\begin{aligned}
&A_{0,2}(z,0)\= (\log y)^2 + \bigl( 8b_1-\frac{4\pi \, a_0}3 +
\frac{2\pi}3\log y\bigr)\, y \\
&\quad\hbox{}
+ \sum_{n=1}^\infty \biggl( -4a_0 \, \s_{-1}(n) (q^n+\bar q^n)
+ 2\s_{-1}(n)\,(q^{-n}+\bar q^{-n})\,\Gf(0,4\pi n y)
- 2 (q^n+\bar q^n) \sum_{d\mid n} \frac{\log(d^2/n)}d \biggr)\,,
\end{aligned}
\end{equation}
which is a complicated, but explicit expression.

A remarkable aspect of this computation that we have used an explicit
computation of the derivatives of the Eisenstein series in weight
zero to compute the second derivative in the $r$-direction of the
more complicated Eisenstein family in two variables. The basic
observation is \eqref{E-spec}, which shows that the Eisenstein family
has easy derivatives in two directions. The Taylor expansion of $E$
at $\bigl(0,-\frac12\bigr)$ has three monomials in order~$2$. So it
suffices to compute a second order derivative in one more direction
to get hold of all terms. Higher order terms in the Taylor expansion
have too many monomials for this method to work. We do not know how
to compute all harmonic perturbations of~$1$ of higher order.

\section{Higher order Fourier expansions}\label{sect-Ft}
This section is needed for the constructions on which the proofs of
Theorems \ref{thm-tMf-mp} and~\ref{thm-thol-mp} are based, but it is
also of independent interest. It provides a higher-order analogue of
the classical Fourier expansions.
\subsection{Fourier expansion of Maass forms}\label{sect-Fe}
If $f$ is in $\tilde E_r(\tGm, \lambda)$, then for each cusp~$\k$ of
$\Gm$ there is a Fourier expansion
\begin{equation}
\label{Ftdef}
f(\tilde g_\k g)\= \sum_\nu \Ft{\k,\nu} f(g)\,,
\qquad F_{\k,\nu}f(g) \= \int_0^1 e^{-2\pi i \nu x}\, f\bigl(\tilde
g_\k n(x)g\bigr)\, dx\,,
\end{equation}
 where $\nu$ runs through a class in $\CC\bmod\ZZ$ determined by~$\ch$
and the cusp~$\k$. The function $F_\nu f$ satisfies
$\Ft{\k,\nu} f(z,\th) = e^{2\pi i\nu
x}\, \Ft {\k,\nu} f(iy ,0)\, e^{ir\th}$ and $\om \Ft{\k,\nu} f =
\ld\; \Ft{\k,\nu} f$.

For each given $\nu$, $r$ and $s$ set
\begin{equation} \label{W_r}
\W_r(\nu,s)\,:=\,\{f: \tG \to \CC\;;\; \om f=(\frac{1}{4}-s^2)f,\;
f(z, \theta)=e^{2 \pi i \nu x+i r \theta}f(iy, 0)\}.
\end{equation}
Because of the second relation in the definition, $f \in \W_r(\nu, s)$
can be thought of as a function of~$y$. 
Therefore the space $\W_r(\nu,s)$ is isomorphic to the space of
$f: \RR \to \CC$ satisfying
\begin{equation}
\label{diff-eq-Wh}
-y^2 h''(y) + \bigl( 4\pi^2\nu^2 y^2-2\pi \nu r \,
y-\txtfrac14+s^2\bigr)\, h(y)\=0.
\end{equation}
It is convenient to write $\ld=\ld_s=\frac14-s^2$ with $s\in \CC$. We
can choose a fixed $s$ with $\re s\geq 0$ corresponding to the
eigenvalue $\ld=\ld_s$ under consideration. The spaces $\W_r(\nu,s)$
are two-dimensional. We will use the basis elements in \S4.2 of
\cite{Br94}.

$\bullet$ For $\re\nu\neq0$ a basis of $\W_r(\nu,s)$ is formed by
\begin{equation}\label{om-def-old}
\begin{aligned}
\om_r\bigl(\nu,s;z,\th\bigr)&\= e^{2\pi i \nu x} \, W_{r \sign(\re\nu
)/2,s}(4\pi \nu \sign(\re \nu)
y) \, e^{ir\th}\,,\\
\hat\om_r\bigl(\nu,s;z,\th\bigr)&\= e^{2\pi i \nu x} \, W_{-r
\sign(\re \nu)/2,s}(-4\pi \nu \sign(\re \nu)
y) \, e^{ir\th}\,.
\end{aligned}
\end{equation}
Here $W_{\mu,s}(t)$ is the Whittaker function that decreases
exponentially as $t\rightarrow\infty$. We use the branch of
$W_{\k,s}(z)$ that is holomorphic for $-\frac\pi 2<\arg z<
\frac{3\pi}2$. The asymptotic behaviour as $y \to \infty$, by \S4.2.1
in~\cite{Sl60} is:
\begin{equation}
\label{om-asympt}
\om_r(\nu,s;z,\th) \;\sim\; (4\pi \nu \e y)^{r\e/2}\, e^{2\pi
\nu(ix-\e y)
+ir\th},
\end{equation}
\begin{equation}\label{hat-om-asympt}
\hat\om_r(\nu, s;z,\th)\;\sim\; e^{-\pi i r\e/2}\, (4\pi \e\nu
y)^{-r\e/2}\, e^{2\pi \nu (ix+\e y)+ir\th}\,,
\end{equation}
where $\e$ denotes $\sign (\re \nu)$. The subspace of $\W_r(\nu,s)$
generated by $\omega_r(\nu, s)$ is denoted by $\W_r^0(\nu,s)$.

$\bullet$ For $\nu=0$, a basis is given by $\{y^{\frac12+s} e^{i r
\th}$, $y^{\frac12-s} e^{i r \th} \}$ if $s \ne 0$ and
 $\{y^{\frac12}e^{i r \th}$, $y^{\frac12}\log y e^{i r \th} \}$ if
$s=0$.\medskip

The following proposition characterises functions with exponential
growth in terms of Fourier series.
\begin{prop}\label{prop-Fe} Let $k\in 2\ZZ$, $\re s\geq 0$. Suppose
that the function 
$f\in C^\infty(\tGm\backslash\tG)$ satisfies $\om f= \ld_s f$ and
$\WW f=i k f$. Then it has at each cusp $\k$ an absolutely converging
Fourier expansion
\begin{equation}\label{Fe}
f\bigl( \tilde g_\k g \bigr) \= \sum_{n\in \ZZ} \Ft{\k,n}f(g)
\end{equation}
with $\Ft{\k,n} f \in \W_k(n,s)$. Moreover, $f\in
\tilde E_k(\tGm,\ld_s)$ if and only if there exists 
$N>0$ such that all Fourier terms $\Ft{\k,n}f $ with $|n|\geq N$ are
in $\W_k^0(n,s)$ for all cusps~$\k$.
\end{prop}
\begin{proof} The existence of such a Fourier expansion is a standard
result. A detailed proof in a more general setting can be found in
\cite{Br94}, \S4.1--3.

Fourier terms of automorphic forms inherit the growth behaviour of the
automorphic form. So if $f$ is in $\A_k(\tGm, \ld_s)$, all Fourier
terms satisfy $F_{\k,\nu}f(z,\th) = \oh(e^{ay})$ as
$y\rightarrow\infty$ for some $a$ depending on~$f$. Each Fourier term
of non-zero order is a linear combination of $\om_k(n,s)$ and
$\hat\om_k(n,s)$. From \eqref{hat-om-asympt} we conclude that
$F_{\k,n} f$ is a multiple of $\om_k(n,s)$ for all but finitely
many~$n$.

Conversely, suppose that for the cusp $\k$ we have $F_{\k,n} f = c_n\,
\om_k(n,s)$ for all $n$ with $|n|\geq N$. Then \eqref{om-asympt} and
the convergence of the Fourier expansion at $(z,\th)
= (iy_0,0)$ with $y_0>0$ implies that $c_n = \oh\bigl(
y_{0}^{-k\sign(n)/2} e^{2\pi |n| y_0} \bigr)$. This in turn shows
that the sum over $|n|\geq N$ gives a bounded contribution
in~\eqref{Fe} for all $y$ large enough. The terms with $|n|<N$ cannot
give a growth at the cusp~$\k$ larger than $\oh\bigl(y^a\,
e^{2\pi(N-1)y}\bigr)$ for some $a>0$.
\end{proof}
\begin{rmk}\label{rem-Fe}The $\tGm$-invariance in
Proposition~\ref{prop-Fe} is not necessary. Invariance under only the
parabolic elements of $\tGm$ suffices. If we work with functions~$f$
on $\bigl\{(z,\th)\;:\; y\geq y_0 \bigr\}$ for some $y_0>0$ that
satisfy $\om f = \ld_s f$, $\WW f =ikf$ and are left-invariant under
$\bigl\{ n(l)\;:\; l\in \ZZ\bigr\}$, then there is an expansion like
in~\eqref{Fe} on the set $y\geq y_0$, and exponential growth of such
a function is equivalent to the statement that all Fourier terms of
sufficiently large order are in $\W^0_k(n,s)$.
\end{rmk}

\subsection{Higher order Fourier terms}\label{sect-hoFt}
The higher order invariants of $\V_k(n, s)$ that we will define now
are the higher-order analogues of the classical Fourier terms.
\begin{defn} \label{V_k}
Let $k\in 2\ZZ$, $n\in \ZZ$, and $s\in \CC$. By $\V_k(n, s)$ we denote
the space of functions $f$ on $\tG$ that satisfy $\om f=\ld_s f$,
have generalised weight~$k$, and satisfy $\bigl(
\partial_x-\nobreak 2\pi i n \bigr)^m f=0$ for some $m\in \NN$ (which
may depend on~$f$).

For $n\neq 0$ we denote by $\V_k^0(n,s)$ the subspace of $f\in
\V_k(n,s)$ that satisfy $f(z,\th) = \oh(y^a\,e^{-2\pi|n|y}\bigr)$ as
 $y\rightarrow\infty$ for some $a\in \RR$.
\end{defn}

The free commutative group $\tDt$ generated by $\tau=n(1)$ and
$\z=k(\pi)$ acts on these spaces by left translation.
\begin{prop}Let $k,n,s$ be as above. The $\tilde \Delta$-modules
$\V_k(n, s)$ and $\V_k^0(n, s)$ are maximally perturbable.

For each $q\in \NN$ the elements $f\in \ivv{\V_k(n,s)}\tDt q$ satisfy,
for each $\dt>0$,
\begin{equation}\label{Fest} f(z,\th) \;\ll_\dt\; e^{(2\pi
|n|+\dt)y}\qquad(y\rightarrow\infty)
\end{equation}
uniformly for $x$ and $\th$ in compact sets. If $n\neq 0$ then for
each $q\in \NN$ the elements $f\in \ivv{\V_k^0(n,s)}\tDt q$ satisfy,
for each $\dt>0$,
\begin{equation}\label{F0est} f(z,\th) \;\ll_\dt \; e^{(\dt-2\pi
|n|)y}\qquad(y\rightarrow\infty)
\end{equation}
uniformly for $x$ and $\th$ in compact sets.
\end{prop}
\begin{proof} To prove that $\V_k(n,s)$ is maximally perturbable, we
start with a characterisation of the space $\V_k(n,s)^\tDt$. We first
note that $\W_k(n,s)\subset\V_k(n,s)^\tDt$. Conversely, if
$f\in \V_k(n,s)^\tDt$, then the reasoning in the proof of
Theorem~\ref{thm-tAi} shows that the weight of~$f$ is strict, and
also that $\partial_x f =
2\pi i n f$, hence $f(z,\th) = e^{2\pi i n x} f(iy,\th)$. So $f\in
\W_k(n,s)$. If, for $n\neq 0$, the function $f$ is also exponentially
decreasing it has to be a multiple of $\om_k(n,s)$. Therefore,
$\V_k^0(n,s)^\tDt=\W_k^0(n, s).$

Let $f$ be an arbitrary element of $\W_k(n, s).$ Since each of the
basis elements of $\W_k(n, s)$ is a specialisation of a holomorphic
family of elements of $\W_r(\nu,s)$, there is a holomorphic family of
$h(r, \nu) \in \W_r(\nu, s)$ such that $h(k,n)=f$. We have $h
\bigl(r,\nu;n(\x)k(\ell\pi)(z,\th) \bigr) \= e^{2\pi i\nu\x+ \pi ir\ell}\,
h(r,\nu;z,\th)$ for $\x\in \RR$ and $\ell\in\ZZ$.

Next consider the polynomials $Q_q \in \QQ[X]$ of degree $q$ defined
by
\begin{equation}\label{Qndef}
Q_0\=1\,,\quad Q_{q+1}(X+1)-Q_{q+1}(X)\=Q_q(X) \text{ and }
Q_q(0)\=0\text{ for }q\geq 1\,.
\end{equation}
Then for each $\mm=(m_1,m_2)$, $m_j\geq 0$ set
\begin{equation}\label{ztmm}
h^\mm_k(n,s) \= Q_{m_1}\bigl( \txtfrac1{\pi i} \partial_r\bigr)\,
Q_{m_2}\bigl( \txtfrac1{2\pi i} \partial_\nu \bigr)\, h(r,\nu)
\bigr|_{\nu=n,\, r=k}\,.
\end{equation}
Upon applying the differential operator $ \frac{1}{2\pi
i}\partial_\nu^a$ on
$h(r,\nu) |(\tau-1)\= \bigl( e^{2\pi i \nu}-1\bigr)
h(r,\nu)$ we obtain
\begin{multline}
(2\pi i)^{-a} \frac{\partial ^a h (r, \nu)}{\partial \nu^a}
\bigm|(\tau-1)\= \sum_{b=0}^{a-1} \binom {\;a\;} b \, (2\pi i)^{-b}
\frac{\partial^b h (r,\nu)}{\partial \nu^b}\=
\Bigl ( \bigl(\txtfrac{1}{2 \pi i} \partial_{\nu}+1\bigr)^a-
\bigl(\txtfrac{1}{2 \pi i} \partial_{\nu}\bigr)^a \Bigr )
h(r,\nu) \,.
\end{multline}
Therefore, 
\begin{multline}
Q_{m_2}\Bigl( \txtfrac{1}{2 \pi i} \partial_{\nu}\Bigr) h(r,\nu)
|(\tau-1)=\Bigl( Q_{m_2}\bigl( \txtfrac{1}{2 \pi i} \partial_{\nu}+1
\bigr )-Q_{m_2}\bigl( \txtfrac{1}{2 \pi i} \partial_{\nu}\bigr ) \Bigr
)
h(r,\nu)
=Q_{m_2-1}\bigl( \txtfrac{1}{2 \pi i} \partial_{\nu} \bigr )
h(r,\nu)\,.
\end{multline}
Since $\tau, \zeta$ commute, this implies $h^\mm_k(n,
s)|(\tau-\nobreak 1) = h^{(m_1,m_2-1)}_k(n, s)$. Likewise, we obtain
the transformation law $h^\mm_k(n, s)|(\z-\nobreak 1)$
$= h^{(m_1-1,m_2)}_k(n, s)$. Therefore, for $l_1+l_2=m_1+m_2$
($l_1, l_2 \ge 0$),
\begin{equation}
\label{typem}
 h^{(m_1, m_2)}_k(n, s) |(\z-1)^ {l_1}\, (\tau-1)^{l_2}=
\dt_{m_1,l_1}\dt_{m_2,l_2} f\,,
\end{equation}
thus obtaining the maximal perturbability of $\V_k(n,s)$. For
convenience, we shall call perturbations statisfying the
transformation law \eqref{typem} \emph{perturbations of type~$\mm$}.

Based on $\V_k^0(n,s)^\tDt=\W_k^0(n, s)$, we deduce in an analogous
way the maximal perturbability of $\V_k^0(n,s)$.

To prove \eqref{Fest} and \eqref{F0est}, we first note that the
maximal perturbability we have just shown implies that the functions
$h^{\mm}$ constructed from $f$'s ranging over a basis of
$\W_k(n,  s)$ (resp. $\W_k^0(n, s)$) induce a basis of the quotients
$\V^{\tilde  \Delta, q+1}/ \V^{\tilde \Delta, q}$. Therefore, it
suffices to show \eqref{Fest} and \eqref{F0est} for $h^{\mm}$ only.
In the case $n\neq0$, the family $h$ may be taken to be
$\om_r(\nu,s)$ or $\hat \om_r(\nu,s)$ in~\eqref{om-def-old}. For
these functions the question reduces to the asymptotic behaviour of
$\partial_t^{j}\partial_\k^{l} W_{\k,s}(t)$, since the factors
$e^{2\pi i\nu x}$ and $e^{ir\th}$ produce polynomials in $x$ and
$\th$, which yield constants when they vary through compact sets. The
differentiation of $4\pi\sign(\re\nu)\,\nu\, y$ yields only a power
of~$y$, which can be absorbed by the factor $e^{\dt y}$.

Differentiation of $W_{\k,s}(t)$ with respect to~$t$ does not change
the exponential part of the asymptotic behaviour, since derivatives
of $W_{\k,s}(t)$ are linear combinations of $W_{\k,s}(t)$ and
$W_{\k+1,s}(t)$ with powers of $t$ in the factors. See
(2.4.24) in~\cite{Sl60}. So we have to look only at differentiation
with respect to~$\k$.

For $t \in \mathbb R$ with $t>0$,
$\k-\frac{1}{2}-s \ne -1, -2, \dots$, we shall use the integral
representation (3.5.18)
in~\cite{Sl60}:
\begin{equation}
\label{W-int1} W_{\k,s}(t) \= \frac{-1}{2 \pi i}
\Gamma(\k+\frac12-s)e^{-t/2}t^{\k} \int_{(0+)}^{\infty} e^{-x}
(-x)^{s-\k-\frac12}(1+\frac{x}{t})^{s+\k-\frac12}\,dx
\end{equation}
where the contour comes from $\infty$ along a line slightly above the
positive real axis, encircles $0$ with radius $\delta<1$ and then
goes back to $\infty$ on a line slightly below the positive real
axis. By a routine computation we see that the part of the
integral over the circular part is $O(e^{\delta |t|})$. The integral
over the remaining part of the contour is $O(|t|^A)$
($A \in \mathbb R$). In all cases, the implied constants does not
depend on~$t.$ Differentiation in terms of $\k$ on $W_{\k,s}(t)$
leads to the appearance of additional factors $\log(-x)$ and
$\log(1+x/t)$ in the integrand. The arguments used in the last
paragraph imply the same estimate. Thus we get the desired
exponential decay of the perturbations of $\om_k(n,s)$.

The representation~\eqref{W-int1} is valid as long as
$-t = e^{-\pi i t}t$ is outside the path of integration. If we 
tilt the path of integration anti-clockwise by an angle $\phi$
we get a representation of
$W_{\k,s}(t)$ for $e^{-\pi i}t$ outside the new path of integration,
provided we keep $\ph\in (-\frac\pi2,\frac\pi2)$ to have convergence.
For $0<\ph<\frac \pi 2$ this gives a representation that can be used
for $\arg\bigl(e^{-\pi i}t\bigr)=0$ with $|t|>\dt$, which leads to
the desired growth of perturbations of $\hat \om_k(n,s)$.

If $\k-\frac{1}{2}-s = -1, -2, \dots$ we take $0<\ph<\frac \pi 2$ and
transform the integral representation \eqref{W-int1} into 
\begin{equation}\label{W-int2}
W_{\k,s}(t) \= \frac{e^{-\frac12 t}\, t^\k\,
e^{i\ph (s-\k+1/2)}}{\Gf(s+\frac12-\k)}\, \int_0^\infty
e^{-e^{i\ph}u}\, u^{s-\k-\frac12} \, \bigl(
1+e^{i\ph}u/t\bigr)^{s+\k-\frac12}\, du\,,
\end{equation}
Proceeding as before we obtained the same estimates.

All these estimates taken together prove \eqref{Fest}, \eqref{F0est}
(when $n \ne 0$). They further show that the derivatives of a family
with exponential decay have exponential decay and thus $\V^0_k(n,s)$
is also maximally perturbable. \smallskip

In the case $n=0$ we might use the same method. However, many
families of special functions have to be considered to cover all
cases. Instead we argue directly that we can find functions
$h^\mm_k(0,s)$ in $\V_k(0,s)$ of the form $p_\mm(x,y,\th)\,
y^{\frac12\pm 2} \, e^{ik\th}$ where $p_\mm$ is a polynomial in three
variables with degree $m_1$ in~$\th$ and degree $m_2$ in~$x$. If the
coefficient of $\th^{m_1} x^{m_2}$ in this polynomial does not depend
on~$y$, this leads to a perturbation of $y^{\frac12\pm s} e^{ik\th}$
of type~$\mm$. Such functions satisfy the required estimates, with a
polynomial factor $y^A$ instead of $e^{\dt y}$. The remaining task is
to check that they can be chosen to satisfy
$(\om-\nobreak\frac14+\nobreak s^2)\, h^\mm_k(0,s)=0$. We do this by
induction in the degrees in $\th$ and~$x$. We check that
\[ \bigl(\om-\nobreak\frac14+\nobreak s^2\bigr) x^{m_2} y^{\frac12\pm
s+a}\th^{m_1}e^{ik\th} \= -a(a\pm 2s) x^{m_2} y^{\frac12\pm
s+a}\th^{m_1}e^{ik\th} + \text{terms of lower degree in $x$ or
$\th$}\,. \]
With $a=0$ this gives the top coefficient of $p_\mm$. Moreover, the
terms of lower degree all are multiples of
$ x^{\tilde m_2} y^{\frac12\pm s+a}\th^{\tilde m_1}e^{ik\th}$ with
$\tilde m_j \leq m_j$, $\tilde m_1<m_1$ or $\tilde m_2 < m_2$, and
$a\in \ZZ_{\geq 0}$. Successively we can determine the lower degree
terms, and arrange that $h^\mm_k(0,s)$ is an eigenfunction of the
Casimir operator $\om$ with eigenvalue $\frac 14-s^2$.

This takes care of the case $n=0$, except if $s=0$. It that case
we also have to perform a computation involving
$y^{\frac12+a}\log y$, which we leave to the reader. 
\end{proof}

Holomorphic Fourier terms on $\tG$ are multiples of
\begin{equation} \label{eta1}
\eta_r(\nu;z,\th) \= y^{r/2} \, e^{2\pi i \nu z}\, e^{ir\th}\,.
\end{equation}
Thus we have the spectral parameter $s=\pm \frac{r-1}2$. For real
values of $\nu$ and $r$ we have
\begin{equation} \label{eta2}
\eta_r(\nu) \= \begin{cases}
(4\pi\nu)^{-r/2} \, \om_r\bigl(\nu,\pm \frac{r-1}2 \bigr)&\text{ if
}\nu>0\,,\\
\mu_r\bigl( 0,\frac{r-1}2\bigr)&\text{ if }\nu=0\,,\\
e^{-\pi i r}\,(4\pi|\nu|)^{-r/2}\, \hat\om_r\bigl(\nu,\pm
\frac{r-1}2\bigr)&\text{ if }\nu<0\,,
\end{cases}
\end{equation}
with notations as in \eqref{om-def-old} and~\eqref{om-mu-expl}. The
functions
\begin{equation} \label{eta3}
 \eta_k^\mm(n;z,\th) \= Q_{m_1}\bigl( \txtfrac{2i\th+\log y}{2\pi
i}\bigr)\,Q_{m_2}(z)\,\eta_k(n;z,\th)\,
\end{equation}
satisfy
\begin{equation}\label{type} \ml_{m_1+m_2}\eta_k^\mm : (\z-1)^ {l_1}\,
(\tau-1)^{l_2} \mapsto \dt_{m_1,l_1}\dt_{m_2,l_2}\,\eta_k(n)
\end{equation}
for $l_1+l_2=m_1+m_2$, and as $y\rightarrow\infty$ their growth is of
order~$\oh(e^{(\dt-2\pi n)y})$. For the commutative group $\tDt$ and
for a fixed $\mm$ they yield a basis of the space of forms of order
$m_1 + m_2 + 1$ modulo lower order forms.

As an example we note that the Fourier expansion~\eqref{FeLdef}
can be written in the following way: 
\begin{equation} L(z,\th) \= \pi i \,\eta^{(1,0)}_0(0;z,\th) +
\frac{\pi i }6\, \eta^{(0,1)}_0(0;z,\th) - 2 \sum_{n\geq 1}
\s_{-1}(n)\, \eta^{(0,0)}_0(n;z,\th)\,.
\end{equation}

\section{Proofs of Theorems \ref{thm-tMf-mp} and~\ref{thm-thol-mp}
}\label{sect-hoM}
The method of the proof is highly inductive. At each step we use the
maximal perturbability of other spaces {which has been proved in a
previous step}. The starting point for this process is the space
$\map(\tGm,\CC)$ whose maximal perturbability is proved based on
general algebraic principles in Proposition~\ref{prop-tGmfct-mp}.
This implies directly the maximal perturbability of the $\tGm$-module
$\map(\uhp\times\nobreak\RR,\CC)$. We proceed by imposing
increasingly stringent regularity conditions on the functions
$\uhp\times\RR\rightarrow \CC$. We consider
$C^\infty(\uhp\times\nobreak\RR)=C^\infty(\tG)$, the subspace
$C^\infty_k(\tG)$ of functions in $C^\infty(\tG)$ with generalised
weight~$k$ and the subspace $\C_k$ of $C^\infty_k(\tG)$ of functions
that have compact support modulo~$\tGm$. In~\S\ref{sect-Ft} we have
considered higher order invariant functions for the group $\tDt$
generated by $n(1)$ and $k(\pi)$. These functions are related to the
Fourier expansions of Maass forms. After proving that some more
auxiliary subspaces of $C^\infty_k(\uhp\times\nobreak\RR)$ are
maximally perturbable, we finally prove in~\S\ref{sect-fp} the
maximal perturbability of~$\tilde\A_k(\tGm,\ld)$ and $\Hol_k(\tGm)$.

\subsection{Higher order invariants in maps on $\tGm$} A
generalisation of Proposition~\ref{prop-Gmfct-mp} is the following:
\begin{prop}\label{prop-tGmfct-mp}If $\tGm$ is a discrete cofinite
subgroup of~$\tG$ with cusps, then the $\tGm$-module $\map(\tGm,\CC)$
(with the action by left translation) is maximally perturbable.
\end{prop}
\begin{proof} We first define (similarly to
Proposition~\ref{prop-Gmfct-mp})
$\gg_\ii$ on the free subgroup $\tGm_0$ of $\tGm$ generated by
$\al_1,\ldots,\allowbreak\al_{\ngen
-1}$ for $\ii \in \{1,\ldots,\ngen\!-\!\nobreak1\}^q$ by the relations
in~\eqref{gidef}, with $A_j$ replaced by~$\al_j$.

Let $\ph_0: \tGm\rightarrow\tGm_0$ be the surjective group
homomorphism given by $\ph_0(\al_j)=\al_j$ for $1\leq j \leq
\ngen\!-\!1$, $\ph_0(\z)=1$ and $\ph_0(\e_j)=1$ for $1\leq j \leq
\nell$. For $1\leq j \leq \ngen$ we define $\ps_j \in \hom(\tGm,\CC)$
such that $\ps_j(\al_{j'})=\dt_{j,j'}$. This determines $\ps_j$
completely, since values on elliptic generators are given by
$\ps_j(\e_j) = \frac1{v_j} \ps_j(\z)$. For
$\ii=(\ii',\ngen,\ldots,\ngen)$ with $m$~coordinates $\ngen$ at the
end and $\ii'\in \{1, \dots, \ngen\!-\!\nobreak1\}^{q-m}$, we put
\begin{equation}\label{fidef}
\ff_\ii(\gm) \= \gg_{\ii'}\bigl(\ph_0(\gm)\bigr) \, Q_{m}\bigl(
\ps_\ngen(\gm)
\bigr)
\end{equation}
where $Q_n$ are the polynomials defined in~\eqref{Qndef}. Now we can
check the following properties of~$\ff_\ii$:
\begin{align}
\label{rec0}
\ff_{()} &\= 1\,,\qquad(\text{empty tuple, $q=0$})\,;\\
\label{cst}
\ff_\ii(1)&\=0\qquad\text{if $|\ii|\geq 1$}\,\,;\\
\label{rec1}
\ff_\ii|(\z-1) &\= \begin{cases}
\ff_{\ii'}&\text{ if } \ii \= (\ii',\ngen)\,,\\
0&\text{ if $\ii$ does not end with a $\ngen$}\,;
\end{cases}\\
\label{rec2}
\ff_\ii|(\al_j-1)&\= \begin{cases}
\ff_{\ii'}&\text{ if $\ii=(j,\ii')$ with $j<\ngen$}\,,\\
0&\text{ if $j<\ngen$, $j\neq \ii(1)$}\;.
\end{cases}
\end{align}

Using this we can see that
\begin{equation}(\ml_q \ff_\ii )\bigl(\bb(\jj) \bigr)\=
\dt_{\ii,\jj}\,.
\end{equation}

Now, the choice of the basis $\bb(\ii)$ in~\eqref{bi-def} for
$\tGm$-$q$-tuples~$\ii$ shows that to prove that Map$(\tGm, \CC)$ is
maximally perturbable it suffices to prove that for each $\ii$ and
for each function $f$ on $\tGm\backslash \tG$ a function $h_\ii
\in \map(\tG,\CC)$ such that for all $\tGm$-$q$-tuples~$\jj$:
\begin{equation} \label{hrel}
 h_\ii | (\al_{\jj(1)}-1) \cdots (\al_{\jj(q)} - 1) \= \dt_{\ii,\jj}
 \cdot f\,.\end{equation}

To construct such functions we choose a \emph{strict fundamental
domain} $\fd_\tGm \subset \tG$ for $\tGm\backslash \tG$, {\sl i.e.},
a set meeting each $\tGm$-orbit exactly once. Such a fundamental
domain can be constructed from a strict fundamental domain $\fd_\uhp$
for $\Gm\backslash \uhp$, by taking
\begin{align*}
\fd_\tGm &\= \bigl\{ (z,\th) \;:\; z\in \fd_\uhp\,, 0\leq \th \leq
\pi/n_z\bigr\}\,,\\
n_z&\= \min\bigl\{n\in \NN\;:\; \text{ there is $\gm\in\tGm$ fixing
$z$ in $\uhp$ conjugate to }k(\pi/n)\bigr\}\,.
\end{align*}
So $n_z=1$ for all $z\in \fd_\uhp$, except for the elliptic fixed
points $z_1,\ldots,z_\nell$ in~$\fd_\uhp$. These are conjugate to a
fixed point of $\e_j$ and $n_{z_j} = v_j$.

A choice for the sought function $h_\ii$ is then
\begin{equation}
\label{hidef-gen} h_\ii(\gm g) \= \ff_\ii(\gm) f(g)\qquad\gm\in
\Gm\,,\; g\in \fd_\tGm\,. \end{equation}
With the characteristic function $\ps$ of $\fd_\tGm$, we can write
this as
\begin{equation}\label{hdef-pu}
h_\ii(g) \= \sum_{\gm\in \tGm} \ff_\ii(\gm) \, f(g)\,
\ps(\gm^{-1}g)\,.
\end{equation}
\end{proof}

\subsection{Higher order invariants in smooth functions on~$\tG$} We
will use essentially the same construction as in the last section to
prove that
\begin{prop}The $\tGm$-module $C^\infty(\tG)$ is maximally
perturbable.
\end{prop}
\begin{proof}In order to show that $C^\infty(\tG)$ is a maximally
perturbable $\tGm$-module, we need to have \eqref{hrel} with
$h_\ii\in C^\infty(\tG)$ for each
$f\in C^\infty(\tGm\backslash \tG)$. Lemma~\ref{lem-pu} in Appendix A
shows that we can find functions $\ps\in C^\infty(\uhp\times\RR)$
such that $\sum_{\gm\in \tGm} \ps\bigl(\gm^{-1}(z,\th) \bigr)=1$ for
all $(z,\th)\in \uhp\times\RR$ as a locally finite sum. If we define
\eqref{hdef-pu} with such a function~$\ps$ and $f\in 
C^\infty(\tGm\backslash \tG)$, then the sum is locally finite, and
the $h_\ii$ are smooth.
\end{proof}

\subsection{Higher order invariants and generalised weight} Set
\begin{equation} \label{C_k}
C^\infty_k(\tG)=\{f \in C^\infty(\tG) \, , \text{of generalised weight
$k$} \} \,.
\end{equation}
\begin{prop}\label{prop-smk-mp}Let $k\in 2\ZZ$. Then the $\tGm$-module
$C^\infty_k(\tG)$ is maximally perturbable.
\end{prop}
\begin{proof}
As with the previous proofs, our approach is to show that for every
$\tGm$-$q$-tuple $\ii = (\ii', \ngen,\ldots,\ngen)$ with exactly $m$
occurrences of~$\ngen$ at the end and for every $f\in
C^\infty_k(\tGm\backslash \tG)$ there exists $h_\ii \in
C^{\infty}_k(\tG)$ satisfying equation \eqref{hrel} for all
$\tGm$-$q$-tuples $\jj$. We note that, by Theorem~\ref{thm-tAi}, the
$\tGm$-invariance of~$f$ implies that its weight~$k$ is strict, {\sl
i.e.}, $f(gk(\th))= f(g)e^{ik\th}$.

We will define the function $h_\ii$ by an analogue of~\eqref{hdef-pu}.
We first define for each $g\in \tG$ the point
$w(g)=\pr(g)\, i\in \uhp$ and the real number $\Th(g)\in \RR$ such
that $g=\bigl(w(g) ,\Th(g)\bigr)\in \tG = \uhp\times\RR$. We also
recall that $\Gm=\tGm/\tZ$. Since the group homomorphism $\phi_0$
defined in the proof of Proposition~\ref{prop-tGmfct-mp} is trivial
on $\tZ=\langle\zeta\rangle$, it induces a homomorphism on $\Gamma$.
Now we take $\ps(z,\th)=\ps_0(z)$, with $\ps_0$ as in Part~ii) of
Lemma~\ref{lem-pu}. So the function $(z,\th) \mapsto
\ps\bigl(\gm^{-1}(z,\th)
\bigr)$ obtained by left translation depends only on the image of
$\gm\in \tGm$ in $\Gm\cong\tGm/\tZ$. Let, as in the proof of
Proposition~\ref{prop-tGmfct-mp}, $\ps_{\ngen}$ be the function
$\tGm\rightarrow\RR$ such that $\ps_{\ngen}(\al_{j'}) = \dt_{\ngen,j'}$.
For a given $\gm\in \tGm$ we have
$\ps_\ngen(\z\gm) = \ps_\ngen(\gm)+1$ and $\Th\bigl(
(\z\gm)^{-1}g) = \Th(\gm^{-1}g)-\pi$. So $\ps_\ngen(\gm) +
\Th(\gm^{-1}g)/\pi$ is well-defined on $\Gm=\tGm/\tZ$. We can
therefore set \rm
\begin{equation}\label{hidef-wt}
h_\ii(g) \= \sum_{\gm\in \Gm} \gg_{\ii'} \bigl(\ph_0(\gm)\bigr)\,
Q_{m} \bigl( \ps_\ngen(\gm) + \Th(\gm^{-1}g)/\pi\bigr)\, f(g)\,
\ps\bigl(\gm^{-1}g\bigr)\,.
\end{equation}
The support property of the partition of unity~$\ps$ ensures
convergence; it is even a locally finite sum with a bounded number of
non-zero terms. All factors depend smoothly on~$g$. So $h_\ii \in
C^\infty( \tG)$.

We consider $(\WW-\nobreak ik) h_\ii$. Since $\WW \ps=0$, we need only
consider
\begin{equation}\
\begin{aligned} \Bigl( \partial_\th &- ik \Bigr) Q_{m}\Bigl(
\psi_{t(\Gm)} (\gm) + \Th(\gm^{-1}g k(\th))/\pi\Bigr) \, f\bigl(g
k(\th)\bigr)\\
&\= Q_{m}\Bigl( \psi_{t(\Gm)} (\gm) + \Th(\gm^{-1}g k(\th))/\pi\Bigr)
\,\Bigl(
\partial_\th - i k \Bigr) f\bigl(g k(\th)\bigr)\\
&\qquad\hbox{}
+ f\bigl(g k(\th)\bigr) \; \partial_\th Q_{m}\Bigl( \psi_{t(\Gm)}
(\gm) + \Th(\gm^{-1}g)/\pi+\th/\pi\Bigr)\\
&\= 0 + \pi^{-1} Q_{m}'\Bigl( \psi_{t(\Gm)} (\gm) +
\Th(\gm^{-1}g)/\pi+\th/\pi\Bigr)\,f\bigl(g k(\th)\bigr)\,.
\end{aligned}
\end{equation}Repeating this we obtain
\begin{equation}
\begin{aligned}
(\WW-\nobreak i k)^{m +1 } &Q_{m}\Bigl( \ps_\ngen(\gm) + \Th(\gm^{-1}g
)/\pi\Bigr) \,f(g)
= \pi^{-m -1 } Q_{m}^{(m +1 )}\bigl( \cdots\bigr)\,\cdots\=0
\end{aligned}
\,,\end{equation}
since the degree of $Q_{m}$ is~$m$. So $h_\ii \in C^\infty_k(  \tilde
G )$.
\end{proof}
In a similar (but much simpler) way, one shows that, if $\Gm$
acts on $C^{\infty}(\uhp)$ via \eqref{action} and $f \in 
C^{\infty}(\uhp)^{\Gm}$, then the function in $C^{\infty}(\uhp)$
$h_\ii(z):= \sum_{\gm\in \Gm} \gg_{\ii}  \bigl(\ph_0(\gm)\bigr)\,
f(z)\, \ps_0\bigl(\gm^{-1}z\bigr)$ satisfies \eqref{hrel} for all
$(t(\Gm)-1)$-tuples of elements of $\{1, \dots, t(\Gm)-1\}$. This
gives an alternative proof of
\begin{prop} [Prop.~4.1, \cite{DSr}] Let $k\in 2\ZZ$. Then the
$\Gm$-module $C^\infty(\uhp)$ is maximally perturbable.
\end{prop}
In fact, since $\psi_0$ is bounded, if $f$ has polynomial growth at
all cusps, then so does $h_{\ii}$ thus proving that the submodule of
$C^{\infty}(\uhp)$ of functions with polynomial growth is also
maximally perturbable. 

\subsection{Higher order invariants with support conditions} We shall
first discuss the motivation for the introduction of the invariants
we will be dealing with. If Definition~\ref{defTa} of the space
$\tilde\A_k(\tGm,\ld)$ did not include a growth condition at the
cusps, we could consider $\tilde\A_k(\tGm,\ld)$ as the
kernel~$\mathcal{K}$ in the exact sequence \[ 0 \longrightarrow
\mathcal{K} \longrightarrow C^\infty_k(\tG)
\stackrel {\om-\ld} \longrightarrow C^\infty_k(\tG) \]
With exponential growth, one might want to try to replace
$C^\infty_k(\tG)$ by its subspace $C^\infty_l(\tGm)^{\mathrm{eg}}$ of
functions with exponential growth at the cusps of~$\tGm$. This would
lead to an exact sequence
\[ 0 \longrightarrow \tilde \A_k(\tGm,\ld) \longrightarrow
C^\infty_k(\tGm)^{\mathrm{eg}} \stackrel {\om-\ld} \longrightarrow
C^\infty_k(\tGm)^{\mathrm{eg}} \]
for which we might try to show that for each~$q\in \NN$
\[ 0 \longrightarrow \tilde \A_k(\tGm,\ld)^{\tGm, q} \longrightarrow
\ivt{\bigl(C^\infty_k(\tGm)^{\mathrm{eg}}\bigr)}q \stackrel {\om-\ld}
\longrightarrow \ivt{\bigl(C^\infty_k(\tGm)^{\mathrm{eg}}\bigr)}q \]
is exact. For this to be of use it seems that we need surjectivity of
the map $\om-\ld: \bigl(C^\infty_k(\tGm)^{\mathrm{eg}}\bigr)^{\tGm}
\rightarrow \bigl(C^\infty_k(\tGm)^{\mathrm{eg}}\bigr)^{\tGm}$, which
we did not succeed in proving, and which may not hold. For this
reason we will instead work with other better behaved subspaces of
the spaces appearing in the exact sequence. We will therefore
define subspaces $\C_k, \D_k(\ld) \subset C_k^{\infty}(\tG)$ and
$  \A_k'(\ld) \subset  \A_k(\tGm, \ld)$ related by an exact sequence
\begin{equation}\label{es}
0 \longrightarrow \tilde \A_k'(\ld) \longrightarrow \D_k(\ld)
\stackrel{\om-\ld }\longrightarrow \C_k.
\end{equation}

\subsubsection{The spaces $\C_k$}\label{sect-spC} For each cusp
$\k=\tilde g_\k \infty$, and each $a>0$ we call
\begin{equation}\label{Dk-def}
D_\k(a) \= \tilde g_\k \bigl\{ (z,\th)\;:\; \im z\geq a,\, \th \in \RR
\bigr\}
\end{equation}
a \emph{horocyclic set}. There is a number $A_\Gm$ such that for each
$a\geq A_\Gm$ the $D_\k(a)$ are disjoint for different cusps. The
sets
\begin{equation} \label{tGa-def}
\tG_a = \bigl\{ (z,\th) \in \uhp\times\nobreak\RR \;:\; \forall \k \;
(z, \th) \not\in D_\k(a)\bigr\}
\end{equation}
satisfy $\tilde\Gm\tG_a=\tG_a$. This follows from the fact that the
$g_\k$ have been chosen so that 
\begin{equation} \label{commut}
\gm \tGm_\k \tilde g_\k =\tilde g_{\gm\k} \tGm_{\infty}
\end{equation} for all cusps $\k$ and for $\gm\in \tGm$. Here
$\tGm_\k: = \pr^{-1}\Gm_\k = \bigl\{ \gm\in
\tGm\;:\; \gm\k=\k\bigr\}$.

\begin{defn} \label{defn-ck}
Let $k\in 2\ZZ$. The space $\C_k$ consists of the $f\in
C^\infty_k(\tG)$ supported in $\tG_a$ for some $a\geq A_\Gm$.
(The $a$ may depend on~$f$).
\end{defn}
So $\C_k$ consists of the smooth functions with generalised weight~$k$
whose supports project to compact subsets of $\Gm\backslash\uhp$.
Clearly, the space~$\C_k$ is $\tGm$-invariant. If we apply the
construction of $h_\ii$ in the proof of Proposition~\ref{prop-smk-mp}
to functions $f\in \C_k^{\tGm} \subset
C^\infty_k(\tGm\backslash \tG)$ then the support of each $h_\ii$ is
contained in the same set $\tG_a$ that contains $\supp(f)$. This
implies:
\begin{prop}\label{prop-Ck-mp}Let $k\in2\ZZ$. Then the $\tGm$-module
$\C_k$ is maximally perturbable.
\end{prop}

\subsubsection{The spaces $\D_k(\ld)$} The construction of $\D_k(\ld)$
and the proof of its maximal perturbability is much lengthier that
those for $\C_k$. We will define $\D_k(\ld)$ essentially as the space
of functions that accept higher-order analogues of Fourier expansions
at the cusps. To make this formal we study spaces of functions
defined on regions of the form
\begin{equation} \label{S(y)}
S(y_0) \= \bigl\{ (x+iy,\th) \in \uhp\times\RR\;:\; y>
y_0\bigr\}\,,
\end{equation}
with $y_0>0$.
\begin{defn}\label{Xdef}Let $k\in 2\ZZ$, $\ld\in \CC$, and $y_0>0$. We
denote by $\A_k(y_0,\ld)$ the space of $f\in
 C^\infty \bigl( S(y_0) \bigr)$ that satisfy $\om f=\ld f$,
$(\WW-\nobreak ik)^nf=0$ for some $n\in \NN$, and have at most
exponential growth as $y\rightarrow\infty$, uniform for $x$ and $\th$
in compact sets. We denote by $\A\hol_k(y_0)$ the space of
holomorphic functions on~$S(y_0)$ with generalised weight~$k$ and at
most exponential growth as $y\rightarrow\infty$
\end{defn}

\begin{prop}\label{prop-Amp}Let $k\in 2\ZZ$, $s\in \CC$ and $y_0>0$.
The spaces $\A_k(y_0,\ld_s)$ and $\A\hol_k(y_0)$ are maximally
perturbable $\tDt$-modules.

Let $q\in \NN$. Each $f \in \ivv{\A_k(y_0,\ld_s)}\tDt q$ has an
absolutely convergent expansion
\begin{equation}\label{hoexp}
f(z,\th) \= \sum_{n\in \ZZ} f_n(z,\th)\end{equation}
on $S(y_0)$ with $f_n \in \ivv{\V_k(n,s)}\tDt q$ for all $n$, 
and $f_n\in
\ivv{\V_k^0(n,s)}\tDt q$ for almost all~$n$.

Each $f\in \ivv{\A\hol_k(y_0)}\tDt q$ has an absolutely convergent
expansion on $S(y_0)$ of the form
\begin{equation}
\label{hoholFexp} f(z,\th) \= \sum_{\mm, \, m_1+m_2 < q} \sum_n
c^n_\mm\, \eta^\mm_k(n; z,\th)\,\end{equation}
where the inner sum ranges from some, possible negative, 
integer to infinity.
\end{prop}
\begin{proof}
We start with the holomorphic case. Let $f\in \A_k\hol(y_0)^\tDt$.
Then the function $z\mapsto y^{-k/2} \, f(z,0)$ is holomorphic on
$\{z\in \uhp\;:\; y>
y_0\}$ with period~$1$. So it has a finite to the left expansion
of the form $\sum_n a_n e^{2\pi i n z}$ converging absolutely on $y>
y_0$. For each $y_1>y_0$ we have
$a_n = \oh\bigl( e^{2\pi n y_1}\bigr)$ as $n\rightarrow\infty$.

Hence $f(z,\th) = \sum_n a_n \eta_k(n; z,\th)$ converges absolutely
on~$y> y_0$, and
\[ f^\mm(z,\th) \;:=\; \sum_{n\geq -N} a_n \eta^\mm_k(n; z,\th)\]
converges absolutely on $S(y_0)$, and the convergence is uniform on
any set $y\geq y_1$ with $y_1>y_0$, with $x$ and $\th$ in compact sets. 
These functions satisfy $f^\mm
|(\tau-\nobreak 1) = f^{(m_1,m_2-1)}$, $f^\mm|(\z-\nobreak 1) =
f^{(m_1-1,m_2)}$ and $f^{(0,0)}=f$, since all $\eta^\mm_k$ have this
property. Thus $f^\mm$, with $\mm$ such that $m_1+m_2<q$ is a
perturbation of type $\mm$ and we deduce that $\A_k\hol(y_0)$ is
maximally perturbable.
An arbitrary element $h \in \ivv{\A_k\hol(y_0)}\tDt q$ can be written as 
a finite linear combination of such $f^{\mm}$,
which all have
expansions of the type given in~\eqref{hoholFexp}. 
\smallskip

For $f\in \A_k(y_0,\ld_s)^\tDt$ we proceed similarly. By
Theorem~\ref{prop-Fe} in combination with Remark \ref{rem-Fe}
and the integrality of $k$, there is an absolutely converging
Fourier expansion
\[ f(z,\th) \= \sum_{n\in \ZZ} f_n(z,\th)\]
on $S(y_0)$ with $f_n \in \W_k(n,s)$. By the exponential growth,
$f_n \in \W_k^0(n,s)$ for $|n|>N$, for some $N\in \NN$.

For $|n|>N$ we have $f_n = a_n \om_k(n,s)$, and from~\eqref{om-asympt}
we conclude that $a_n = \oh \bigl( e^{2\pi|n|y_1})$ as
$|n|\rightarrow\infty$ for each $y_1>y_0$. So by \eqref{om-asympt}
the series
\[ \sum_{n,\,|n|>N} a_n \om^\mm_k(n,s)\]
converges absolutely on~$S(y_0)$, uniformly on each set $y\geq y_1$
with $y_1>y_0$, and gives an exponentially decreasing function as
$y\rightarrow\infty$. It is a $\ld_s$-eigenfunction of~$\om$, since
the decay allows differentiation inside the sum. To produce a
perturbation $f^\mm$ of~$f$ we pick $f^\mm_n \in
\ivv{\V_k(n,s)}\tDt {m_1+m_2+1}$ such that $f^\mm_n |(\tau-\nobreak
1) = f_n^{(m_1,m_2-1)}$, $f^\mm_n|(\z-\nobreak 1) =
f^{(m_1-1,m_2)}_n$ and $f^{(0,0)}_n=f_n$ for the finitely many $n$
with $|n|\leq N$. The estimate \eqref{Fest} shows that the growth of
these terms is at most of the order $\oh(e^{(2\pi N+\dt)y})$ as
$y\rightarrow\infty$ for each $\dt>0$. Thus we get (non-uniquely) a
perturbation of type~$\mm$
\[ f^\mm \= \sum_{n,\, |n|\leq N} f^\mm_n + \sum_{n,\,|n|>N} a_n
\om^\mm_k(n,s)\]
in $\A_k(y_0,\ld_s)$. Thus we get \eqref{hoexp} and the maximal
perturbability of~$\A_k(y_0,\ld_s)$.
\end{proof}

We are now ready to define $\D_k(\ld)$ and $\D\hol_k$.
\begin{defn}\label{Ddef}Let $k\in 2\ZZ$, and $\ld\in \CC$. We define
$\D_k(\ld)$ as the space of functions $f\in C^\infty_k(\tilde G)$
(hence with generalised weight~$k$) for which there exist
$b\geq A_\Gm$, $a\in \RR$, and $q\in\NN$ such that for each cusp $\k$
of $\tGm$ the function $(z,\th) \mapsto f\bigl(\tilde g_\k(z,\th)
\bigr)$ is an element of $\ivv{ \A_k(b,\ld)}\tDt q$, and satisfies a
bound $\oh(e^{ay})$ as $y\rightarrow\infty$.

We define $\D\hol_k$ similarly, with $(z,\th) \mapsto f\bigl(\tilde
g_\k(z,\th)
\bigr)$ in $\A_k\hol(b)^{\tDt, q}$, with bound $\oh(e^{ay})$.
\end{defn}

\begin{rmk} The numbers $a$,
$b$ and $q$ may depend on the function~$f$. 
\end{rmk}

\begin{rmk} \label{eigenproperty}
Definition~\ref{Xdef} of $\A_k(b,\ld)$ implies that elements of
$\D_k(\ld)$ are $\ld$-eigenfunctions of~$\om$ on the set
$\bigsqcup_\k D_\k(b)$. Similarly, elements of $\D\hol_k$ are
holomorphic functions on $\bigsqcup_\k D_\k(b)$. In both cases we
have exponential growth at each cusp. The definition requires that
the order of this exponential growth stays bounded when we vary the
cusp.
\end{rmk}

\begin{rmk}
In the definition we impose $\tDt$-invariance of bounded order near
all cusps. This is a bit artificial, but serves our purpose.
\end{rmk}

The space $\C_k$ is contained in $\D_k(\ld)$ and in $\D\hol_k$.
Indeed, for given $f\in \C_k$ we can take $b$ large so that
$\bigsqcup_\k D_\k(b)$ is outside the support of~$f$.

Elements $f$ of $ \tilde\A_k(\tGm,\ld)^\tGm$ restricted
to $D_\k(b)$ induce elements $(z,\th) \mapsto f\bigl(
\tilde g_\k(z,\th)
\bigr)$ in $\A_k(b,\ld)^\tDt$ for each cusp~$\k$, and similarly in the
holomorphic case. Hence
\begin{equation}
\label{inclusion}
\tilde\A_k(\tGm,\ld)^\tGm \;\subset\; \D_k(\ld)^\tGm\,,\qquad
\Hol_k(\tGm)^\tGm\;\subset\; (\D\hol_k)^\tGm\,.
\end{equation}

\it Maximal perturbability of $\D_k(\ld)$ and $\D\hol_k$. \rm We first
need a technical lemma in order to relate $\tDt$-invariants to
$\tGm$-invariants.

We first note that if $\infty$ is a cusp of~$\tGm$ and if
$\tilde g_\infty=1$, then $\tDt = \tGm_\infty$. In general the group
$\tGm_\k$ can be conjugated to $\tilde g_\k^{-1}\tGm_\k \tilde g_\k =
\tDt$ in $\tilde g_\k^{-1}\tGm\tilde g_\k$. So we can assume here
that $\tDt \subset\tGm$.

The abelian group $\tDt$ is  free on the 
generators $\tau=n(1)$ and
$\z=k(\pi)$. The dimension of $\ivv{\map(\tDt,\CC)}{\tDt}{q+1}$ is
$(q+1)(q+2)/2$ 
with an explicit basis described as
follows. Define a sequence of maps on $\tDt$ by setting
\begin{equation} \label{philmrel}
\begin{aligned}
\ph^{(l,m)}|(\z-1)&\= \ph^{(l-1,m)} \\
\ph^{(l,m)}|(\tau-1)&\=\ph^{(l,m-1)} \\
\end{aligned}
\end{equation}
and 
\begin{align*}\ph^{(0,0)}&\=1, \quad \ph^{(l,m)}=0 \, \text{ for $l$ or
$m$ negative} \\
\text{and} \, \ph^{(l,m)}(1)&\=0\quad\text{for $l,m\geq 0$, $l+m>0$}.
\end{align*}
 Then
$$(\ml_{q} \ph^{(l, m)}) \left ((\z^r-1)(\tau^s-1)
\right ) \=\dt_{l,r} \dt_{m, s}$$ for $l+m=r+s=q$, and therefore the
$\ph^{(l, m)}$ with $l, m \ge 0, l+m \le q$ is a basis of
$\ivv{\map(\tDt,\CC)}{\tDt}{q+1}$.

Let $R $ be a system of representatives of $\tGm/\tDt$; so
$R\subset\tGm$. Consider the system
$\{\ff_\jj\}_{|\jj|=q} \subset \ivt{\map(\tGm,\CC)}{q+1}$ in the
proof of Proposition~\ref{prop-tGmfct-mp}. If $|\jj|=q$, then, for
every $\gm \in \tGm$, $\dt\mapsto \ff_\jj(\gm\dt)$ is a function
on~$\tDt$ of order at most~$q+1$. Hence there are functions
$a^\jj_{l,m}$ on~$R$ such that for all $\rho\in R$ and $\dt\in \tDt$
\begin{equation}\label{aiidef}
\ff_\jj(\rho\dt) \= \sum_{l,m\geq 0\,,\; l+m\leq q}
a^\jj_{l,m}(\rho)\, \ph^{(l,m)}(\dt)\,.
\end{equation}

\begin{lem}\label{lem-lo}Let $a^\jj_{l,m}$ be as in~\eqref{aiidef},
and suppose that we have functions $\ps^{(l,m)}\in \map(\tDt,\CC)$
satisfying
\begin{equation}\label{psilmrel}
\begin{aligned}
\ps^{(0,0)}&\=0\,,\\
\ps^{(l,m)}|(\tau-1)&\= \ps^{(l-1,m)}\quad \text{ for }l\geq 1\,,\\
\ps^{(l,m)}|(\z-1)&\=\ps^{(l,m-1)}\quad\text{ for }m\geq 1\,.
\end{aligned}
\end{equation}
Then
\begin{equation}\label{oq} f(\rho\dt)\= \sum_{l,m\geq 0\,,\; l+m\leq
q} a^\jj_{l,m}(\rho)\, \ps^{(l,m)}(\dt) \quad(\rho\in R\,,\; \dt\in
\tDt)
\end{equation}
defines an element of $\ivt{\map(\tGm,\CC)}q$.
\end{lem}
\begin{proof} We proceed by induction in $q=|\jj|.$ If $q=0$, then
$m=n=0$, so $f(\rho\dt)= a^\jj_{0,0}(\rho)
\cdot
\ps^{(0,0)} \=0\in \ivt{\map(\tGm,\CC)}{0}=\{0\}$.

It is clear that \eqref{oq} gives a well-defined map on $\tGm$. It
suffices to prove that, for every generator $\al_j$ of $\tGm$,
$f|(\al_j-1) \in\map(\tGm, \CC)^{q-1}.$ Suppose first that
$\jj=(j,  \jj')$. For each $\rho\in R$ there are unique $\rho_1\in R$
and $\dt_1\in \tDt$ such that $\al_j\, \rho = \rho_1 \dt_1$.
From~\eqref{aiidef} it follows that
\begin{equation}
\begin{aligned} \label{decomp}
\ff_\jj|&(\al_{\jj(1)}-1)(\rho\dt)\\
&\= \sum_{l,m\geq 0\,,\; l+m\leq q}
a^\jj_{l,m}(\rho_1) \; \ph^{(l,m)}|(\dt_1-1) (\dt)
\quad\hbox{} + \sum_{l,m\geq 0\,,\; l+m\leq q} \bigl(
a^\jj_{l,m}(\rho_1) - a^\jj_{l,m}(\rho) \bigr) \; \ph^{(l,m)}
(\dt)\,.
\end{aligned}
\end{equation}
By \eqref{rec2}, the left-hand side equals
$\sum_{l,m\geq 0\,,\; l+m\leq q-1} a^{\jj'}_{l,m}(\rho)\,
\phi^{(l,m)}(\dt) $.  The function $\ph^{(l,m)}|(\dt_1-\nobreak 1)$
is a linear combination, depending on~$\rho$, of $\ph^{(a,b)}$ with
$0\leq a\leq l$, $0\leq b \leq m$ and $a+b\leq q-1$. Thus we get an
expression that expresses the $a^{\jj'}_{l,m}(\rho)$ in the
$a^\jj_{l,m}(\rho)$. The form of this expression depends on the
relations \eqref{philmrel} but not on the specific value of the
constant basis element $\ph^{(0, 0)}$. The relations of
\eqref{philmrel} hold for $\ps^{(l,m)}$ too. Therefore, the right
hand side of \eqref{decomp}, upon replacement of $\phi$ by
$\psi$, equals
$$\sum_{l,m\geq 0\,,\; l+m\leq q-1} a^{\jj'}_{l,m}(\rho)\,
\ps^{(l,m)}(\dt) \quad(\rho\in R\,,\; \dt\in \tDt)\,,
$$
which, by induction, is in $\map(\tGm, \CC)^{\tGm, q-1}$. Since, from
\eqref{oq}, it follows that the right hand side of \eqref{decomp}
with $\ph$ replaced by $\psi$ equals $f|(\al_j-1)$ too, we
deduce that
$f|(\al_j-1) \in \allowbreak \map (\tGm, \CC)^{\tGm, q-1}$.

In the same way, we deduce that $f|(\al_j-1)\in \map (\tGm,
\CC)^{\tGm, q-1}$ when $j=t(\Gm)$ or $j<t(\Gm)$ and $j \neq \jj(1)$.
\end{proof}

\begin{prop}\label{prop-Dk-mp}The $\tGm$-modules $\D_k(\ld)$ and
$\D\hol_k$ are maximally perturbable for all $k\in 2\ZZ$ and $\ld\in
\CC$.
\end{prop}
\begin{proof}
It suffices to construct for a given $f\in \D_k(\ld)^\tGm$, a given
$q\in \NN$ and a given $\tGm$-$q$-tuple~$\ii$, an element
$\eta_\ii\in \D_k(\ld)$ such that
$\eta_\ii|(\al_{\ii'(1)}\nobreak-1)\cdots
(\al_{\ii'(q)}-\nobreak 1)=\delta_{\ii, \ii'}f$ for all
$\tGm$-$q$-tuples~$\jj$.

We will write $f=f_\cp + \sum_\k f_\k$, with $\k$ running over a set
$C$ of representatives of the $\tGm$-orbits of cusps, where
$f_\cp  \in (\C_k)^\tGm$, $f_\k \in \D_k(\ld)^\tGm$, and will produce
perturbations for each of these components.

We choose a strict fundamental domain $\fd_\tGm$ for 
$\tGm\backslash \tG$ 
so that
\[ \fd_\tGm \cap D_\infty(b) \= \bigl\{ (x+iy,\th)\;:\; 0\leq x<1\,,\;
y\geq b\,,\; 0\leq \th<\pi\bigr\}\,.\]
Definition~\ref{Ddef} provides $b\geq A_\Gm$ and $r\in \NN$ such that
$v_\k(z,\th) = f\bigl(\tilde g_{\k} (z,\th)
\bigr)$ is in $\ivv{\A_k(b,\ld)}\tDt r$ for each cusp~$\k$.
Furthermore, $b$ can be chosen large enough for the sets 
$\fd_\tGm
\cap D_{\k}(b)$ ($\k \in C$) to be pairwise disjoint. Since $f$ is
$\tGm$-invariant, we even have $v_\k\in \A_k(b,\ld)^\tDt$. We choose
a function $\ch\in C^\infty(0,\infty)$ that is equal to $0$ on
$(0,b+\frac{1}{2}]$ and equal to~$1$ on $[b+\nobreak 1,\infty)$, and
define for $\k\in C$
\begin{equation}
\label{f-kap-def}
 f_\k \bigl(z,\th \bigr) \= \begin{cases} 0 & (z, \th) \in \fd_\tGm-
D_{\k}(b) \\
\ch(\text{Im}(z_1))\, v_\k(z_1,\th_1) & (z, \th)=\tilde g_{\k}(z_1,
\th_1)
\in \fd_\tGm \cap D_{\k}(b)
\end{cases}
\end{equation}
Extend to $\tilde G$ by $\tGm$-linearity. So $f_\k =0$ outside $\tGm
D_\k(b)$ and equal to $f$ on $\tGm D_\k(b\!+\!\nobreak 1)$. We check
in Definition~\ref{Ddef} that $f_\k\in \D_{\k}(\ld)$. The function
$$f_\cp = f - \sum_{\k\in C} f_\k$$
is $\tGm$-invariant and vanishes on $D_\k(b\!+\!\nobreak 1)$ for all
cusps~$\k$, hence $f_\cp\in \C_k^\tGm$.

Proposition~\ref{prop-Ck-mp} implies that there is $h_\ii \in \C_k
\subset \D_k(\ld)$ satisfying the conditions
$h_\ii|(\al_{\ii(1)}\nobreak-1)\allowbreak\cdots\allowbreak
(\al_{\ii(q)}-\nobreak 1)=f_\cp$ and
$h_{\ii'}|(\al_{\ii'(1)}\nobreak-1)\allowbreak\cdots\allowbreak
(\al_{\ii'(q)}-\nobreak 1)=0$ for $\tGm$-$q$-tuples $\ii'\neq \ii$. So
we can restrict our attention to the $f_\k$.

Since the supports of the $f_\k$ with $\k\in C$ are disjoint, we can
consider each of the $f_\k$ separately. Without loss of generality we
can assume that $\infty$ is a cusp of~$\tGm$ with
$\tilde g_{\k}=1$ and take $\infty \in C$. Conjugation
by 
the original~$\tilde g_\k$ then gives the same result for a general $\k\in C$.

The function $v_\infty$ used in~\eqref{f-kap-def} is an element of
$\A_k(b,\ld)^\tDt$. The proof of Proposition~\ref{prop-Amp} shows
that for each $\mm\in \NN_0^2$ there is a perturbation
$v^\mm_\infty\in \ivv{\A_k(b,\ld)}\tDt{m_1+m_2+1}$ of
$(z,\th) \mapsto f_{\infty}(z,\th)$ of type~$\mm$. We define
$\eta_\ii$ by $\eta_\ii=0$ on~$\tG_b$ and on all $\tGm D_\k(b)$ for
all $\k\in C\setminus\{\infty\}$, and
\begin{equation}
\eta_\ii\bigl(\rho (x+iy,\th) \bigr) \= \sum_{l,m\geq 0\,,\; l+m\leq
q} \ch(y)\, a^\ii_{l,m}(\rho)\, v_\infty^{(l,m)}(x+iy,\th)
\end{equation}
for $y\geq b$ and $\rho$ in a system of representatives $R$ of
$\tGm/\tDt$. The functions $a^\ii_{l,m}$ are as in~\eqref{aiidef}.
Since the sets $\rho D_\infty(b)$ are disjoint, this defines a smooth
function, which can be checked to be an element of $\D_k(\ld)$.

For each fixed $g=(x+\nobreak iy,\th)$ with $y\geq b$ the function
$\dt \mapsto v_\infty^{(l,m)}(\dt g)$ on~$\tDt$ satisfies the same
relations as $\dt\mapsto \ph^{(l,m)}(\dt) \, v_\infty(g)$
in~\eqref{philmrel}. So, their difference, as a function of $\dt$,
satisfies~\eqref{psilmrel}.

Ignoring smoothness for a moment we have $f_\infty\in
\map(\tG,\CC)^\tGm$. Equation \eqref{hidef-gen} gives a function
$h_\ii$ on $\tG$ such that
$h_\ii|(\al_{\ii'(1)}-\nobreak1)\allowbreak\cdots
\allowbreak(\al_{\ii'(q)}-\nobreak1)=\delta_{\ii, \ii'}f_\infty$ for
all $\tGm$-$q$-tuples~$\ii'$. With our choice of fundamental domain,
and using \eqref{aiidef}, we find for $\rho\in R$, $\dt\in \tDt$ and
$g=(x+\nobreak iy,\th)$ with $y\geq b$:
\begin{equation}
h_\ii(\rho\dt g) \= \sum_{l,m\geq0\,,\; l+m\leq q} a^\ii_{l,m}(\rho)\,
\ph^{(l,m)}(\dt)\, \ch(y)\, v_\infty(g)\,.
\end{equation}
Outside $\tGm D_\infty(b)$ the functions $f_\infty$, $h_\ii$ are zero.
With Lem.~\ref{lem-lo} we conclude that the function induced by
\begin{equation}
(\eta_\ii-h_\ii)(\rho\dt g)
= \sum_{l,m\geq0\,,\; l+m\leq q} a^\ii_{l,m}(\rho)
\chi(y)\bigl( v_\infty^{(l,m)}(\dt g) - \ph^{(l,m)}(\dt)\,v_\infty(g)
\bigr)
\end{equation}
is in $\map(\tG,\CC)^{\tGm,q}$. This implies that
\[ \eta_\ii \in \bigl( h_\ii + \ivt{\map(\tG,\CC)}q \bigr) \cap
\D_k(\ld)
= \ivt{\D_k(\ld)}{q+1} \,,\]
and behaves in the desired way under
$(\al_{\ii'(1)}-\nobreak1)\cdots\allowbreak(\al_{\ii'(q)}-\nobreak
1)$ for all $\tGm$-$q$-tuples~$\ii'$. Thus, we have proved that
$\D_k(\ld)$ is maximally perturbable.\smallskip

Everywhere in this proof we can replace $\A_k(b,\ld)$ by
$\A\hol_k(b)$, and $\D_k(\ld)$ by $\D\hol_k$. In that way we also
obtain the maximal perturbability of~$\D\hol_k$, thus completing the
proof of Proposition~\ref{prop-Dk-mp}.
\end{proof}

\subsubsection{Relations between the spaces $\C_k$ and $\D_k(\ld)$} By
Remark~\ref{eigenproperty}, for each $f\in \D_k(\ld)$ the support of
$(\om-\nobreak \ld)f$ is contained in some set $\tG_b$, hence
$(\om-\nobreak\ld)f\in \C_k$. So the differential operator~$\om-\ld$
maps $\D_k(\ld)
$ to~$ \C_k$. Since the operator $\om$ commutes with the action
of~$\tGm$, we have $(\om-\nobreak\ld)
\ivt{\D_k(\ld)}q\subset \ivt{\C_k}q$ for all $q\geq 1$. Similarly,
$\EE^- \ivt{(\D\hol_k)}q \subset \ivt{\C_{k-2}}q$ for all $q\geq 1$.
\begin{prop}\label{prop-DCs}Let $\ld\in \CC$ and $k\in 2\ZZ$. The
following maps are surjective:
\\
i. $\om -\ld:\D_k(\ld)^\tGm \rightarrow \C_k^\tGm$ and \\
ii. $\EE^- :
(\D\hol_k)^\tGm \rightarrow \C_{k-2}^\tGm$
\end{prop}
\begin{proof} \S\ref{i.} and \S\ref{ii.}.
\end{proof}

\begin{cor} \label{corol}For each $q\geq 1$ the maps $\om-\ld:\ 
\ivt{\D_k(\ld_s)}q \longrightarrow \ivt{\C_k}q$ and
$\EE^- :\ivt{(\D\hol_k)}q\rightarrow \ivt{\C_k}q$ are surjective.
\end{cor}
\begin{proof}Proposition~\ref{prop-DCs} gives the case $q=1$. The rows
in the following commutative diagram are exact by
Proposition~\ref{prop-Ck-mp} and~\ref{prop-Dk-mp}. See \eqref{tmldef}
for~$\ml_q$.
\begin{equation}
\label{DCdiag} \xymatrix@R=.5cm{ 0 \ar[r]& \ivt{\D_k(\ld)} q
\ar[r]\ar[d]^{\om-\ld}& \ivt{\D_k(\ld)}{q+1} \ar[r]^{ \ml_q}
\ar[d]^{\om-\ld}& \bigl( \D_k(\ld)^\tGm\bigr)^{n(\tGm,q)} \ar[r]
\ar[d]^{\om-\ld}
& 0\\
0\ar[r]& \ivt{\C_k}q \ar[r] \ar[d] & \ivt{\C_k}{q+1} \ar[r]^{ \ml_q}
\ar[d]
& \bigl( \C_k^\tGm\bigr)^{n(\tGm,q)} \ar[r] \ar[d] & 0\\
& 0 & \coker(\om-\ld) & 0 } \end{equation}
The third column is exact by Proposition~\ref{prop-DCs}. With the
exactness of the first column as induction hypothesis, we obtain the
vanishing of $\coker(\om-\nobreak\ld)$ and thus the surjectivity of
$\om-\ld: \ivt{\D_k(\ld)}{q+1} \rightarrow \ivt{\C_k}{q+1}$ by the
Snake Lemma.

The case of $\EE^-:\ivt{(\D\hol_k)}q\rightarrow \ivt{\C_k}q$ is
similar.
\end{proof}

\subsubsection{Proof of Proposition~\ref{prop-DCs}(i) }\label{i.}
We first note that the spaces $\D_k(\ld)^\tGm$ and $\C_k^\tGm$ are
invariant under~$\tZ$. Hence the weight~$k$ is strict and we are
dealing with functions on $G=\PSL_2(\RR)$. (See the first statement
in Theorem~\ref{thm-tAi}.)
We 
use the spectral theory of automorphic forms to prove
Proposition~\ref{prop-DCs}.

We work with the space of square integrable functions on
$\tGm\backslash \tG= \Gm\backslash G$ of strict weight~$k\in2\ZZ$,
where $G=\PSL_2(\RR)$. We can view the elements of the Hilbert space
$H_k= L^2(\tGm\backslash \tG)_k=L^2(\Gm\backslash G)_k $ as functions
$z\mapsto f(z,0)$ on~$\uhp$, transforming according to weight~$k$ as
indicated in~\eqref{actionb}. The inner product in $H_k$ is given by
\[ (f,f_1) \= \int_\fd f(z,0)\, \overline{f_1(z,0)}\,
\frac{dx\,dy}{y^2}\,.\]
Here $\fd$ can be any fundamental domain for $\Gm\backslash\uhp$. We
take it so that for each $b>A_\Gm$ it has a decomposition
\begin{equation}\label{fd1}
\fd \= \fd_b \sqcup \bigsqcup_{\k\in C}V_\k\,,\quad V_\k\= \bigl\{
g_\k(x+iy) \;:\; x_\k \leq x \leq x_\k+1\,,\; y\geq b\bigr\}\,,
\end{equation}
with $C$ a system of representatives of the $\Gm$-orbits of cusps, and
$x_\k\in\RR$ depending on $\fd$ and on the earlier choice of the $g_\k$. The
set $\fd_b$ has compact closure in~$\uhp$.

The differential operator $\om_k = -y^2\partial_y^2-y^2\partial_x^2+i
k y
\partial_x
$ in~\eqref{efC} determines a densely defined self-adjoint operator
$A_k$ in~$H_k$. The spectral theory of automorphic forms gives the
decomposition of this operator~$A_k$ in terms of Maass forms. One may
consult Chapters 4 and~7 in~\cite{Iw} for weight~$0$. For other
weights the proofs are almost completely similar. (See \cite{Ro}.)
There is a subspace $H_k^\discr$ with an at most countable orthonormal
basis $\{\ps_k^\ell\}$ of Maass forms, indexed by some subset
of~$\ZZ$. The $\ps_k^\ell$ are square integrable elements of the
space of Maass forms $E_k(\Gm,\ld^\ell)$ with
$\ld^\ell\geq \frac k2(1-\nobreak \frac k2)$. We denote the
eigenspace associated to $\ld$
(which is known to be finite-dimensional) by $H_k(\ld)$. If $k=0$ the
eigenvalue $0$ occurs with multiplicity one, corresponding to
constant functions, and all other $\ld^\ell$, if any, are positive.
If $k\neq 0$, then $H_k^\discr$ may be zero. If $k\geq 2$ and the space
$S_{\!k}(\Gm)$ of holomorphic cusp forms of weight~$k$ is non-zero, then there
are $\ps_k^\ell \in H^\discr_k$ of the form $\ps_k^\ell(z,0) = y^{k/2} \, h(z)$
with $h\in S_{\!k}(\Gm)$.
The corresponding
eigenvalues are $\ld^\ell = \frac k2\bigl(1-\nobreak \frac k2\bigr)$,
which is negative if $k\geq 4$. There may also be elements obtained
by differentiation of holomorphic cusp forms of weights between $2$
and $k-2$. Similarly, for negative $k$ there may be eigenfunctions
corresponding to antiholomorphic cusp forms.

The orthogonal complement $H_k^\cont$ of $H_k^\discr$ in~$H_k$ is
isomorphic to a sum of $\npar$ copies of $L^2\bigl(
(0,\infty),dt\bigr)$, where $\npar$ is the number of $\Gm$-orbits of
cusps. The spectral decomposition gives the \emph{Parseval formula}
\begin{equation}\label{parseval}
\bigl(f,f_1) \= \sum_\ell a_k^\ell(f)\,\overline { a_k^\ell(f_1) } +
\sum_\k \frac1{2\pi} \int_0^\infty e_k^\k(f;it)
\,\overline{e_k^\k(f_1;it)}\, dt\,,
\end{equation}
with $\k$ running through a set of representatives of the cuspidal
orbit. For each $f\in H_k$ we have $a_k^\ell(f) = \bigl( f,
\ps_k^\ell\bigr)$. If $f$ is sufficiently regular, then the functions
$e_k^\k(f;\cdot)$
are obtained by integration against the Eisenstein series $E_k^\k(it)$
at the cusp~$\k$.

The space $\C_k^\tGm$ is contained in $H_k$. For $f\in \C_k^\tGm$ the
functions $e_k^\k(f;\cdot)$ are given by
\[e_k^\k(f;s)\= \int_\fd f(z,0) \, \overline{E_k^\k(-\bar s;z)}\,
\frac{dx\,dy}{y^2} \= \int_\fd f(z,0) \,E_{-k}^\k(-s ; z)\,
\frac{dx\, dy}{y^2}\,,\]
for all $s$ at which the meromorphic continuation of the Eisenstein
series
$$E_k^\k(s;z):=\sum_{\gm \in \Gm_\k\backslash \Gm} \im(
g_\k^{-1}\gm z)^{\frac12+s}\, e^{-ik \arg(j(g_\k^{-1}\gm, z))}
$$ is holomorphic. In particular, $e_k^\k(f;s)$ is holomorphic at
points of the line~$i\RR$.

On the square integrable Maass forms and on the Eisenstein series the
self-adjoint operator $A_k$ is given by $\om_k$ in~\eqref{efC}. For
$f\in H_k$ in the domain of $A_k$, the self-adjointness of $A_k$
together with the eigenproperty of $\psi_k^l$ imply
$a_k^\ell(\om_k f)
= \ld^\ell a_k^\ell(f)$ and $e_k^\k(\om_k f;t) = \bigl( \frac14+t^2
\bigr)\,e_k^\k(f;t)$. This implies that the spectral data of elements
$f\in H_k$ such that $A_k^n f$ is well defined for all $n\in \NN$,
are quickly decreasing. The convergence in $L^2$-sense of the
Parseval formula in~\eqref{parseval} is very fast for functions of
this type, since the summands and integrands in the expansion are
those of $(A^n_k f, A^n_k f_1)$ divided by $(\ld^\ell)^n$,
respectively $(\frac14+\nobreak t^2)^n$ for each $n\in \NN$. (If
there is a term with $\ld^\ell=0$ we treat it separately; it does not
influence the convergence.)

The central point of the proof of Proposition~\ref{prop-DCs} is that
we transform the equation $(A_k-\nobreak \ld)f_1=f$ with unknown
$f_1\in H_k$ for a given $f\in \C_k^\tGm$ to the spectral
decomposition. Application of $A_k-\ld$ to $f\in \C_k^\tGm$ amounts
to multiplying $a_k^\ell(f)$ by $\ld^\ell-\ld$ and multiplying
$e_k^\k(f;t)$ by $\frac14+t^2-\ld $. This suggests the following
\begin{defn} Let $\ld\in \CC$. We denote by $ \C_k(\Gm,\ld)$ the space
of $f \in \C_k^\tGm$ such that the following conditions are
satisfied.
\begin{enumerate}
\item[i)] $a_k^\ell(f)=0$ if $\ld_{\ell}=\ld$,
\item[ii)] $e_k^\k(f;it_\ld)=0$ for all~$\k$, if $\ld=\frac14+t_\ld^2$
($t_\ld\in \RR\setminus\{0\})$ ,
\item[iii)] for all $\k$, the map $t\mapsto e_k^\k(f;it)$ has a double
zero at~$t=0$ if $\ld=\frac14$.
\end{enumerate}
\end{defn}
Note that, for each $\ld$, the conditions i), ii), iii) impose
finitely many linear conditions, so $\C_k(\Gm,\ld)$ has finite
codimension in~$\C_k^\tGm$. If $\ld$ is not in the spectrum of~$A_k$,
then $\C_k(\Gm,\ld)$ is equal to~$\C_k^\tGm$.
\smallskip

Case I: $f \in \C_k(\tGm, \lambda)$. In this case, if we have the
spectral decomposition
\begin{equation}\label{spectral}
f \= \sum_\ell a_k^\ell(f) \psi_k^{\ell}
+\sum_\k \frac1{2\pi} \int_0^\infty e_k^\k(f;it)
\overline{E_k^{\k}(it; -)} dt,
\end{equation}
then a solution of $(A_k-\ld)f_1=f$ is given by
\begin{equation} \label{solution}
f_1 := \sum_\ell \frac{a_k^\ell(f)}{\ld^\ell-\ld} \psi_k^{\ell}
+ \sum_\k \frac1{2\pi} \int_0^\infty \frac{e_k^\k(f;it)}
{\frac14+t^2-\ld } \overline{E_k^{\k}(it; -)} dt.
\end{equation}
If $\ld $ is not in the spectrum of~$A_k$, then the convergence of
this $L^2$-expansion is better than that in~\eqref{spectral}, since
$\lambda^{\ell} \to \infty$ and hence the denominators improve the
convergence. If $\ld$ is in the spectrum of~$A_k$, condition~i)
ensures that the $a^\ell_k(f)$ with $\ld^\ell=\ld$ vanish, and that,
by the other conditions, the simple or double zero of
$t\mapsto \frac14+t^2-\ld$ at $t=t_\ld$ is canceled by the zeros at
$s=it_\ld$ of the holomorphic functions $s\mapsto e^\k_\k(f;s)$.
The same reasoning shows that the obtained $f_1$ is in $H_k$
and, in fact, in the domain of $A_k$. We have $(A_k-\nobreak\ld)
f_1=f$. Therefore the relation $(\om-\nobreak\ld)f_1=f$ holds in
distribution sense. There is $b\geq A_\Gm$ such that the support
of~$f$ is contained in $\tG_b$. So on each $D_\k(b)$ we have
$(\om-\nobreak\ld)f_1=0$. Since $\om_k$ determines an elliptic
differential operator on~$\uhp$, elliptic regularity implies that
$(\om-\nobreak\ld)f_1=0$ holds as a relation for real-analytic
functions on each $D_\k(b)$. Further, the square integrability
implies that $f_1$ must have less than exponential growth at the
cusps
and hence it is an element of $\D_k(\ld)^\tGm$. We have shown:
\begin{lem}\label{lem-CG}For each $\ld\in \CC$, the space
$\C_k(\Gm,\ld)$ is contained in $(\om-\nobreak\ld)\, \D_k(\ld)^\tGm$.
\end{lem}

Case II: $f \in \C_k^\tGm \smallsetminus \C_k(\Gm,\ld)$ for $\ld$ in
the spectrum of~$A_k$. The following result enables us to pick
representatives $h$ of $\C_k^\tGm / \C_k(\Gm,\ld)$ for which we can
solve $(\om-\nobreak\ld) f_1=h$ directly. This procedure can be
carried out by singling out one cusp~$\k$, which we fix for the proof
of Case~II.
\begin{lem} \label{lem-hn} Let $\k$ be the cusp that we keep fixed.
Suppose that $\ld$ is in the spectrum of~$A_k$. Then there is a
finite set $X\subset \ZZ$ such that, for each $n\in X$, there exist
$h_n \in \C_k^\tGm$ of the form
\begin{equation}
h_n\bigl( \gm \tilde g_\k (z,\th) \bigr) =\begin{cases}
 e^{2\pi i n x}\, \ch_n(y) \, e^{ik\th} & \text{on} \, \, \tGm
D_\k(A_\Gm)
\\
0 & \text{elsewhere}
\end{cases}
\end{equation}
for some $\ch_n \in C_c^\infty(A_\Gm,\infty)$, such that
$\{h_n+C_k(\Gm, \lambda)\}_n$ spans $\C_k^\tGm / \C_k(\Gm,\ld)$.
\end{lem}
If we can solve $(A_k-\nobreak\ld) f_1= h_n$ in another way for all
$n\in X$, this lemma enables us to reduce the proof of
Proposition~\ref{prop-DCs}(i)
to Lemma~\ref{lem-CG}.
\begin{proof}[Proof of Lemma \ref{lem-hn}] We shall examine each of
the three cases for the eigenvalues of $A_k$ on $H_k$ separately:
\begin{itemize}
\item $\ld=\frac14-s^2\not\in \bigl[\frac14,\infty\bigr)$. Assume
$s > 0$. There are finitely many indices $\ell_1,\ldots,\ell_m$ such
that $\ld_{\ell_j}=\ld$. The $\ps_k^{\ell_j}$ form a basis of
$\ker\left(A_k-\ld\right)$. Each of these $m$ linearly independent
square integrable automorphic forms is given by its Fourier expansion
at the fixed cusp~$\k$. By Proposition~\ref{prop-Fe}, the Fourier
terms of non-zero order are multiples of $\om_k(n,s)$. The Fourier
term of order zero is a multiple of $y^{\frac12-s}\, e^{ik\th}$. We
choose a set $X$ of $m$ elements in $\ZZ$ such that the
$m\times m$-matrix whose columns are the $n$-th Fourier coefficients
of $\ps_k^{\ell_j}$
($1\leq j \leq m$) with $n \in X$ is invertible. We choose the
$\ch_n\in C_c^\infty$, $n\in X$, in the statement of the lemma, so
that $\int_{A_{\Gm}}^\infty \ch_n(y)
\,\overline{\om_k(n,s)(iy,0)}\, \frac {dy}{y^2}\neq 0$, respectively
$\int_{A_{\Gm}}^\infty \ch_n(y)
\,\overline{y^{\frac12 -s}}\, \frac {dy}{y^2}\neq 0$. Consider the
linear form on the space $A_k^2(\ld)$ of square integrable
automorphic forms with eigenvalue $\ld$ given by
\begin{align*}
\ps \mapsto &\bigl(h_n, \ps \bigr)\= \int_{\fd} h_n(z, 0)
\overline{\psi(z, 0)} \frac{dx dy}{y^2}
\\
&\= \int_{A_{\Gm}}^{\infty} \int_{-1/2}^{1/2} \chi_n(y) e^{2 \pi i n
x} \bar a_0 y^{1/2-\bar s}\, \frac{dx dy}{y^2}
+ \sum_{m \ne 0} \bar a_m \int_{A_{\Gm}}^{\infty} \int_{-1/2}^{1/2}
\chi_n(y) e^{2 \pi i n x} \overline{\om_k(m, s)(iy, 0)} \frac{dx
dy}{y^2}\,.
\end{align*}
This depends only on the Fourier coefficient of $\ps$ of order $n$ in
the expansion at~$\k$. Therefore, the $m\times m$-matrix with the
scalar product $\bigl(h_n,\ps^{\ell_j}_k\bigr)$ at position $(j,n)$
is invertible. (Here $j$ runs from $1$ to~$m$, and $n$ runs
through~$X$.) Hence there are complex numbers $b_{j,p}$ (with
$1\leq j \leq m$, $p\in X$) such that
$\sum_{n\in X} b_{j,n} \, \bigl(h_n,\ps^{\ell_{j'}}_k\bigr) =
\dt_{j,j'}$. Setting, for $f \in \C_k^\tGm$, $c_n(f) = 
\sum_{j'=1}^m \bigl( f, \ps_k^{\ell_{j'}} \bigr) \, b_{j ',n}$, we
obtain for $1\leq j \leq  m$:
\[ \sum_{n\in X} c_n(f) \, \bigl( h_n, \ps^{\ell_j}\bigr) \= \bigl(
f,\ps_k^{\ell_{j}}\bigr)\,. \]
So $f- \sum_n c_n(f)\, h_n$ is indeed in $\C_k(\Gm,\ld)$.
\item $\ld=\frac14+t^2$, $t\in \RR\setminus\{0\}$. A basis of
$\ker\left(A_k-\ld\right)$ in this case consists of Eisenstein series
$E^\nu_k(it, \cdot)$ ($\nu \in C$) and possibly cusp forms
$\ps_k^{\ell_j}$ with $\ld_{\ell_j}=\ld$. The proof of the previous
case can be applied with the obvious adjustments (e.g. replacing
scalar products by integrals for the terms corresponding to
$E^\nu_k$) to give the result. The only essential modification is
that we have to use the space $A^\ast_k(\ld)$ of automorphic forms
with polynomial growth and eigenvalue~$\ld$ in place of $A_k^2(\ld)$
because the Eisenstein series are not square integrable. This can be
done because (conjugates of) elements of $A^\ast_k(\ld)$
appear only integrated against elements of $\C_k^\tGm$ which have
compact support modulo~$\tGm$.
\item $\ld=\frac14$. Now we have the condition that $e_k^{\k}(f-\sum_n
h_n; it)$ should have a double zero at $t=0$ or, equivalently, that
the first two terms of the Taylor expansion at $s=0$ should vanish.
Since the first two Taylor terms of $E_k^\k(-; z)$ are linearly
independent from the other functions in $A^\ast_k(1/4)$, a choice of
$\ch_n$ with the desired properties is again possible.\qedhere
\end{itemize}
\end{proof}

Now we turn to the task to solve $(\om-\nobreak\ld_s) f_1 = h_n$ with
$f_1\in \D_k(\ld)^\tGm$ for $h_n$ as in Lemma~\ref{lem-hn}. We aim at
$f_1$ with support in $\tGm D_\k(A_\Gm)$. Writing
$f_1\bigl( \tilde g_\k (z,\th) \bigr) = e^{2\pi
i n x} \, h(y)\, e^{ik\th}$, the differential equation $(\om-\nobreak
\ld)f_1=h_n$ becomes
\[ -y^2 h''(y) + \bigl(4\pi^2n^2 y^2-2\pi n k y -\frac14+s^2 \bigr)
\,h(y) \= \ch_n(y)\,.\]
(Compare~\eqref{diff-eq-Wh}.)
This ordinary differential equation is regular on $y\geq A_\Gm$. It
has a unique solution for the initial conditions
$h(A_\Gm)=h'(A_\Gm)=0$. It is zero below the support of~$\ch_n$.
Since $\ch_n$ has compact support, the function $h$ thus obtained is
a solution of the homogeneous equation \eqref{diff-eq-Wh} on
$(b,\infty)$ for some $b>A_\Gm$ depending on~$\supp(\ch_n)$. Thus we
see that $(z,\th) \mapsto f_1\bigl( \tilde g_\k(z,\th)
\bigr)$ is an element of $\W_k(n,s)$. Hence it may have exponential
growth of order $e^{(2\pi|n|+\dt)y}$. This is the point where the
need to work with exponentially growing functions arises.

We extend $f_1$ by $\tGm$-invariance, and check that it is an element
of~$\D_k(\ld_s)$. This completes the proof of the first statement in
Proposition~\ref{prop-DCs}.\bigskip

\subsubsection{Proof of Proposition~\ref{prop-DCs} ii} \label{ii.}
For the surjectivity of $\EE^- : (\D\hol_k)^\tGm\rightarrow
\C_{k-2}^\tGm$ we first note that, on an eigenfunction of $\om$ in
weight~$k-2$ with eigenvalue $\ld$ the operator $\EE^-_k \EE^+_{k-2}$
acts as multiplication by $-4 \bigl( \ld - \nobreak \frac k2
+\nobreak\frac {k^2}4\bigr)$. See~\eqref{Lie-diff}. We will use
$\EE^+_{k-2}$ to ``invert"~$\EE^-_k$.

Let $H^a_{k-2}$ denote the kernel of $\EE_{k-2}^+$ in
$H_{k-2}^\discr$. It is finitely dimensional and it contains the
constant functions if $k=2$, and the functions corresponding to
antiholomorphic cusp forms if $k\leq 0$.

On the orthogonal complement of $H^a_{k-2}$ in $H_{k 
-2}^\discr$ the factor $-4 \bigl( \ld -\nobreak \frac k2
+\nobreak\frac {k^2}4\bigr)$ is negative and stays away from $0$ for
all $\ld$ in the spectrum of $A_k$. Likewise, we denote by $H^h_k$
the finite dimensional kernel of $\EE_k^-$ in~$H^\discr_k$. Its
elements correspond to square integrable holomorphic automorphic
forms of weight~$k$.

Let $\bigl( \ps_{k-2}^\ell\bigr)_\ell$ be an orthonormal basis of the
orthogonal complement $H_{k-2}^\discr \ominus H^a_{k-2}$ consisting
of eigenfunctions of $\om_{k-2}$ with eigenvalue~$\ld^\ell$. The
relation $(\EE^-_k v_1, v_2)=-(v_1,\EE^+_{k-2} v_2)$ for suitably
differentiable elements of $H_k$ and $H_{k-2}$ (see Lemma 6.1.4 of
\cite{Br86}) implies that $\EE^-_k(H_k^\discr \ominus H_k^h)
\subset H_{k-2}^\discr \ominus H_{k-2}^a$ and hence
$(\ps_k^\ell)_{\ell}$ with $\ps_k^\ell=
\frac1{\sqrt{4\ld^\ell-2k+k^2}}\EE_{k-2}^+\ps_{k-2}^\ell$ is an
orthonormal system spanning $H_k^\discr \ominus H_k^h$.

For a given $f\in \C_{k-2}^\tGm$ orthogonal to $H^a_{k-2}$ we set
\begin{multline*} f_1 \;:=\; -\sum_\ell
\frac{a_{k-2}^\ell(f)}{\sqrt{4\ld^\ell-2k+k^2}}
\frac{\EE^+_{k-2}\ps_{k-2}^{\ell}}{\sqrt{4\ld^\ell-2k+k^2}} \;-\;
\sum_\k \frac1{2\pi} \int_0^\infty \frac{e_{k-2}^\k(f;it)}
{\sqrt{4t^2+(k-1)^2}} \frac{\EE^+_{k-2}E_{k-2}^{\k}(it;
-)}{\sqrt{4\ld^\ell-2k+k^2}} dt\,.
\end{multline*}
We have $f_1 \in H_k\ominus H^h_k$ and $\EE^-f_1=f$. 
A reasoning as in the previous case shows that
$f_1\in \D\hol_k(\ld)^\tGm$.

So we have solved the problem for a subspace of $\C_{k-2}^\tGm$ with
finite codimension. A general element of $\C_{k-2}^\tGm$will not be
orthogonal to $H^a_{k-2}$. We proceed as in the first case in the
proof of Lemma~\ref{lem-CG}. Instead of $\ps_k^{\ell_j}$ we now use
an orthogonal basis of $H^a_{k-2}$, and form functions $h_n$ as in
Lemma~\ref{lem-CG}, corresponding to a set $X$ of Fourier term orders
such that elements of $H^a_{k-2}$ are determined by the Fourier
coefficients in~$X$. Solving $\EE_k^- f_1 = h_n$ leads to the
differential equation
\[ (-2iy\partial_x+2y\partial_y-k) e^{2\pi i n x} \ph(y) \=
\ch(y)\,,\quad \ph(y_0)\=\ph'(y_0)=0\,,\]
with which we proceed as in the previous case.

This establishes the surjectivity of $\EE^-:(\D\hol_k)^\tGm\rightarrow
\C_{k-2}^\tGm$ in Proposition~\ref{prop-DCs}.

\subsection{Higher order invariants and Maass forms}\label{sect-fp} We
now will derive the main results of this paper,
Theorems~\ref{thm-tMf-mp} and~\ref{thm-thol-mp}, from the following
result:
\begin{prop} The $\tGm$-modules
\begin{equation}
\tilde\A_k'(\ld):=\ker\bigl( \om-\ld: \D_k(\ld) \longrightarrow
\C_k\bigr)\end{equation}
and
\begin{equation}
\Hol'_k :=\ker\bigl( \EE^-:\D\hol_k\longrightarrow \C_{k-2}\bigr)
\end{equation}
are maximally perturbable.
\end{prop}
\begin{proof} \ We have the following extension of the commutative
diagram~\eqref{DCdiag}:
\begin{equation} \xymatrix@R=.5cm{
0\ar[r]& \ivt{\tilde \A_k'(\ld)}q \ar[r] \ar[d]& \ivt{\tilde
\A_k'(\ld)}{q+1} \ar[r]^{\ml_q} \ar[d]& \bigl( \tilde
\A_k'(\ld)^\tGm\bigr)^{n(\tGm,q)} \ar[d]& \\
0 \ar[r]& \ivt{\D_k(\ld)} q \ar[r] \ar[d]^{\om-\ld}&
\ivt{\D_k(\ld)}{q+1} \ar[r]^{\ml_q} \ar[d]^{\om-\ld}& \bigl(
\D_k(\ld)^\tGm\bigr)^{n(\tGm,q)} \ar[r] \ar[d]^{\om-\ld}
& 0\\
0\ar[r]& \ivt{\C_k}q \ar[r] \ar[d] & \ivt{\C_k}{q+1} \ar[r]^{\ml_q}
\ar[d]& \bigl( \C_k^\tGm\bigr)^{n(\tGm,q)} \ar[r] \ar[d] & 0
\\
& 0 & 0 & 0 } \end{equation}
The exactness of the columns follows from the definition of $\tilde
\A_k'(\ld)$, \eqref{iv-I}, the left-exactness of the functor
$\hom_{\CC[\Gm]}(I^q\backslash \CC[\Gm], \text{--})$ and Corol.
\ref{corol}. Propositions \ref{prop-smk-mp} and~\ref{prop-Dk-mp}
imply that the second and third row are exact. 
The Snake Lemma then implies that the first row is exact and that
$\ml_q : \ivt{\tilde\A_k'(\ld)}{q+1} \rightarrow \bigl( \tilde
\A_k'(\ld)^\tGm\bigr)^{n(\tGm,q)}$ is surjective.

Replacing in this diagram the space $\tilde \A_l'(\ld)$ by $\Hol'_k$
and the map $\om-\ld$ by $\EE^-$, we obtain the maximal
perturbability of~$\Hol'_k$.
\end{proof}

\begin{proof}[Proof of Theorems \ref{thm-tMf-mp}
and~\ref{thm-thol-mp}] The $\tGm$-module $\tilde\A_k'(\ld)$ is
contained in $\tilde\A_k(\tGm,\ld)$. See Definition~\ref{defTa}. It
is a smaller space than $\tilde \A_k(\tGm,\ld)$ since elements of
$\D_k(\ld)$ have a special structure near the cusps. With
\eqref{inclusion}, $\tilde \A_k(\tGm,\ld)^\tGm$ is a subspace of
$\tilde{\A}_k'(\ld)^\tGm$. Therefore $\tilde 
\A_k(\tGm,\ld)^\tGm =\tilde{\A}_k'(\ld)^\tGm$ and thus
\[\xymatrix@R=.5cm{ 0\ar[r]
& \ivt{\tilde \A_k'(\ld)}q \ar[r] \ar[d] & \ivt{\tilde
\A_k'(\ld)}{q+1} \ar[r] \ar[d] & \bigl( \tilde
\A_k'(\ld)^\tGm\bigr)^{n(\tGm,n)} \ar[r] \ar@{=}[d]
&0
\\
0\ar[r]& \ivt{\tilde{\A}_k(\tGm,\ld)}q \ar[r] & \ivt{\tilde
\A_k(\tGm,\ld)}{q+1} \ar[r] & \bigl( \tilde
\A_k(\tGm,\ld)^\tGm\bigr)^{n(\tGm,q)} }\]
with exact rows. Induction with respect to~$q$ and the Snake Lemma
show that $\ivt{ \tilde \A_k(\tGm,\ld)}q$ is equal to
$\ivt{\tilde \A_k'(\ld)}q$ for all~$q$. Hence the space
$\tilde{\A}_k(\tGm,\ld)$ is maximally perturbable.

The proof of Theorem~\ref{thm-thol-mp} is completely similar.
\end{proof}

\appendix
\section{Partition of unity} The following technical lemma gives
partitions of unity that are adapted to $\tGm\backslash\tG$ and
$\Gm\backslash\uhp$.
\begin{lem}\label{lem-pu}
\begin{enumerate}
\item[i)]For a given cofinite discrete $\tGm\subset\tG$ containing
$\tZ$ there are $\ps\in C^\infty(\tG)$ such that
\begin{enumerate}
\item[a)] $\ps$ is a bounded function.
\item[b)] There is $N\in \NN$ such that for each $g\in \tG$ the number
of $\gm\in \tGm$ with $\ps(\gm^{-1}g) \neq 0$ is bounded by~$N$.
\item[c)] $\sum_{\gm\in \tGm}\ps(\gm^{-1}g)=1$ for all $g\in \tG$.
\end{enumerate}
\item[ii)] For a given cofinite discrete $\Gm\subset\PSL_2(\RR)$ there
are $\ps_0\in C^\infty(\uhp)$ such that
\begin{enumerate}
\item[a)] $\ps_0$ is bounded.
\item[b)] There is $N\in \NN$ such that for each $z\in\uhp$ the number
of $\gm\in \Gm$ with $\ps_0(\gm^{-1}z)\neq0$ is bounded by~$N$.
\item[c)] $\sum_{\gm\in \Gm}\ps_0(\gm^{-1}z) =1$ for all $z\in \uhp$.
\end{enumerate}
\end{enumerate}
\end{lem}
\begin{proof}
We fix a strict fundamental domain $\fd_\uhp$ for
$\tGm\backslash \uhp$ of the following form, based on the choice of a
real number $a>A_\Gm$, as in~\S\ref{sect-spC}. The set  
$\fd_\uhp$ is bounded
by finitely many geodesic segments and half-lines such that
\begin{equation}\label{fdecomp}
\begin{aligned} \fd_\uhp &\= C_a \sqcup \bigsqcup_\k V_\k(a)\,,\\
V_\k(a)&\= \bigl\{ g_\k(x+iy) \;:\; y \geq a \,,\; x_\k \leq x <
x_\k+1\bigr\}\,,
\end{aligned}
\end{equation}
where $C_a$ is relatively compact in~$\uhp$, and is contained in the
image of $\tG_a$ under the projection $\tG\rightarrow\uhp$. The
disjoint union is over the set $C$ of cusps $\k$ in the closure
$\bar \fd_\uhp$ of $\fd_\uhp$ in $\uhp\cup\partial
\uhp$. We take $\fd_\uhp$ such that $C$ forms a system of
representatives for the $\tGm$-orbits of cusps. By taking the
parameter~$a$ sufficiently large we arrange that all orbits of
elliptic fixed points intersect $\fd_\uhp$ in~$C_a$. These points are
necessarily on the boundary of $\fd_\uhp$.

We take a strict fundamental domain for $\tGm\backslash\tG$ of the
form
\[ \fd \= \bigl\{ (z,\th)\;:\; z\in \fd_\uhp\,,\; \th \in
[0,\pi/v(z))\bigr\}\,,\]
where $v(z)\in \NN$ is the order of the subgroup $\Gm_z$ fixing $z$,
or equivalently $\tGm_z $ is conjugate in~$\tG$ to the group $\{
k(n\pi/v(z))\;:\; n\in \ZZ\}$. So $v(z)$ is in general equal to~$1$,
and only larger if $z$ is an elliptic fixed point of~$\Gm$.\smallskip

i. We first define a function on $\tG$ satisfying a) and~c), and a
variant of~b).

Let $\om:\tG\rightarrow \{0,1\}$ be the characteristic function
of~$\fd$. It satisfies conditions a)--c) in part~i) of the lemma, but
is not smooth. To make it smooth we convolve it with a function
$\ps\in C_c^\infty(\tG)$ with $\ps\geq 0$ such that
$\int_{\tG} \ps(g)\, dg=1$ for a choice $dg$ of a Haar measure
on~$\tG$ and such that $\supp(\ps)$ is a compact neighborhood of the
unit element in~$\tG$.

Since $\om$ is measurable, the integral
\[ \ph_0(g) \= \int_{\tG} \om(g_1)\, \ps(g_1^{-1}g)\, dg_1 \= \int_\tG
\om(gg_1^{-1})\, \ps(g_1)\, dg_1 \]
defines a smooth function $\ph_0$ on~$\tG$ with values in $[0,1]$ and
with support contained in the neighborhood $\fd\cdot\supp(\ps)$
(multiplication in $\tG$) of~$\fd$. From the second form of the
convolution integral we see that
$\sum_{\gm\in \tGm} \ph_0(\gm^{-1}g)=1$ for all~$g\in \tG$. This
smooth function $\ph_0$ satisfies conditions a) and~c) in part~i) of
the lemma. Condition~b)
is not satisfied, since although the support of $\ph_0$ is contained
in a neighborhood of $\fd$ of the form $ \fd\,\supp(\ps) $, this
neighborhood may meet near the cusps infinitely many
$\tGm$-translates of~$\fd$. We will construct two
functions, one ``away from the cusps" and another ``close to the
cusps" satisfying all conditions a), b), c) on overlapping regions. A
suitable combination of these two functions will produce the sought
function on $\tG$.

$\bullet$ The first function is simply the restriction of $\ph_0$ to
$\tilde G_b$ for any $b\geq a$. We will show that this function
satisfies condition~b) (and thus all conditions). First we note that
the projections $p_1:\tG\rightarrow \uhp$ and $p_2:\tG\rightarrow\RR$
given by $p_1(z,\th)=z$ and $p_2(z,\th)=\th$ are continuous. Next we
note that $\fd\,\supp(\ps) \cap \tG_b$ is contained in a compact set,
and hence has compact image in $\uhp$ under~$p_1$. So
\[ p_1\bigl(\fd\,\supp(\ps) \cap \tG_b \bigr) \,\subset\,
\bigsqcup_{\dt\in E} \dt \,\fd_\uhp \]
for some finite subset $E$ of~$\Gm$.

Fix a $g\in \tG_b$. We will show that there is a finite number
(independent of $g$) of $\gm\in \tGm$ with $\ph_0(\gm \, g)\neq  0$.
Indeed, for each such $\gm$ we have
$\gm g \in \fd\, \supp(\ps)\cap \tilde G_b$, hence
$p_1(\gm g) = \pr(\gm)\,
p_1(g) \in \bigsqcup_{\dt\in E} \dt\, \fd_\uhp$. This leaves finitely
many possibilities for the image $\pr(\gm)$:
\[ \pr(\gm) = \dt\dt_0^{-1}\qquad\text{ with }\dt\in E\,.\]
for some $\dt_0 \in \Gm$. We conclude that $\gm = 
\widetilde{\dt\dt_0^{-1}} k(\pi m)$ with $m\in \ZZ$.

On the other hand, the image $p_2\bigl( 
\fd\,\supp(\ps)\cap\tG_b\bigr)$ is contained in a compact set, hence
it is contained in a set $[-B,B]\subset \RR$. For the
$\gm=\widetilde{\dt\dt_0^{-1}}k(\pi m)$ with $\ph_0(\gm g)\neq0$ we
conclude from~\eqref{k(th)} that
$p_{2} \bigl( \widetilde {\dt\dt_0^{-1}} k(\pi m) g \bigr)=
p_2 \bigl( \widetilde {\dt\dt_0^{-1}} g \bigr)+m \pi$. This leaves
only finitely many possibilities for the integer~$m$. This shows that
condition~b) is satisfied by the restriction of $\ph_0$ to $\tG_b$
($b\geq a$).
\medskip

$\bullet$ We now start the construction of another function $\ph_1$
with the desired properties near the cusps. We take a compactly
supported smooth partition $\bt$ of unity for $\RR/\ZZ$, {\sl i.e.},
$\bt \in C_c^\infty(\RR)$ with values in $[0,1]$ such that
$\sum_{k\in \ZZ} \bt(x+\nobreak k)=1$ for all $x\in \RR$.
(For instance take a smooth function $\ups$ in $C^\infty(\RR)$ with
value $0$ on a neighborhood of~$0$ and value $1$ on a neighborhood
of~$\frac12$. Then
\[ \bt(x) \= \begin{cases}
0&\text{ if }x<0\,,\\
\ups(x)&\text{ if }0\leq x< \frac12\,,\\
1&\text{ if }\frac12\leq x<1\,,\\
1-\ups(x-1)&\text{ if }1\leq x<\frac 32\,,\\
0&\text{ if }x\geq \frac 32\,.
\end{cases}
\]
defines such a partition of unity.)
We define a function $\ph_1$ on $\tG$ in the following way.
\begin{align*}
\ph_1\bigl( \tilde g_\k (z,\th) \bigr)&\= \bt(x-x_\k)\,
\bt(\th/\pi)\quad\text{ if }\k \in C,\, y>a\,,\\
\ph_1&\=0\quad\text{elsewhere}\,,
\end{align*}
with $x_\k$ as in~\eqref{fdecomp}.

The function $\ph_1$ is smooth on $\tG\setminus \tG_a$ and bounded
there. By the definition and \eqref{commut}, it is clear that the
only $\tGm$-translates of $\fd$ intersecting the support of $\ph_1$
are the $\tGm_\k$-translates. The definition of $\beta$ then implies
that at most four $\tGm_\k$-translates can intersect Supp$(\ph_1)$,
implying (ii). Likewise, the definitions of $\ph_1$, $\beta$ and
\eqref{commut} imply that $\sum_{\gm\in \tGm} \ph_1(\gm^{-1}g)=1$ for
$g\in \tG\setminus \tG_a$.

We choose a (bounded) function $\ch\in C^\infty(\tGm\backslash\tG)$
equal to $0$ on $\tG_a$ and equal to $1$ on $\tG\setminus\tG_{a+1}$.
Put
\[ \ps = \ch\cdot \ph_1 + (1-\ch) \cdot \ph_0\,,\]
where $\ph_0$ is as constructed above with $b$ equal to $a+1$. Since
$\ch$ vanishes on $\tG_a$ the product $\ch\cdot \ph_1$ is smooth
on~$\tG$. Similarly, $(1-\nobreak\ch)\cdot \ph_0$ is smooth. So
$\ps\in C^\infty(\tG)$. Conditions a)--c)
are easily checked to hold for~$\ps$.

\medskip

ii. We turn to $\Gm=\tGm/\tZ$ and start with $\ps$ as in part~i). The
sum $\ps_1(z,\th)= \sum_{m\in \ZZ} \ps(z,\th-\nobreak m\pi)$ is
locally finite and defines a smooth function with values in $[0,1]$
that is invariant under left translation by elements of~$\tZ$. So
$\ps_1\bigl(\gm^{-1}(z,\th)\bigr)=\ps_1\bigl(\tgm^{-1}(z,\th)\bigr)$
is well defined for $\gm\in \Gm$, and
\[\sum_{\gm\in \Gm} \ps_{1}\bigl(\gm^{-1}(z,\th)
\bigr)\= \sum_{\gm\in \Gm} \sum_{m\in \ZZ} \ps\bigl( (\tilde\gm k(m
\pi))^{-1} (z,\th) \bigr)
\=1 \]
for all $(z,\th)$. Since the support of $\ps$ meets only finitely many
$\tGm$ translates of $\fd\subset \fd_\uhp \times [0,\pi)$, the
support $ \supp(\ps)\cdot\tZ $
of $\ps_1$ meets only finitely many $\tGm$-translates of
$\fd_\uhp\times \RR$. Set
$$\psi_0(z)=\frac{1}{\pi} \int_0^{\pi} \psi_1(z, \th) d \th.$$
It clearly satisfies (i). For condition~c) we note that
\begin{eqnarray*} \sum_{\gm\in \Gm} \ps_0(\gm\, z) &=& \frac1\pi
\sum_{\gm \in \Gm} \int_0^{\pi} \ps_1\bigl(\gm z, \th \bigr)\, d\th
\\
&=& \frac1\pi \sum_{\gm \in \Gm} \int_{\arg(j(\gm, z))}^{\pi+
\arg(j(\gm, z))} \ps_1\bigl(\gm z, \th \bigr)\, d\th
\qquad \text{by the $\pi$-periodicity of $\ps_1$}\\
&=& \frac1\pi \sum_{\gm \in \Gm} \int_0^{\pi} \ps_1\bigl(\gm (z, \th)
\bigr)\, d\th \\
&=&1
\end{eqnarray*}
The support of $\ps_0$ is contained in the image $p_1\bigl(
\supp(\ps_1)\bigr)\subset\uhp$. Since $\supp(\ps_1)$ is contained in
finitely many $\tGm$-translates of $\fd_\uhp\times\RR$, we conclude
that condition~b) is satisfied as well.
\end{proof}

\section{Index of commonly used notation}
\begin{multicols}{3}
\begin{trivlist}
\item $a(y)$ \ \hfill \S\ref{Ucg}
\item $a_k^\ell(f)$ \ \hfill \eqref{parseval}
\item $\alpha$ \ \hfill \S\ref{sect-cdsg}
\item $\alpha_i$ \ \hfill \S\ref{cgen}
\item $\bb(\ii)$ \ \hfill \eqref{bi-def0}, \eqref{bi-def}
\item $C^\infty_k(\tG)$ \ \hfill \eqref{C_k}
\item $\C_k$ \ \hfill Defn. \ref{defn-ck}
\item $D_{\k}(a)$ \ \hfill \eqref{Dk-def}
\item $\D_k(\ld),$ $\D\hol_k$ \ \hfill Defn. \ref{Ddef}
\item $\e_i$ \ \hfill \S\ref{cgen}
\item $E_i$ \ \hfill \S\ref{canon-gener}
\item $E_k(\Gm, \ld)$ \ \hfill Defn.~\ref{MAkdef}
\item $E_k\hol(\Gm,\ld_k)$ \ \hfill \S\ref{haf}
\item $\EE^{\pm}$ \ \hfill \S\ref{Lie}
\item $e_k^\k(f;it)$ \ \hfill \eqref{parseval}
\item $\A_k(y_0,\ld), \A\hol_k(y_0)$ \ \hfill Defn.~\ref{Xdef}
\item $\tilde \A_k(\tGm,\ld)$ \ \hfill Defn.~\ref{defTa}
\item $\tilde E_r(\tGm,\ch,\ld)$ \ \hfill Defn.~\ref{defTa}
\item $E_k(\Gm, \ld)$ \ \hfill Defn.~\ref{MAkdef}
\item $\ff_\ii$ \ \hfill \eqref{fidef}
\item $g_\k$ \ \hfill \eqref{scal}
\item $\tG$ \ \hfill Defn.~\ref{ucg}
\item $\tG_a$ \ \hfill \eqref{tGa-def}
\item $h_{\ii}$\ \hfill \eqref{hidef-gen}, \eqref{hdef-pu}
\item $H_i$ \ \hfill \S\ref{canon-gener}
\item $\HH$ \ \hfill \S\ref{Lie}
\item $h^\mm_k(n,s)$ \ \hfill \eqref{ztmm}
\item $\Hol_k(\tGm), \Hol^p_k(\tGm), \Hol^c_k(\tGm)$ \ \hfill Defn.
\ref{holom} \
\item $\eta_i$ \ \hfill \S\ref{cgen}
\item $\eta_r(n;z,\th)$ \ \hfill \eqref{eta1}
\item $\eta_k(n)$ \ \hfill \eqref{eta2}
\item $\eta_k^\mm(n;z,\th)$ \ \hfill \eqref{eta3}
\item $k(\th)$ \ \hfill \S\ref{Ucg}
\item $\k_i$ \ \hfill \S\ref{cgen}
\item $L_k$ \ \hfill \eqref{Lkdef}
\item $L(z,\th)$ \ \hfill \eqref{Ldef}
\item $\ml_q \quad$ \ \hfill \eqref{exs-iv}
\item $\M_k(\Gm, \ld)$ \ \hfill Defn.~\ref{MAkdef}
\item $M_k(\Gm, \ld)$ \ \hfill Defn.~\ref{MAkdef}
\item $\M_k\hol(\Gm,\ld_k)$ \ \hfill \S\ref{haf}
\item $\mu_f$ \ \hfill \eqref{pert}
\item $\nell, \npar$ \ \hfill \S\ref{canon-gener}
\item $n(x)$ \ \hfill \S\ref{Ucg}
\item $n(\Gm, q)$ \ \hfill \eqref{dimen}
\item $P_i$ \ \hfill \S\ref{canon-gener}
\item pr, pr$_2$ \ \hfill \S\ref{Ucg}
\item $\pi_i$ \ \hfill \S\ref{cgen}
\item $Q_n$ \ \hfill \eqref{Qndef}
\item $s$ \ \hfill \S\ref{sect-cdsg}
\item $S(y_0)$ \ \hfill \eqref{S(y)}
\item $t$ \ \hfill \S\ref{sect-cdsg}
\item $t(\Gm)$ \ \hfill \S\ref{canon-gener}
\item $\V_k(n, s), \V_k^0(n, s)$ \ \hfill Defn.~\ref{V_k}
\item $\W_r(\nu,s)$ \ \hfill \eqref{W_r}
\item $\WW$ \ \hfill \S\ref{Lie}
\item $\XX$ \ \hfill \S\ref{Lie}
\item $\zeta$ \ \hfill \S\ref{sect-cdsg}
\item $\omega$ \ \hfill \S\ref{Lie}
\item $\omega_r, \hat \omega_r$ \ \hfill \eqref{om-def-old}
\end{trivlist}
\end{multicols}

\raggedright
\newcommand\bibit[4]{\bibitem {#1}#2: {\em #3;\/ } #4}

\end{document}